\documentclass[11pt]{article}          
\usepackage{latexsym}
\usepackage{amssymb}
\usepackage{bbm}
\makeatletter
\def\@begintheorem#1#2{\trivlist%
 \item[\hskip \labelsep{\sffamily\bfseries #2\ #1}]\itshape}
\newtheorem{teo}{Theorem}[section]
\newtheorem{defi}[teo]{Definition}
\newtheorem{cor}[teo]{Corollary}
\newtheorem{lem}[teo]{Lemma}
\newtheorem{pro}[teo]{Proposition}
\newtheorem{_rem}[teo]{Remark}
\newtheorem{_eje}[teo]{Example}
\newtheorem{_conj}[teo]{Conjecture}
\newenvironment{rem}{\def\@begintheorem##1##2{\trivlist%
 \item[\hskip\labelsep{\sffamily\bfseries ##2\ ##1}]}\begin{_rem}}{\end{_rem}}

\makeatother

\newenvironment{beweis}{{\em Proof:}}{\hfill $\rule{2mm}{2mm}$
\vspace{3mm}

}

\DeclareMathAlphabet{\Ma}{U}{msa}{m}{n}
\DeclareMathAlphabet{\Mb}{U}{msb}{m}{n}
\DeclareMathAlphabet{\Meuf}{U}{euf}{m}{n}
\DeclareSymbolFont{ASMa}{U}{msa}{m}{n}
\DeclareSymbolFont{ASMb}{U}{msb}{m}{n}
\DeclareMathSymbol{\hrist}{\mathord}{ASMa}{"16}
\DeclareMathSymbol{\varkappa}{\mathalpha}{ASMb}{"7B}
\DeclareMathSymbol{\CrPr}{\mathord}{ASMb}{"6F}

\def\got#1{\Meuf{#1}}
\def\wz{\widetilde{Z}}
\def\tz{\widetilde{\al Z.}}
\def\tg{\widetilde{\Gamma}}

\def\mr #1.{\mathrm{#1\,}}
\def\mrt #1.{\mathrm{\mbox{\tiny #1\,}}}   
\def\mt #1.{{\mbox{\tiny $#1$}}}
\def\ms #1.{{\mbox{\small $#1$}}}
\def\ol #1.{\overline{#1}}
\def\mb #1.{\mathbf{#1\,}}
 
\def\1{\mathbbm 1} 
\def\wg{\wwh {\al G.}.}

\def\restriction{{\mathchoice{
 \mbox{\unitlength1cm\begin{picture}(.2,.4)%
  \bezier{5}(.07,.3)(.1,.27)(.13,.24)%
  \put(.07,.35){\line(0,-1){.5}}\end{picture}}}{
 \mbox{\unitlength1cm\begin{picture}(.2,.4)%
  \bezier{5}(.07,.3)(.1,.27)(.13,.24)%
  \put(.07,.35){\line(0,-1){.5}}\end{picture}}}{
  \hrist}{\hrist}}}

  \def\al #1.{{\mathcal{#1}}}
  \def\ot #1.{{\got{#1}}}
  \def\C{\Mb{C}}
  \def\c{\mt {\C}.}
  \def\kl #1.{{\mbox{\tiny {$\Mb{#1}$}}}}
  \def\c{\kl C.}
  \def\A{\Mb{A}}
  \def\D{\Mb{D}}

  \def\N{\Mb{N}}
  \def\Z{\Mb{Z}}
  \def\P{\Mb{P}}
  \def\Q{\Mb{Q}}
  
  \def\t #1.{\tilde{#1}} 
  \def\T #1.{\widetilde{#1}}

\def\bt{\mt {\boxtimes}.}

\def\Wort#1{\mbox{\fontfamily{cmr}\selectfont\mdseries\upshape #1}}
\def\In#1{{\fontsize{8pt}{10pt}\selectfont\Wort{#1}}}
\def\KIn#1{{\fontsize{6pt}{8pt}\selectfont\Wort{#1}}}

\DeclareMathSymbol{\hsemi}{\mathord}{ASMb}{"6E}
\newcommand{\semi}[2]{\mbox{$#1\kern.1em\hsemi\kern.1em#2$}}

\def\LA{\left\langle\bgroup}
\def\LE{\left[\bgroup}
\def\LG{\left\{\bgroup}
\def\LR{\left(\bgroup}
\def\RA{\egroup^{\rule{0mm}{2mm}}\right\rangle}
\def\RE{\egroup^{\rule{0mm}{2mm}}\right]}
\def\RG{\egroup^{\rule{0mm}{2mm}}\right\}}
\def\RR{\egroup^{\rule{0mm}{2mm}}\right)}
\def\Ldummy{\left.\bgroup}
\def\Rdummy{\egroup^{\rule{0mm}{2mm}}\right.}

\def\Kbegin{\begin{equation} \left. \begin{array}{rcl}}
\def\Kend{\end{array} \right\} \end{equation}}

\newcommand{\onb}{ortho\-normal basis}
\def\l2{\Lambda^{\mbox{\tiny $(2)$}}}

  \def\ccr #1,#2.{\overline{\Delta(#1,\,#2)}}

  \def\b #1.{{\bf #1}}
  \def\cross#1.{\mathrel{\mathop{\times}\limits_{#1}}}
  \def\C{\Mb{C}}
  \def\N{\Mb{N}}
  \def\P{\Mb{P}}
  
  \def\Z{\Mb{Z}}
  \def\T{\Mb{T}}

  \def\wh{\widehat}
  \def\dg{\widehat{\al G.}}
  \def\wwh #1.{\widehat{#1}}
  \def\wt #1.{\widetilde{#1}}

  \def\cross #1.{\mathrel{\raise 3pt\hbox{$\mathop\times\limits_{#1}$}}}
\def\set #1,#2.{\left\{\,#1\;\bigm|\;#2\,\right\}}
\def\b #1.{{\bf #1}}

\def\aut{{\rm Aut}\,}

\def\ob{{\rm Ob}\,}
\def\endo{{\rm End}\,}
\def\spec{{\rm spec}\;}

\def\Ad{{\rm Ad}\,}
\def\tr{{\rm Tr}\,}

\def\ol #1.{\overline{#1}}
\def\rn#1.{\romannumeral{#1}}
\def\chop{\hfill\break}
\def\rest{\restriction}
\def\un{\1}
\def\spa{{\rm span}}

\def\s #1.{_{\smash{\lower2pt\hbox{\mathsurround=0pt $\scriptstyle #1$}}\mathsurround=3pt}}
\def\bra #1,#2.{{\left\langle #1,\,#2\right\rangle_{\al A.}}}

\def\XP#1!{\renewcommand{\baselinestretch}{.7}\marginpar{{\footnotesize #1}\hfil}
\renewcommand{\baselinestretch}{1.5}}
\def\XB{\marginpar{
{\footnotesize\bf Change~starts----}\lower 11pt\hbox{\mathsurround=0pt$
\!\!\displaystyle{
\Bigg\downarrow}$\mathsurround=3pt}}}
\def\XE{\marginpar{{\footnotesize\bf Change~ends-----}\raise 10pt
\hbox{\mathsurround=0pt$ 
\!\!\displaystyle{
\Bigg\downarrow}$\mathsurround=3pt}}}
\def\HS{{\{\al F.,\,\al G.\}}}

\title{\bf Duality of compact groups and
           Hilbert C*-systems for C*-algebras with a nontrivial center}
           
\author{
 {\sc Hellmut Baumg\"artel}\\[2mm] 
 {\footnotesize Mathematical Institute, University of Potsdam,}\\
 {\footnotesize Am Neuen Palais 10, PF 601 553,}  \\       
 {\footnotesize D-14415 Potsdam, Germany.}                       \\[1mm]
 {\footnotesize baumg@rz.uni-potsdam.de}    \\
\and
 {\sc Fernando Lled\'o}  \\[2mm] 
 {\footnotesize Institute for Pure and Applied Mathematics,}       \\
 {\footnotesize RWTH-Aachen, Templergraben 55,}                   \\ 
 {\footnotesize D-52062 Aachen, Germany.}                         \\[1mm]
 {\footnotesize lledo@iram.rwth-aachen.de}}

\date{\today{}}
\begin{document}
\maketitle

\begin{center}
{\sl Dedicated to Detlev Buchholz on his 60th birthday.}
\end{center}
\vspace{.5cm}

\begin{abstract} 
In this paper we present a duality theory for compact 
groups in the case when the C*-algebra $\al A.$, the fixed 
point algebra of the corresponding Hilbert C*-system 
$(\al F.,\al G.)$, has a nontrivial
center $\al Z.\supset\C\1$ and the relative commutant 
satisfies the minimality condition
\[
 \al A.'\cap\al F.=\al Z.\,,
\]
as well as a technical condition called regularity.
The abstract characterization of the 
mentioned Hilbert C*-system is expressed by means 
of an inclusion of
C*-categories $\al T._\c < \al T.$, where $\al T._\c$ is a suitable
DR-category and $\al T.$ a full subcategory of the category of 
endomorphisms of $\al A.$. Both categories have the same objects
and the arrows of $\al T.$ can be generated from the arrows of 
$\al T._\c$ and the center $\al Z.$.

A crucial new element that appears in the present analysis is an
abelian group $\ot C.(\al G.)$, which we call the chain group
of $\al G.$, and that can be constructed from certain equivalence
relation defined on $\wh{\al G.}$, the dual object of $\al G.$.
The chain group, which is isomorphic 
to the character group of the center of $\al G.$,
determines the action of irreducible endomorphisms
of $\al A.$ when restricted to $\al Z.$.
Moreover, $\ot C.(\al G.)$ encodes the possibility of defining 
a symmetry $\epsilon$ also for the larger category $\al T.$ of the previous
inclusion.
\end{abstract}

\section{Introduction}

The superselection theory in algebraic quantum field theory, 
as stated by the
Doplicher-Haag-Roberts~(DHR) selection criterion
\cite{bHaag92,DHR69a,DHR69b},
led to a profound body of work, culminating in
the general Doplicher-Roberts~(DR) duality theory for compact groups
\cite{Doplicher89b}.
The DHR criterion selects a distinguished class of ``admissible''
representations of a quasilocal algebra $\al A.$ of observables,
which has trivial center
$\al Z.:=\al Z.(\al A.)=\C\1$.
This corresponds to the selection of a 
so-called DR-category $\al T.$,
which is a full subcategory of the category of 
endomorphisms of the C*-algebra $\al A.$ (see Definition~\ref{DRCat} below). 
Furthermore, from this endomorphism category 
$\al T.$
the
DR-analysis constructs a C*-algebra $\al F.\supset\al A.$
together with a compact group action 
$\alpha:\al G.\ni g\to\alpha_{g}\in\aut\al F.$
such that:
\begin{itemize}
\item
$\al A.$ is the fixed point algebra of this action
\item
$\al T.$ coincides with the category of all ``canonical
endomorphisms" of $\al A.$, associated with the pair
$\{\al F.,\alpha_{\al G.}\}$ (cf.~Subsection~\ref{SubCanEn}).
\end{itemize}
$\al F.$  
is called a Hilbert extension of $\al A.$ in \cite{bBaumgaertel92}.
Physically, $\al F.$ is identified as a field algebra and
$\al G.$ with a global gauge group of the system.
The pair $\{\al F.,\alpha_{\al G.}\}$, which we call 
{\em Hilbert C*-system} (cf.~Definition~\ref{defs2-1};
the name {\em crossed product} is also used),
is uniquely determined by
$\al T.$
up to
$\al A.$-module isomorphisms. Conversely,
$\{\al F.,\alpha_{\al G.}\}$
determines uniquely its category of all canonical endomorphisms.
Therefore $\{\al T.,\al A.\}$ can be seen as the abstract side
of the representation category of a compact group, while 
$\{\al F.,\alpha_{\al G.}\}$ corresponds
to the concrete side of the representation category of 
$\al G.$, and, roughly, any irreducible representations
of $\al G.$ is explicitly realized within 
the Hilbert C*-system.
One can state the equivalence of the ``selection
principle", given by
$\al T.$
and the ``symmetry principle", given by the compact group
$\al G.$. This is one of the crucial theorems of the
Doplicher-Roberts theory.

In the DR-theory the center $\al Z.$ of the C*-algebra $\al A.$ 
plays a peculiar role:
as stated above, if $\al A.$ corresponds
to the inductive limit of a net of local C*-algebras indexed by
open and bounded regions of Minkowski space, then the triviality
of the center of $\al A.$ is a consequence of standard assumptions
on the net of local C*-algebras. But, in general, the C*-algebra appearing 
in the DR-theorem does not need to be a quasilocal algebra and, in fact,
one has to assume explicitly that $\al Z.=\C\1$
in this context (see \cite[Theorem~6.1]{Doplicher89b}).
Finally, we quote from the introduction of the article
\cite{Doplicher89b}: ``There is, however, no
known analogue of Theorem~4.1 of \cite{Doplicher89a} for a C$^*$-algebra
with a non-trivial center and hence nothing resembling a ``duality''
in this more general setting.''

The aim of the present paper is to show that
{\em there is} a duality theory for compact groups in the
case of a nontrivial center, if the relative commutant of the 
corresponding Hilbert C*-system satisfies the following minimality
condition:
\begin{equation}\label{1afz}
\al A.'\cap\al F.=\al Z. 
\end{equation}
(cf.~Theorem~\ref{Teo2}). The essence of the 
previous result is that now the abstract characterization of the 
mentioned Hilbert C*-system is expressed by means 
of an inclusion of
C*-categories $\al T._\c < \al T.$, where $\al T._\c$ is a suitable
DR-category and $\al T.$ a full subcategory of the category of 
endomorphisms of $\al A.$. Both categories have the same objects
and the arrows of $\al T.$ can be generated from the arrows of 
$\al T._\c$ and the center $\al Z.$.

Several new elements appear in the generalization of the DR-theory 
studied here. The crucial one is an 
abelian group $\ot C.(\al G.)$, which we call the {\em chain group}
of $\al G.$, and that can be constructed from certain equivalence
relation defined on $\wh{\al G.}$, the dual object of the compact
group $\al G.$. The chain group, which is interesting
in itself and isomorphic to the character group of the center
of $\al G.$, determines the action of irreducible endomorphisms
of $\al A.$ when restricted to the center $\al Z.(\al A.)$.
Moreover, $\ot C.(\al G.)$ appears explicitly in the construction
of a family of examples realizing the inclusion of categories
$\al T._\c<\al T.$ mentioned above (cf.~Section~\ref{RealInclu}).
Finally, the chain group encodes also the possibility of defining 
a symmetry $\epsilon$ also for the larger category $\al T.$ of the previous
inclusion.

There are several reasons that motivate the generalization of
the DR-theory for systems satisfying the minimality condition 
(\ref{1afz}) for the relative commutant:
\begin{itemize}
\item[(i)] 
 In this context there is a nice intrinsic characterization
 of the Hilbert C*-systems satisfying (\ref{1afz}) and a further
 technical condition called regularity (cf.~Theorems~\ref{Teo1} 
 and \ref{Teo2}). One can also prove
 several results in the spirit of the DR-theory: for example,
 the category $\al T.$ is isomorphic to a
 subcategory $\al M._\al G.$ of the category of free Hilbert 
 $\al Z.$-bimodules generated by the algebraic Hilbert spaces in 
 $\al T._\al G.$ (cf.~Proposition~\ref{prop0}).
\item[(ii)] 
 In the context of compact groups, the equation (\ref{1afz}) 
 is also convenient for technical reasons. The minimality of
 the relative commutant implies that irreducible endomorphisms are
 mutually disjoint (cf.~Proposition~\ref{disj}) and this fact is
 crucial to have a nice decomposition of objects in terms of irreducible
 ones (cf.~Proposition~\ref{DecompEnd}).
\item[(iii)]
 The nontriviality of the center gives also the possibility to a more
 geometrical interpretation of the DR-theory. Indeed, from
 Gelfand's theorem we have $\al Z.\cong C(\Gamma)$, $\Gamma$ a compact
 Hausdorff space, and in certain situations the
 Hilbert C*-system $\{\al F.,\alpha_{\al G.}\}$ is a
 direct integral over $\Gamma$,
 where the Hilbert C*-system corresponding to a.e.~base point 
 $\lambda\in\Gamma$ is of a DR-type with the same group $\al G.$. 
 Here the chain group 
 plays again an important role.
 This more geometrical line of research has lead to recent developments
 in the context of vector bundles 
 (cf.~\cite{pVasselli03a,pVasselli03b,Vasselli03}).
\item[(iv)]
 There are physically relevant examples that satisfy the
 condition (\ref{1afz}).
 For example, this equation is presented in \cite{Mack90}  
 as a ``new principle". Moreover, the elements of the center $\al Z.$ 
 of $\al A.$ may be interpreted as classical observables contained
 in the quasilocal algebra.
\item[(v)] The present generalization of the DR-theory in the context
 minimal and regular Hilbert C*-systems has also found application 
 in the context of superselection theory for systems carrying quantum
 constraints (see \cite{pBaumgaertel03} as well as 
 \cite{Grundling85,Lledo00} for a C*-algebraic formulation of 
 the theory of quantum constraints).
\end{itemize}

The paper is structured in 9 sections: in Section~\ref{BasicHCS}
we introduce the notion of a Hilbert C*-system 
(cf.~Definition~\ref{defs2-1})
and give a detailed account of its properties. 
Hilbert C*-systems are special 
types of C*-dynamical systems $\{\al F.,\alpha_{\al G.}\}$ that,
in addition, contain the information of the representation category
of $\al G.$. They also satisfy important properties,
which are interesting in themselves, as for example: 
the fixed point algebra $\al A.$ is simple if $\al F.$ is simple 
(cf.~Subsection~\ref{SimpleAF} for further 
results on the ideal structure of these algebras); 
one can naturally introduce spectral
subspaces of $\al F.$ and prove Parseval-type equations for a suitable
$\al A.$-valued scalar product on $\al F.$ (cf.~Proposition~\ref{prop1}). 
Finally, Hilbert C*-systems provide a natural and 
concrete frame to describe the 
DR-theory as well as the generalization 
to the nontrivial center situation that we study here.
In Section~\ref{TwoEx} we study the 
important relation between two C*-categories
$\al T._\al G.$ and $\al T.$
that are naturally associated with a Hilbert C*-system.
In general, $\al T._\al G.$ is a subcategory of $\al T.$ 
and this inclusion turns out to be characteristic for the inverse 
result stated in Theorem~\ref{Teo2}. In 
Section~\ref{MinRegSect}
the main duality theorems are stated in the context of minimal
and regular Hilbert C*-systems. The next section defines 
the notion of an irreducible object and introduces the 
chain group of $\al G.$, denoted by $\ot C.(\al G.)$.  
We give examples of chain groups for several finite and
compact Lie groups and show that 
the chain group is isomorphic to the character group of the
center of $\al G.$ (see also \cite{pMueger03}). 
There is a close relation between
the chain group and the set of irreducible canonical
endomorphisms: an irreducible canonical endomorphism of $\al A.$ 
restricted to the center $\al Z.$ turns out to be an 
automorphism of $\al Z.$. We show that
there is a group homomorphism between the chain group and
the subgroup of $\mr aut.\al Z.$ generated by irreducible 
objects (cf.~Theorem~\ref{EndoChain}). One of the typical
difficulties in the context of a nontrivial center is that
$\al Z.$ is not stable under the action of a general canonical
endomorphism $\sigma$, i.e.
\[
 \sigma(\al Z.)\not\subset\al Z.\,.
\]
In this section we also give  
an explicit formula in terms of isotypical projections
that describes the action of reducible endomorphisms 
restricted to the center (cf.~Theorem~\ref{GeneralZMap}).
In Section~\ref{RealInclu} we construct a family of examples
that satisfy the requirements of the pair of
categories $\al T._\c<\al T.$ considered in Theorem~\ref{Teo2}.
In Section~\ref{TrivialChainHom} we analyze the 
situation where the homomorphism between the chain group and
the subgroup of $\mr aut.\al Z.$ generated by irreducible 
objects is trivial. In this case $\al Z.$ becomes the common 
center of $\al A.$ {\it and} $\al F.$. We can therefore 
decompose these algebras, which in this section are assumed
to be separable, w.r.t.~$\al Z.$. Then the 
Hilbert C*-system $\{\al F.,\alpha_{\al G.}\}$ becomes
a direct integral over $\Gamma:=\mr spec.\al Z.$ and 
the fibre Hilbert C*-system corresponding to the base point 
$\lambda\in\Gamma$ is of a DR-type with the same group $\al G.$. 
That means, in particular, that
the fixed point algebra associated with a.e.~$\lambda$ has a trivial 
center. Another simplifying condition of the present situation
is the fact that any canonical endomorphism acts trivially 
on the center, i.e.~$\rho\rest\al Z.=\mr id.\rest\al Z.$. Moreover, 
we show that in this case the minimality condition already implies
the regularity of the corresponding Hilbert C*-system
(cf.~Corollary~\ref{Mintoreg}). The special situation studied
in this section is also related to the notion of extention of C*-categories
by abelian C*-algebras (cf.~\cite{Vasselli03}).

Some conclusions connecting the present analysis to related
lines of research are stated in Section~\ref{Conclu}. Finally,
the paper contains 
an appendix recalling the decomposition
of a C*-algebra w.r.t.~its center.


\section{Basic properties of Hilbert C*-systems}\label{BasicHCS}

In this section  we summarize the structures from superselection
theory which we need. For proofs, we refer to the literature
if possible, otherwise proofs are included in this paper.

Below $\al F.$ will always denote a unital C*-algebra. 
A Hilbert space $\al H.\subset
\al F.$ is called {\it algebraic} if the scalar product $\langle\cdot,\cdot
\rangle$ of $\al H.$ is given by
$\langle A,B\rangle\un := A^{\ast}B$ for $A,\; B\in\al H.\,.$ Henceforth,
we consider only finite-dimensional algebraic Hilbert spaces. The 
support 
$\hbox{supp}\,\al H.$
of $\al H.$ is defined by
$\hbox{supp}\,\al H.:=\sum_{j=1}^{d}\Phi_j\Phi_{j}^{\ast}$,
where $\{\Phi_j\,\big|\,
j=1,\ldots,\,d\}$ is any orthonormal basis of $\al H..$ 
Unless otherwise specified, we assume below that
each considered algebraic Hilbert space $\al H.$
satisfies  ${\rm supp}\,\al H.
=\un.$

We also fix a compact 
C*-dynamical system
$\{\al F.,\al G.,\alpha\}$,
i.e.
$\al G.$
is a compact group and
$\alpha:\al G.\ni g\to\alpha_{g}\in\aut\al F.$
is a pointwise norm-continuous morphism.
For $D\in\wh{\al G.}$
(the dual of $\al G.$) its {\it spectral projection}
$\Pi_{\mt D.}\in\al L.(\al F.)$
is defined by
\begin{eqnarray} \label{PiD}
\Pi_{{\mt D.}} (F)&:=&\int_{\al G.}\ol\chi_{\mt D.} (g).\,\alpha_{g}(F)\,dg
\quad\hbox{for all}\quad F\in\al F.,  \\[1mm]
\nonumber \hbox{where:}\quad\qquad
\chi_{\mt D.} (g)&:=&\dim{D}\cdot\tr\pi(g),\quad\pi\in D\,,
\end{eqnarray}
and $dg$ is the normalized Haar measure of the compact group $\al G.$.
The spectrum of 
$\alpha_{\al G.}$ 
can then be defined by
\[
\spec\alpha_{\al G.}:=\set D\in\wh{\al G.}, \Pi_{\mt D.}\not=0.\,.
\]
Note that
$\spec\alpha_{\al G.}$
coincides with the so-called
Arveson spectrum of
$\alpha_{\al G.}$ 
(see e.g.~\cite{Baumgaertel95}).
 
Our central object of study is:

\begin{defi}\label{defs2-1}
The compact C*-dynamical system
$\{\al F.,\al G.,\alpha\}$
is called a {\bf Hilbert C*-system} if
for each $D\in\wh{\al G.}$
there is an algebraic Hilbert space 
$\al H._{\mt D.}\subset\Pi_{\mt D.}\al F.,$
such that
$\alpha_{\al G.}$
acts invariantly on
$\al H._{\mt D.},$  
and the unitary representation
$\alpha_\al G.\rest\al H._{\mt D.}$ 
is in the equivalence class 
$D\in\wh{\al G.}.$
\end{defi}

We are mainly interested in Hilbert C*-systems whose fixed point
algebras coincide such that they appear as extensions of it.

\begin{defi}\label{defs2-2}
A Hilbert C*-system 
$\{\al F.,\al G.,\alpha\}$ 
is called a {\bf Hilbert
extension} of a C*-algebra $\al A.\subset\al F.$ if $\al A.$ is the
fixed point algebra of ${\al G.}.$
Two Hilbert extensions 
$\{\al F._i,\,\al G.\,,\alpha^{i}\},\;i=1,\,2$
of $\al A.$
(w.r.t.~the same group $\al G.$)
are called 
$\al A.\hbox{\bf-module isomorphic}$ 
if there is an isomorphism 
$\tau:\al F._1\to\al F._2$ 
such that
$\tau(A)=A$ 
for 
$A\in\al A.,$ 
and 
$\tau$ 
intertwines the group actions, i.e.
$\tau\circ\alpha^{1}_g=\alpha^{2}_g\circ\tau$, $g\in\al G.$.
\end{defi}

\begin{rem}
\begin{itemize}
\item[(i)]
For a Hilbert C*-system
$\{\al F.,\al G.,\alpha\}$
one has
$\hbox{spec}\,\alpha_{\al G.}=\wh{\al G.}$
and the morphism
$\alpha:\al G.\to \hbox{Aut}\,\al F.$
is necessarily faithful. So, since
$\al G.$
is compact and
$\hbox{Aut}\,\al F.$
is Hausdorff w.r.t.~the topology of pointwise norm-convergence,
$\alpha$
is a homeomorphism of
$\al G.$
onto its image. Thus
$\al G.$
and
$\alpha_{\al G.}$
are isomorphic as topological groups.
\item[(ii)]
Group automorphisms of
$\al G.$
lead to $\al A.$-module isomorphic Hilbert extensions of
$\al A.$,
i.e.~if
$\{\al F.,\al G.,\alpha\}$
is a Hilbert extension of
$\al A.$
and
$\xi$
an automorphism of
$\al G.$,
then the Hilbert extensions
$\{\al F.,\al G.,\alpha\}$
and
$\{\al F.,\al G.,\alpha\circ\xi\}$
are $\al A.$-module isomorphic.

Therefore, the Hilbert C*-system
$\{\al F.,\al G.,\alpha\}$
depends, up to $\al A.$-module isomorphisms, only on
$\alpha_{\al G.}$,
which is isomorphic to
$\al G.$.
In other words, up to $\al A.$-module isomorphism we may identify
$\al G.$
and
$\alpha_{\al G.}\subset\aut\al F.$
neglecting the action
$\alpha$
which has no relevance from this point of view. Therefore in the
following, unless it is otherwise specified, we use the notation
$\{\al F.,\al G.\}$
for a Hilbert extension of
$\al A.$,
where
$\al G.\subset\aut\al F.$.
\item[(iii)]
 As mentioned above,
Hilbert C*-systems arise in DHR-superselection theory
(cf.~\cite{bBaumgaertel92,bBaumgaertel95}). Mathematically, 
there are constructions by means of
tensor products 
$\al B.$
of Cuntz algebras
$\al O._{\al H._{u}},\;\al B.=\otimes_{u\in\ob\,\al R.}\al O._{\al H._{u}},$
where
$\al R.$
is a category whose objects $u$ are finite-dimensional continuous 
unitary representations 
of a compact group
$\al G.$
on Hilbert spaces
$\al H._{u}$
with
$\dim\,\al H._{u}>1$
and whose arrows are the corresponding intertwining operators 
(cf.~\cite[Section~7]{Doplicher88}). In these examples the center 
$\al Z.$ of the fixed point algebra 
$\al A.$
is trivial. 

Further examples
in the context of the CAR-algebra 
with an abelian group
$\al G.=\T$
and nontrivial
center $\al Z.$ are given
in \cite{Baumgaertel01}. In Section~\ref{RealInclu} we construct
a family of examples of minimal and regular Hilbert C*-systems
for nonabelian groups and with nontrivial $\al Z.$.
\end{itemize}
\end{rem}
\begin{rem}
\label{remark1}
A Hilbert C*-system is a very highly structured object;-
below we list some important properties (for details, consult
~\cite{bBaumgaertel95,bBaumgaertel92}):
\begin{itemize}
\item[(i)]
Given two $\al G.$-invariant
algebraic Hilbert spaces 
$\al H.,\al K.\subset\al F.,$
then
$\spa(\al H.\cdot\al K.)$
is also a $\al G.$-invariant
algebraic Hilbert space 
which we will 
briefly denote by $\al H.\cdot\al K..$
It carries the tensor product of the representations of $\al G.$
carried by $\al H.$ and $\al K..$
\item[(ii)]
Let
$\al H.,\al K.$
as before but not necessarily of support $\1$:
There is a natural isometric embedding of
$\al L.(\al H.,\al K.)$
into
$\al F.$
given by
\[
\al L.(\al H.,\al K.)\ni T\to\al J.(T):=\sum_{j,k}
t_{jk}\Psi_{j}\Phi^{\ast}_{k},\quad t_{jk}\in \C,
\]
where
$\{\Phi_{k}\}_{k}$ resp.~$\{\Psi_{j}\}_{j}$
are orthonormal basis of
$\al H.$ resp.~$\al K.$ and where
\[
T(\Phi_{k})=\sum_{j}t_{jk}\Psi_{j},
\]
i.e.
$(t_{j,k})$
is the matrix of $T$ w.r.t.~these orthonormal basis. One has
\[
T(\Phi)=\al J.(T)\cdot\Phi,\quad \Phi\in\al H..
\]
For simplicity of notation we will often put $\wwh T.:=\al J.(T)$.
Moreover, we have
$\wwh T.\in\al A.$
iff
$T\in\al L._{\al G.}(\al H.,\al K.)$,
where
$\al L._{\al G.}(\al H.,\al K.)$
denotes the linear subspace of
$\al L.(\al H.,\al K.)$
consisting of all intertwining operators of the representations
of 
$\al G.$
on 
$\al H.$ and $\al K.$ (cf.~\cite[p.~222]{bBaumgaertel92}).

\item[(iii)]
Generally for a Hilbert C*-system, the assignment $D\to\al H._{\mt D.}$
is not unique.
If
$U\in \al A.$
is unitary then also
$U\al H._{\mt D.}\subset\Pi_{{\mt D.}}\al F.$
is an $\al G.$-invariant algebraic Hilbert space carrying 
a representation in
$D$. Note that each 
$\al G.$-invariant algebraic
Hilbert space 
$\al K.$
which carries a representation of
$D$
is of this form, i.e.~there is a unitary
$V\in\al A.$
such that
$\al K.=V\al H._{{\mt D.}}.$
\item[(iv)]
There is a useful partial order on the $\al G.$-invariant
algebraic Hilbert spaces. We define
$\al H.<\al K.$
to mean that there is an orthoprojection
$E$
on
$\al K.$
such that
$E\al K.$
is invariant w.r.t.
$\al G.$
and the representation
$\al G.\rest\al H.$
is unitarily equivalent to
$\al G.\rest E\al K.$.

Note that
$\al H.<\al K.$
iff there is an isometry
$V\in\al A.$
such that
$VV^{\ast}=:E$
is a projection of
$\al K.$,
i.e.
$V\al H.=E\al K.$
(use (ii)). 
\item[(v)]
Given a Hilbert C*-system 
$\HS$ a useful
*-subalgebra of
$\al F.$ is
\[
\al F._{\rm fin}:= \set F\in\al F.,\Pi_{\mt D.} F\not=0
\quad\hbox{for only finitely many $D\in\wh{\al G.}$}.
\]
which is dense in $\al F.$ w.r.t.~the C*-norm (cf.~\cite{Shiga55}).
\item[(vi)]
The spectral projections  
satisfy:
\begin{eqnarray*}
\Pi_{{\mt D.}_1}\Pi_{{\mt D.}_2} 
                    &=& \Pi_{{\mt D.}_2}\Pi_{{\mt D.}_1}=
                        \delta\s{\mt D.}_1{\mt D.}_2.\Pi_{{\mt D.}_1}\\[1mm]
\|\Pi_{\mt D.}\| &\leq& d({D})^{3/2}\;,
\qquad d({D}):=\dim(\al H._{\mt D.})\;,      \\[1mm]
\Pi_{\mt D.}\al F.  &=&  \spa(\al AH._{\mt D.})\;,    \\[1mm]
\Pi_{{\mt D.}}(AFB) &=&  A\cdot\Pi_{{\mt D.}}(F)\cdot B,\quad
A,B\in\al A.,\,F\in\al F.\;,   \\[1mm]
\al A. &=& \Pi_{\iota}\al F.\;,
\end{eqnarray*}
where
$\iota\in\wh{\al G.}$
denotes the trivial representation of $\al G..$
\item[(vii)]
In $\al F.$ there is an $\al A.$-scalar product given by
${\langle F,\, G\rangle_{\al A.}:=\Pi_\iota FG^*},$ w.r.t.~which
the spectral projections are symmetric, i.e.
$\bra\Pi_{\mt D.} F,G.=\bra F,\Pi_{\mt D.} G.$ 
for all $F,\; G\in\al F.,$ $D\in\wh{\al G.}$.
Using the $\al A.$-scalar product one can define a norm on
$\al F.$,
called the $\al A.$-norm
\[
\vert F\vert_{\al A.}:=\Vert\langle F,F\rangle\Vert^{1/2},\quad
F\in \al F..
\]
Note that
$\vert F\vert_{\al A.}\leq \Vert F\Vert$
and that
$\al F.$
in general is not closed w.r.t.~the $\al A.$-norm.
\end{itemize}
\end{rem}

The following result confirms the importance and naturalness of the 
previously defined norm $|\cdot|_\al A.$ in the context of Hilbert 
C*-systems. This norm plays also a fundamental role in the 
so-called inverse
superselection theory which reconstructs the Hilbert C*-system from
the data $\al A.$ and a suitable family of endomorphisms of $\al A.$
(cf.~\cite{Baumgaertel97,bBaumgaertel95,Lledo01a}).

\begin{pro}
\label{prop1}
Let 
$\HS$ 
be a Hilbert C*-system, then for
each $F\in\al F.$ we have
\begin{equation}\label{FG}
F=\sum_{{\mt D.}\in\wh{\al G.}}\Pi_{\mt D.} F
\end{equation}
where the sum on the right hand side is convergent w.r.t.~the 
$\al A.$-norm 
and we have Parseval's equation:
\begin{equation}\label{PEq}
\langle F,F\rangle_{\al A.}
=\sum_{{\mt D.}\in\wh{\al G.}}\langle\Pi_{\mt D.} F,\Pi_{\mt D.} 
  F\rangle_{\al A.}
\;.
\end{equation}
\end{pro} 
\begin{beweis}
Let
$\Gamma\subset\wh{\al G.},\;\hbox{card}\,\Gamma<\infty$.
The set
$\{\Gamma\}$
of all such subsets of
$\wh{\al G.}$
is a directed net. The assertion (\ref{FG}) means
\[
\sum_{{\mt D.}\in\wh{\al G.}}\Pi_{{\mt D.}}F:=
\lim_{\Gamma\to\wh{\al G.}}F_{\Gamma},
\]
where
\[
F_{\Gamma}:=\sum_{{\mt D.}\in\Gamma}\Pi_{{\mt D.}}F,
\]
and ``lim'' means convergence w.r.t.~the $\al A.$-norm. On the
other hand, if
$\Gamma$
is fixed, we put
\[
G_{\Gamma}=G_{\Gamma}(C_{{\mt D.}},\;D\in\Gamma):=
\sum_{{\mt D.}\in\Gamma}C_{{\mt D.}},\quad
C_{{\mt D.}}\in\Pi_{{\mt D.}}\al F..
\]
Then
$G_{\Gamma}\in\al F._{\rm fin}$.
By a simple calculation one obtains
\[
\langle F-G_{\Gamma},F-G_{\Gamma}\rangle_{\al A.}=
\langle F,\,F\rangle_{\al A.}-\sum_{{\mt D.}\in\Gamma}\langle\Pi_{{\mt D.}}F,
\Pi_{{\mt D.}}F\rangle_{\al A.} +
\sum_{{\mt D.}\in\Gamma}\langle \Pi_{{\mt D.}}F-C_{{\mt D.}},
\Pi_{{\mt D.}}F-C_{{\mt D.}}\rangle_{\al A.}.
\]
Since
\[
\sum_{{\mt D.}\in\Gamma}\langle\Pi_{{\mt D.}}F-C_{{\mt D.}},
\Pi_{{\mt D.}}F-C_{{\mt D.}}\rangle_{\al A.}\geq 0
\]
we obtain
\begin{eqnarray}
\langle F-G_{\Gamma},F-G_{\Gamma}\rangle_{\al A.} &\geq&
\langle F,\,F\rangle_{\al A.}-\sum_{{\mt D.}\in\Gamma}
\langle\Pi_{{\mt D.}}F,\,
\Pi_{{\mt D.}}F\rangle_{\al A.}\label{lasteq}  \\[1mm]
&=& \Big\langle F-\sum_{{\mt D.}\in\Gamma}\Pi_{{\mt D.}}F,\,
F-\sum_{{\mt D.}\in\Gamma}\Pi_{{\mt D.}}F\Big\rangle_{\al A.}\geq 0.
\nonumber
\end{eqnarray}

Therefore
\[
\langle F,\,
F\rangle_{\al A.}\geq \sum_{{\mt D.}\in\Gamma}
\langle\Pi_{{\mt D.}}F,\Pi_{{\mt D.}}F\rangle_{\al A.}
\]

Since
$\Vert X\Vert\geq\vert X\vert_{\al A.}$
for all
$X\in\al F.$
we have
\[
\Vert F-G_{\Gamma}\Vert\geq\vert F-G_{\Gamma}\vert_{\al A.}\geq
\vert F-F_{\Gamma}\vert_{\al A.}.
\]
According to Shiga's theorem (see ~\cite{Shiga55}) the left hand side
can be chosen arbitrary small for suitable
$\Gamma$
and suitable coefficients
$C_{{\mt D.}}$. Hence
$\vert F-F_{\Gamma}\vert_{\al A.}\to 0$
for
$\Gamma\to\wh{\al G.}$
follows. This is (\ref{FG}) and this implies
\[
\lim_{\Gamma\to\wh{\al G.}}\Vert\langle F,F\rangle_{\al A.}-
\sum_{{\mt D.}\in\Gamma}\langle\Pi_{{\mt D.}}F,\Pi_{{\mt D.}}F\rangle_{\al A.}
\Vert=0,
\]
which proves (\ref{PEq}).
\end{beweis}

Note that (\ref{FG}) does not 
in general converge
w.r.t.~the C*-norm $\Vert\cdot\Vert$.
\begin{cor}\label{Compo}
(i) Each $F\in\al F.$ is uniquely determined by its projections 
$\Pi_{\mt D.} F,$ $D\in\wh{\al G.},$ i.e.~$F=0$
iff $\Pi_{\mt D.} F =0$ for all $ D\in\wh{\al G.}.$
\chop
(ii) We have that $\big|\Pi_{\mt D.}\big|_{\al A.}=1$
for all $D\in\wh{\al G.}$,
where $\vert\cdot\vert_{\al A.}$
denotes the operator norm of $\Pi_{\mt D.}$ w.r.t.~the norm
$\vert\cdot\vert_{\al A.}$ in $\al F.$.
\end{cor}


\section{Two natural examples of C*-categories associated with a Hilbert
         C*-system}\label{TwoEx}

In the following we introduce two important examples of 
C*-categories that naturally appear in the context of Hilbert 
C*-systems. For the general definition and further 
properties of tensor C*-categories we refer to
\cite{Doplicher89b,bMaclane98}. We mention only that the 
notion of an {\bf irreducible object} introduced in
\cite[Section~5]{Lledo97b} (see also \cite{Lledo01a})
can be defined for arbitrary tensor C*-categories $\ot T.$:
$\rho\in\mr Ob.\ot T.$
is called {\bf irreducible} if 
\begin{equation}\label{AllgIrr}
 (\rho,\rho)=1_\rho\times (\iota,\iota) \,,
\end{equation}
where $\iota$ denotes the unit for the tensor product of objects,
$1_\rho$ is the unit of the unital C*-algebra $(\rho,\rho)$
and $\times$ is the tensor product of arrows.
We denote the set of all irreducible objects in $\ot T.$ by 
$\mb Irr.{\ot T.}$.

\subsection{The category $\al T._\al G.$ of all $\al G.$-invariant algebraic
            Hilbert spaces}\label{GCanEnd}

The $\al G.$-invariant algebraic Hilbert spaces
$\al H.$ 
of
$\{\al F.,\al G.\}$, satisfying $\mr supp.\al H.=\1$,
form the objects of a C*-category
$\al T._{\al G.}$
whose arrows are given by 
$(\al H.,\al K.):=\al J.(\al L._{\al G.}(\al H.,\al K.))\subset\al A.$.
The tensor product of objects is given by the product in $\al F.$,
the unit object is $\iota:=\C\1$ and $(\iota,\iota)=\C\1$. The composition
of arrows $\al J.(T)\in(\al H.,\al H.')$, $\al J.(S)\in(\al K.,\al K.')$
is given by $\al J.(T\otimes S)\in(\al H.\al K.,\al H.'\al K.')$, where
$T\otimes S\in\al L._{\al G.}(\al H.\otimes\al K.,\al H.'\otimes\al K.')$.
$\al H.$ is irreducible iff $(\al H.,\al H.)=\C\1$ (Schur's lemma).

We will focus next on the additional structure of $\al T._\al G.$. 
For this recall the partial order in
$\hbox{Ob}\,\al T._{\al G.}$
given in Remark~\ref{remark1}~(iv).
If
$\al K.\in\hbox{Ob}\,\al T._{\al G.}$
is given, an object
$\al H.<\al K.$
is called a {\it subobject} of
$\al K.$.
If
$E\in\al J.(\al L._{\al G.}(\al K.))$
is an orthoprojection
$0<E<\un,$
i.e.~$E$
is a reducing projection for the representation of
$\al G.$
on
$\al K.$,
then the question arises whether there is an object
$\al H.$
such that the representations on
$\al H.$
and
$E\al K.$
are unitarily equivalent. This suggests the concept of {\it
closedness} of
$\al T._{\al G.}$
w.r.t.~subobjects.
\begin{defi}
The category
$\al T._{\al G.}$
is {\bf closed w.r.t.~subobjects} if to each
$\al K.\in\hbox{Ob}\,\al T._{\al G.}$
and to each nontrivial orthoprojection
$E\in\al J.(\al L._{\al G.}(\al K.))$
there is an isometry
$V\in\al A.$
with
$VV^{\ast}=E.$
(In this case
$\al H.:=V^{\ast}\al K.$
is a subobject
$\al H.<\al K.$
assigned to
$E$.)
\end{defi}

Second, if
$V,W\in\al A.$
are isometries with
$VV^{\ast}+WW^{\ast}=\un$
and
$\al H.,\al K.\in\hbox{Ob}\,\al T._{\al G.}$
then 
we call the algebraic Hilbert space
$V\al H.+W\al K.$
of support $\1$ a {\it direct sum} of
$\al H.$ and $\al K.$.
It is $\al G.$-invariant and carries the direct sum of the
representations on
$\al H.$ and $\al K.$.
Therefore we define

\begin{defi}
The category
$\al T._{\al G.}$
is {\bf closed w.r.t.~direct sums} if to each
$\al H._{1},\al H._{2}\in\ob\,\al T._{\al G.}$
there is an object
$\al K.\in\ob\,\al T._{\al G.}$
and there are isometries
$V_{1},V_{2}\in\al A.$
with
$V_{1}V_{1}^{\ast}+V_{2}V_{2}^{\ast}=\un$
such that
$\al K.=V_{1}\al H._{1}+V_{2}\al H._{2}.$
\end{defi}

Since
$\ob\,\al T._{\al G.}$
contains {\it all} $\al G.$-invariant algebraic Hilbert spaces,
$\al T._{\al G.}$
is always closed w.r.t.~direct sums provided that
$\al A.$
contains a pair
$V,W$
of isometries with
$VV^{\ast}+WW^{\ast}=\un$.
This condition for $\al A.$ we call {\em Property~B}. 
It will play an important role in the rest of the paper.
We have

\begin{pro}
If $\al A.$ satisfies Property~B then
$\al T._{\al G.}$
is closed w.r.t.~direct sums.
\end{pro}

\begin{pro}
\begin{itemize} 
\item[(i)] 
If the category $\al T._{\al G.}$ is closed w.r.t.~direct sums, then
it is closed w.r.t.~subobjects.
\item[(ii)] Let $\al G.$ be nonabelian. 
If the category $\al T._{\al G.}$ is closed w.r.t.~subobjects, then
it is closed w.r.t.~direct sums.
\end{itemize}
\end{pro}
\begin{beweis}
(i) First we assume that
$\al T._{\al G.}$
is closed w.r.t.~direct sums. Note that in this case for any
$n\in\N$
there are isometries
$W_{j}\in\al A.,j=1,2,...,n$,
such that
$\sum_{j=1}^{n}W_{j}W_{j}^{\ast}=\un.$
Then to each finite-dimensional representation
$U$
there is a $\al G.$-invariant algebraic Hilbert space
$\al K.$
such that
$U$
is realized on
$\al K.$,
because
$U=\oplus_{{\mt D.}\in\wh{\al G.}}\;m_{\mt D.}\cdot U_{\mt D.}$,
$U_{\mt D.}\in D$,
and in
$\al A.$ there are isometries
$W_{{\mt D.},\,l}, l=1,2,...,m_{{\mt D.}}$ ($m_{{\mt D.}}$ being the 
multiplicity of $U_\mt D.$ in the decomposition of $U$)
such that
$\sum_{{\mt D.},\,l}W_{{\mt D.},\,l}W_{{\mt D.},\,l}^{\ast}=\un.$
Therefore
$\al K.:=\sum_{{\mt D.},\,l}W_{{\mt D.},\,l}\al H._{{\mt D.}}$
is an object from
$\al T._{\al G.}$
and carries exactly the representation
$U$.

Now let
$\al K.\in\ob\,\al T._{\al G.}$
and
$E\in\al J.(\al L._{\al G.}(\al K.)),\,0<E<\un,$
a reducing projection, 
i.e.~$E\al K.\subset\al K.$
is a reducing subspace that carries a certain representation of
$\al G.$. 
Note that $\mr supp. E\al K.\not=\1$. 
Nevertheless there is an object
$\al H.\in\ob\,\al T._{\al G.}$
which carries this representation: choose in
$E\al K.$
and
$\al H.$
the (orthonormal) basis
$\{\Phi_{j}\}_{j}$
of
$E\al K.,\;\{\Psi_{j}\}_{j}$
of
$\al H.$
in such a way that the representation matrices coincide. Put
$A:=\sum_{j}\Phi_{j}\Psi_{j}^{\ast}$.
Then
$A\in\al A.,\;A^{\ast}A=\un$
and
$AA^{\ast}=E$
follows, i.e.
$\al H.$
is a subobject of
$\al K.$
w.r.t.
$E$.

(ii) Now we assume that
$\al T._{\al G.}$
is closed w.r.t.~subobjects. Then choose
$\al H._{\mt D._{1}},\al H._{\mt D._{2}}\in\ob\,\al T._{\al G.}$,
whose dimensions are larger than 1.
Then
$\al K.:=\al H._{\mt D._{1}}\cdot\al H._{\mt D._{2}}\in\ob\,\al T._{\al G.}$
and it carries the reducible representation
$U_{\mt D._{1}}\otimes U_{\mt D._{2}}$,
i.e.~there is a projection
$E,\,0<E<\un,\,E\in\al J.(\al L._{\al G.}(\al K.)).$
Then to
$E$
and
$\un-E$
there correspond isometries
$V,W\in\al A.$
with
$VV^{\ast}+WW^{\ast}=\un$,
hence
$\al T._{\al G.}$
is closed w.r.t.~direct sums.
\end{beweis}

\begin{rem}\label{ChangeAb}
Note that if the group
$\al G.$
is a compact abelian then
$\wh{\al G.}$
is a discrete abelian {\it group}, the character group. 
Pontryagin's duality theorem shows that in this case the 
notions of direct sums and subobjects are irrelevant for
the duality theory (see also Remark~\ref{commut}).
If the compact group is non abelian the duality theory
changes radically and closure under direct sums and 
subobjects become essential properties.
\end{rem}

\subsection{The category $\al T.$ of all canonical endomorphisms}
\label{SubCanEn}

\begin{defi}
To each $\al G.$-invariant algebraic Hilbert space $\al H.\subset
\al F.$  
there is assigned a corresponding 
inner endomorphism $\rho\s{\al H.}.\in\endo\al F.$ given by
\[
\rho\s{\al H.}.(F):=\sum_{j=1}^{d(\al H.)}\Phi_jF\Phi_j^*\;,
\]
where $\{\Phi_j\,\big|\,j=1,\ldots,\,d(\al H.)\}$ is any orthonormal basis of
$\al H.$. We call {\bf canonical endomorphism} the restriction of
$\rho_{\al H.}$ to
$\al A.$,
i.e.
$\rho_{\al H.}\rest\al A.\in\mbox{end}\,\al A.$.
\end{defi}  

\begin{rem}
\label{remark2}
\begin{itemize}
\item[(i)]
Note that the definition of the canonical endomorphisms
uses terms of
$\al F.$
explicitly. Therefore, the question arises whether the 
inner endomorphisms $\rho_\al H.$ can be
characterized by intrinsic properties of their restriction to
$\al A.$ 
(see the beginning of Section~\ref{MinRegSect} below).
This interplay between the inner and the canonical
endomorphisms $\rho_{\al H.}$
resp.~$\rho_{\al H.}\rest\al A.$
plays an essential role in the DR-theory. 
Below, we omit
the restriction symbol and regard the
$\rho_{\al H.}$
also as endomorphisms of
$\al A.$.
We will identify the set of canonical endomorphisms
of $\al A.$ as the objects of a  very important category with interesting
closure properties.
\item[(ii)]
If the emphasis is only on the class
$D\in\wh{\al G.}$
and not on its corresponding algebraic Hilbert space
$\al H._{\mt D.},$
 we will write
$\rho_\mt D.$
instead of
$\rho_{\al H._\mt D.}$.
\item[(iii)]
Note that $\Phi A=\rho\s{\al H.}.(A)\Phi$ for all 
$\Phi\in\al H.$ and $A\in\al A..$
\item[(iv)]
Note that the identity endomorphism $\iota$ is assigned to
$\al H.=\C\un,$
i.e.~$\rho\s{\C\un}.:=\iota.$
\item[(v)]
Let
$\al H.,\al K.$
be as before, then 
$\rho\s{\al H.}.\circ\rho\s{\al K.}.=\rho\s{\al H.\cdot\al K.}..$
\item[(vi)]
Whilst an invariant algebraic Hilbert space
uniquely determines its canonical endomorphism, in general
the converse does not hold.
\end{itemize}
\end{rem}
\begin{pro}
\label{HSnonunique}
Let
$\al H.,\al K.$
be $\al G.$-invariant algebraic Hilbert spaces.
Then: 
$\rho_{\al H.}\rest\al A.=\rho_{\al K.}\rest\al A.$
iff
$\Psi^{\ast}\Phi\in\al A.'\cap\al F.$
for all
$\Phi\in\al H.,\Psi\in\al K.$.
\end{pro}
\begin{beweis}
It is straightforward to check the condition for orthonormal
basis of 
$\al H.$
and
$\al K.$.
\end{beweis}

\begin{defi}
Let $\HS$ be a Hilbert C*-system with fixed point algebra $\al A.$.
The {\bf intertwiner space} of canonical endomorphisms 
$\sigma,\;\tau$ is:
\[
(\sigma,\,\tau):=\set X\in\al A.,X\sigma(A)=\tau(A)X\quad\hbox{for all}
\;{A\in\al A.}.
\]
and this is a complex Banach space.
We will say that $\sigma,\;\tau\in\endo\al A.$ are
{\bf mutually disjoint} if  ${(\sigma,\,\tau)}=\{0\}$ when
$\sigma\not=\tau$.

We denote by $\al T.$ the category  
with objects consisting of
the canonical endomorphisms
$\rho\s{\al H.}.$ for 
$\al G.$-invariant algebraic Hilbert spaces $\al H.\subset\al F.$
with $\mr supp.\al H.=\1$
and with arrows given by the intertwiner spaces.
\end{defi}

\begin{rem}
\label{remark3}
\begin{itemize}
\item[(i)]
$\al T.$ is the second example of a tensor C*-category. The 
tensor product of objects
is given by composition of endomorphisms (see Remark~\ref{remark2}(v))
and $\iota=\mr id.$. The composition of arrows is defined as follows:
For $A\in(\sigma,\sigma'),$ $B\in(\tau,\tau'),$ we put
${A\times B}:=A\sigma(B)\in{(\sigma\tau,\,\sigma'\tau')}$.

\item[(ii)]
We have
$(\iota,\,\iota)=\al Z. :=$ center of $\al A.$ and from the
Definition in (\ref{AllgIrr}) we have $\rho_\al H.\in\mr Ob.\al T.$ is 
{\bf irreducible} if $(\rho_\al H.,\rho_\al H.)=\rho_\al H.(\al Z.)$. 
Note that this corresponds precisely to the case where $\al G.$ acts  
irreducibly on $\al H.$ (see \cite[Subsection~3.1]{Lledo01a} and
\cite[Section~5]{Lledo97b}). We denote the set of irreducible 
objects in $\al T.$ by {\bf Irr}$\,\al T.$.

\item[(iii)] Recall the isometry 
$\al J.:\al L.\s{\al G.}.(\al H.,\al K.)
\longrightarrow \al A.$ encountered in Remark~\ref{remark1}(ii). 
We claim that
its image is in fact contained in ${\left(\rho\s{\al H.}.,\,
\rho\s{\al K.}.\right)}.$ To see this, let $\Phi\in\al H.,\;
A\in\al A.$ 
and
$T\in\al L._{\al G.}(\al H.,\al K.)$.
Then putting $\wwh T.:=\al J.(T)$ we have
\[
\wwh T.\rho\s{\al H.}.(A)\Phi=\wwh T.\Phi\cdot A=
T(\Phi)\cdot A=\rho\s{\al K.}.(A)T(\Phi)=
\rho\s{\al K.}.(A)\wwh T.\cdot\Phi
\]
hence 
\[
\wwh T.\rho\s{\al H.}.(A)=\rho\s{\al K.}.(A)\wwh T.
\]
i.e.~$\wwh T.\in {\left(\rho\s{\al H.}.,\,
\rho\s{\al K.}.\right)}$
or
\[
(\al H.,\al K.)=\al J.(\al L._{\al G.}(\al H.,\al K.))\subseteq
(\rho_{\al H.},\,\rho_{\al K.})\,.
\]
In general, the inclusion is proper. Note finally, that if 
$A=\al J.(T)$, $B=\al J.(S)$ for $T\in\al L._\al G.(\al H.,\al H.')$,
$S\in\al L._\al G.(\al K.,\al K.')$, then $A\times B=\al J.(T\otimes S)$,
i.e.~$\times$ restricted to the intertwiner spaces 
$(\al H.,\al K.)$ etc., coincides with the composition of arrows
of $\al T._\al G.$.

\item[(iv)]
Recall that 
 $\al H.<\al K.$
 iff there is an isometry 
 $V\in\al A.$
 such that
 $VV^{\ast}=:E$
 is a projection of
 $\al K.$
 i.e.
 $V\al H.=E\al K.$. 
 In this case we have
 $V\in (\rho_{\al H.},\rho_{\al K.})$
 and
 $E\in (\rho_{\al K.},\rho_{\al K.}).$
 Moreover, $E$ does not belong to the center of $\al A.$.
\end{itemize}
\end{rem}
 
There is an important connection between the categories
$\al T._{\al G.}$
and
$\al T. $ in the case of a trivial relative commutant
$\al A.'\cap\al F.=\C\un$
(see \cite[Lemma~2.4]{Doplicher88}).
Note that in this case $\al A.$ must have a trivial center
$\al Z.=\C\1$.

\begin{pro}\label{Isom}
There is a faithful functor from the categories
$\al T._{\al G.}$
to
$\al T. $
which is a bijection of objects.
In general the functor is not full, but
if the relative commutant satisfies
$\al A.'\cap\al F.=\C\un$,
then
$\al T._{\al G.}$
and
$\al T. $
are isomorphic categories.
\end{pro}
\begin{beweis}
For the objects the functor is specified by 
$\al H.\to \rho_{\al H.}$ and 
for the arrows by
$(\al H.,\al K.)\ni A\to A\in(\rho_{\al H.},\rho_{\al K.})$.
Note that the compatibility of functor w.r.t.~the composition of
arrows follows from Remark~\ref{remark3}~(iii).
For the second assertion use Proposition~\ref{HSnonunique} and that
$\Phi^{\ast}A\Psi\in\C\un$
for
$\Phi\in\al K.,\Psi\in\al H.$.
\end{beweis}

We now want to exhibit closure properties of
$\al T. $
similarly as for
$\al T._{\al G.}.$

\begin{defi}
\begin{itemize}
\item[(i)]
$\tau\in\ob\al T. $ is a {\bf subobject} 
 of $\sigma\in\ob\al T. ,$
denoted ${\tau<\sigma,}$ if
there there is an isometry $V\in(\tau,\sigma)$.
In this case $\tau(\cdot)=V^{\ast}\sigma(\cdot)V$
and $VV^{\ast}=:E\in (\sigma,\sigma)$ follow.
\item[(ii)]
$\rho\in\ob\,\al T. $
is a direct sum of
$\sigma,\tau\in\ob\,\al T. $,
if there are isometries
$V\in(\sigma,\rho),\,W\in(\tau,\rho)$
with
$VV^{\ast}+WW^{\ast}=\un$
such that
\[
\rho(\cdot)=V\sigma(\cdot)V^{\ast}+W\tau(\cdot)W^{\ast}=:\sigma
\oplus\tau.
\]
\end{itemize}
\end{defi}

\begin{rem}
\begin{itemize}
\item[(i)] 
The subobject relation ${\tau<\sigma}$ 
is again a partial order: let $\tau<\sigma$ and $\sigma<\mu$,
so that there are isometries $V\in(\tau,\sigma)$ and $W\in(\sigma,\mu)$.
Then $WV\in(\tau,\mu)$ is also an isometry, i.e.~$\tau<\mu$.

\item[(ii)] 
Note that if $\tau=\rho_{\al H.},\; \sigma=\rho_{\al K.}$
for $\al G.$-invariant algebraic Hilbert spaces 
$\al H.\,,\al K.$ satisfying $\al H.<\al K.$, then a 
fortiori $\tau<\sigma$. However if $\tau,\sigma$ are given and one
only knows that there are algebraic Hilbert spaces $\al H.<\al K.$,
then the transitivity property may not hold in general.

\item[(iii)]
If
$\sigma:=\rho_{\al H.},\,\tau:=\rho_{\al K.}$
then
$\rho=\rho_{\al L.}$
where
$\al L.:=V\al H.+W\al K.$.

\item[(iv)]
A direct sum
$\sigma\oplus\tau,$
defined above (where a priori $\rho$ is not necessarily an
object of
$\al T. $)
with isometries
$V,W\in\al A.\,,VV^{\ast}+WW^{\ast}=\un$
is only unique up to unitary equivalence, i.e.~if
$\rho,\rho'$
are direct sums of
$\sigma$
and
$\rho$,
then there is a unitary
$U\in(\rho,\rho')$.

\end{itemize}
\end{rem}

The closedness of
$\al T. $
w.r.t.~direct sums is defined by the closedness of
$\al T._{\al G.}$
w.r.t.~direct sums.
The closedness of $\al T.$ w.r.t.~subobjects 
is defined by the closedness of $\al T._{\al G.}$ w.r.t.~subobjects 
in the following sense:
if 
\begin{equation}\label{rhoH}
 \rho=\rho_{\al H.}\in\hbox{Ob}\,\al T. 
\end{equation}
is given, then for all $\al H.$ satisfying (\ref{rhoH}) and 
to each nontrivial projection
$E\in\al J.(\al L._\al G.(\al H.))$
there is an isometry
$V_\al H.\in\al A.$
with $V_\al H.V_\al H.^*=E$.

This means

\begin{pro}
If
$\al A.$
satisfies Property~B then
$\al T. $
is closed w.r.t.~direct sums and subobjects.
\end{pro}

\subsection{Permutation and conjugation structures on $\al T._\al G.$:
            DR-categories} \label{Conjperm}

To complete the analysis of the categories $\al T._\al G.$ and 
$\al T.$ we will recall briefly their permutation and conjugation
structure. First, we will consider these structures on 
$\al T._\al G.$ (cf.~Remark~\ref{ConstrER}~(ii) below). 
We assume in this subsection 
that the fixed point algebra $\al A.$ of $\HS$ satisfies Property~B.

\begin{pro} \label{PermuStructu}
{\bf (Permutation structure)} $\al T._\al G.$ has
a permutation structure, i.e.~a map 
\[
 \ob\al T._\al G.
      \ni\al H.,\al K.\to \epsilon(\al H.,\al K.)
      \in (\al H.\al K.,\,\al K.\al H.)\,,
\]
where $\epsilon(\al H.,\al K.)$ is unitary and satisfies
\begin{itemize}
\item[(i)] $\epsilon(\al H.,\al K.)\cdot\epsilon(\al K.,\al H.)=\1$.
\item[(ii)] $\epsilon(\C\1,\al H.)=\epsilon(\al H.,\C\1)=\1$.
\item[(iii)] $\epsilon(\al H._1\al H._2,\al H._3)
               =\epsilon(\al H._1,\al H._3)
                \cdot\rho_{\al H._1}
                     \big(\epsilon(\al H._2,\al H._3)\big)\quad$,
             $\al H._i\in\ob\al T._\al G.$, $i=1,2,3$.
\item[(iv)]  $\epsilon(\al H.',\al K.')\,A\times B
               = B\times A\,\epsilon(\al H.,\al K.)$,
                 $A\in(\al H.,\al H.')$,
                 $B\in(\al K.,\al K.')$.
\end{itemize}
\end{pro}

Tensor categories that have a map $\epsilon(\cdot,\cdot)$ satisfying
the properties (i)-(iv) adapted from above are called 
{\it symmetric} (cf.~\cite[p.~160]{Doplicher89b}). 
The map $\epsilon(\cdot,\cdot)$ is also
called a {\it permutator} or {\it symmetry}.

\begin{pro} \label{ConjuStructu}
{\bf (Conjugation structure)}
$\al T._\al G.$ has
a conjugation structure, i.e.~to each $\al H.\in\ob\al T._\al G.$
there is a conjugated algebraic Hilbert space 
$\ol {\al H.}.\in\ob\al T._\al G.$ carrying the corresponding conjugated
representation and there are {\em conjugates}
$R_\al H.\in (\C\1,\ol {\al H.}.\al H.)$,
$S_\al H.=\epsilon(\ol {\al H.}.,\al H.)\,R_\al H.$ such that
\[
 S_\al H.^*\, \rho_\al H.(R_\al H.)=\1\quad\mr and.\quad
 R_\al H.^*\, \rho_{\ol {\al H.}.}(S_\al H.)=\1\,.
\]
\end{pro}

\begin{rem}\label{ConstrER}
\begin{itemize}
\item[(i)]
A permutator (or symmetry) $\epsilon(\cdot,\cdot)$ as in
Proposition~\ref{PermuStructu} is given by
\[ 
 \epsilon(\al H.,\al K.):=\al J.(\theta(\al H.,\al K.))\,,
\]
where $\theta$ denotes the flip operator of the tensor product
$\al H.\otimes\al K.$.
Let $\{\Phi_i\}_i$,$\{\Psi_k\}_k$ 
be orthonormal basis of the algebraic Hilbert spaces 
$\al H.$ and $\al K.$, respectively. Then
\[
 \epsilon(\al H.,\al K.)=\sum_{i,k} \Psi_k\,\Phi_i\,\Psi_k^*\,\Phi_i^* \,.
\]
If $\al H.$ carries the direct sum of irreducible representations 
in $\{D_j\}_j \subset\wg$, then $\ol{\al H.}.$ carries the corresponding
direct sum of the conjugated representations in $\{\ol {D_j}.\}_j\subset\wg$.
Denote by $\{\ol \Phi._i\}_i$ the conjugated basis of $\ol {\al H.}.$
w.r.t.~$\{\Phi_i\}_i$, then we have the relation
\[
 R_\al H.=\sum_i \ol \Phi._i\,\Phi_i\,.
\]
\item[(ii)] 
It is possible to use the functor in Proposition~\ref{Isom} 
to transfer the corresponding permutation and conjugation 
structure to $\al T.$. Note, nevertheless, that if the 
inclusion $(\al H.,\al K.)\subseteq (\rho_{\al H.},\,\rho_{\al K.})$
is proper (cf.~Remark~\ref{remark3}~(iii)),
then the property corresponding to 
(iv) in Proposition~\ref{PermuStructu} is valid only for 
a smaller set of arrows.
\end{itemize}
\end{rem}

We can now sum up the rich structure of the category $\al T._\al G.$ 
in the notion of a (Doplicher/Roberts) DR-category 
(cf.~\cite{Doplicher89b}).

\begin{defi}\label{DRCat}
An (abstract) tensor C*-category 
$\al T._\c$ with $(\iota,\iota)=\C\1$,
closed w.r.t.~direct sums and subobjects, equipped with a
permutation and a conjugation structure is called an (abstract)
{\bf DR-category}. 
\end{defi}

The category $\al T._\al G.$ introduced in Subsection~\ref{GCanEnd} 
is an example of a DR-category.
It is a special {\em Tannaka-Krein} category 
for the group $\al G.$, where the objects
and the arrows are embedded in the algebra $\al F.$. 
Moreover, if $\al A.'\cap\al F.=\C\1$ (which implies $\al Z.=\C\1$), the 
category $\al T.$ of canonical endomorphisms is another
example of a DR-category (cf.~Proposition~\ref{Isom}). 

\begin{rem}\label{AssocD}
   In the context of the DR-theory we can associate with any
   $\rho\in\mr Irr.\al T._\c$ a unique element $D\in\widehat{\al G.}$, 
   where $\al G.$ is the group associated with the DR-category $\al T._\c$. 
   We denote by $\mr Irr._0\al T._\c$ a complete system of irreducible
   and mutually disjoint objects.
\end{rem}

One of the most fundamental results associated with DR-categories
is the existence of an integer-valued {\bf dimension function} on the 
objects of $\al T._\c$. It is defined as follows:
Let $\rho\in\mr Ob.\al T._\c$ and $R_\rho\in(\iota,\ol\rho.\rho)$
a conjugate. Then
\begin{equation}\label{DimFkt}
 d(\rho)\1:=R_\rho^* R_\rho\in(\iota,\iota)=\C\1\,.
\end{equation}
The dimension function $d(\cdot)$ is independent of the choice of
conjugates and gives the same value on unitarily equivalent objects.
Moreover, it satisfies the following properties 
(cf.~\cite[Sections~2]{Doplicher89b} or 
\cite[Subsection~11.1.6]{bBaumgaertel92}):

\begin{pro}\label{Dec}
Let $\mr Ob.\al T._\c\ni\rho\to d(\rho)$ be the dimension function 
defined above. Then for $\rho,\rho_1,\rho_2\in\mr Ob.\al T._\c$
we have
\begin{itemize}
\item[(i)] $d(\rho)\in\N$.
\item[(ii)] $d(\iota)=1\quad$ and $\quad d(\ol \rho.)=d(\rho)$.
\item[(iii)] $d(\rho_1\circ\rho_2)=d(\rho_1)\,d(\rho_2)\quad$ and 
             $\quad d(\rho_1\oplus\rho_2)=d(\rho_1)+ d(\rho_2)$
\item[(iv)]$\lambda= \bigoplus_{j=1}^{r}\bigoplus_{l=1}^{m(\rho_{j},\lambda)}
            \rho_{jl}$, with $\rho_{jl}:=
            \rho_{j}\in\mathrm{Irr}_0\,\mathcal{T}_\c$
            (recall Remark~\ref{AssocD}),$
            \,l=1,2,\ldots ,m(\rho_{j},\lambda)$ \\[2mm] and  
            $d(\lambda)=\sum_{j=1}^{r}m(\rho_{j},\lambda)d(\rho_{j})$, where
            $(\rho,\lambda)_\c$ are algebraic Hilbert spaces and 
            $m(\rho,\lambda):=\dim\,(\rho,\lambda)_\c$.
\end{itemize}
\end{pro}

\subsection{The ideal structure of Hilbert C*-Systems}
\label{SimpleAF}

Given a Hilbert C*-system $\HS$,
we will analyze in the present subsection
the relation between the ideal structures of $\al F.$ and 
of the fixed point algebra $\al A.$.
It is clear that these must be closely related since $\al F.$ can be
generated from $\al A.$ and 
$\{\al H._\mt D.\}_{\mt {D\!\in\!\wh {G}}.}$,
the latter being a complete system of irreducible algebraic
Hilbert spaces with support $\1$ 
(cf.~Definition~\ref{defs2-1} and Remark~\ref{remark1}~(v),(vi)).

First we introduce the following weaker notions of simplicity which
are natural in the context of Hilbert C*-systems.

\begin{defi}  \label{GSimple}
$\HS$ denotes a Hilbert C*-system.
  \begin{itemize}
  \item[(i)] Let $\al E.\subset\al F.$ be a closed two-sided ideal in 
    $\al F.$, i.e.~$\al E.\lhd\al F.$. We say that $\al E.$ is 
    {\bf $\al G.$-invariant} if $g(\al E.)\subseteq\al E.$, $g\in\al G.$.
    $\al F.$ is {\bf $\al G.$-simple} if it has no nontrivial
    $\al G.$-invariant closed two-sided ideals.
   \item[(ii)] Let $\al I.\subset\al A.$ be a closed two-sided ideal in 
    $\al A.$, i.e.~$\al I.\lhd\al A.$. We say that $\al I.$ is 
    {\bf $\rho$-invariant} if $\rho_\mt D.(\al I.)\subseteq\al I.$, 
    $D\in\wh{\al G.}$, where $\{\rho_\mt D.\}_{\mt {D\!\in\!\wh {G}}.}$,
    are the canonical endomorphisms associated to a complete system
    $\{\al H._\mt D.\}_{\mt {D\!\in\!\wh {G}}.}$ of irreducible
    algebraic Hilbert spaces with support $\1$.
    $\al I.$ is {\bf $\rho$-simple} if it has no nontrivial
    $\rho$-invariant closed two-sided ideals.
  \end{itemize}
\end{defi}

\begin{rem}
Note that the notion of $\rho$-simplicity is independent of the 
particular choice of the system 
$\{\al H._\mt D.\}_{\mt {D\!\in\!\wh {G}}.}$. In fact, 
if $\al I.\lhd\al A.$ is $\rho$-invariant 
w.r.t.~$\{\rho_\mt D.\}_{\mt {D\!\in\!\wh {G}}.}$, 
then any other unitary equivalent endomorphism 
$\rho_\mt D.'(\cdot)=V\rho_\mt D.(\cdot)V^*$, 
$D\in\dg$ and $V$ a unitary in $\al A.$,
still satisfies
\[
 \rho_\mt D.'(\al I.)=V\rho_\mt D.(\al I.)V^*\subseteq
                      V\,\al I.\,V^*\subseteq\al I.\;.
\]
\end{rem}

\begin{lem}\label{IdealGInv}
Let $\al E.\lhd\al F.$ be $\al G.$-invariant and $\Pi_\mt D.$ the 
spectral projections defined in (\ref{PiD}). Then 
$\al E.$ is $\Pi_\mt D.$-invariant, 
i.e.~$\Pi_\mt D.(\al E.)\subseteq\al E.$ for all $D\in\wh{\al G.}$.
\end{lem}
\begin{beweis}
Let $E\in\al E.$. By the definition of spectral projection in
(\ref{PiD}) we have 
\[
 \Pi_{{\mt D.}}(E):=\int_{\al G.}\ol\chi_{\mt D.}(g).\,\alpha_{g}(E)\,dg\,.
\]
Since $\al E.$ is closed, we obtain from the definition of the integral and 
the pointwise norm continuity of the group action that
$\Pi_{{\mt D.}}(E)\in\al E.$.
\end{beweis}

\begin{pro}
Let $\HS$ be a Hilbert C*-system with fixed point algebra $\al A.$. Then
we have the following four implications:\\ \phantom{mmmm}
$\al F.$ is simple $\Rightarrow$ $\al F.$ is $\al G.$-simple
                   $\Rightarrow$ $\al A.$ is simple
                   $\Rightarrow$ $\al A.$ is $\rho$-simple
                   $\Rightarrow$ $\al F.$ is $\al G.$-simple.
\end{pro}
\begin{beweis}
The first and third implication above are trivial.

i) To show that the $\al G.$-simplicity of $\al F.$ implies
   the simplicity of $\al A.$ assume that $\al A.$ is not simple: 
   let $\al I.\lhd\al A.$ be a nontrivial closed 2-sided ideal and consider
\[
 \al E._r:=\mr clo._{\|\cdot\|}\;
         \mr span.\{ \al I.\,\al H._\mt D.\mid D\in\wh{\al G.}\}\,.
\]
$\al E._r$ is a closed right ideal in $\al F.$: indeed, recall that
\[
\al F.:=\mr clo._{\|\cdot\|}\;
         \mr span.\{ \al A.\,\al H._\mt D.\mid D\in\wh{\al G.}\}\,
\]
and take $A\in\al A.$. Then
\[
 \al I.\,\al H._\mt D.\;A\al H._{\mt D.'}
                    =\underbrace{\al I.\rho_\mt D.(A)}_{^{\in\al I.}}\;
                         \al H._\mt D.\,\al H._{\mt D.'} \subset\al E._r\,,
                     \quad D,D'\in\wg\,,
\]
where the latter inclusion follows from the fact that the 
tensor product $\al H._\mt D.\,\al H._{\mt D.'}$ can be decomposed in
terms of irreducible algebraic Hilbert spaces: indeed, for 
$H\in\al H._\mt D.$, $H'\in\al H._\mt {D'}.$, there are
$A_{\mt {D,k}.}\in\al A.$ such that
$H\cdot H'=\sum_{\mt {D,k}.} A_{\mt {D,k}.} \Phi_{\mt {D,k}.}$.
Thus we have shown that
\[
 \mr span.\{ \al I.\,\al H._\mt D.\mid D\in\wh{\al G.}\}
 \cdot \mr span.\{ \al A.\,\al H._{\mt D.'}\mid D'\in\wh{\al G.}\}
 \subset\al E._r\,.
\]
Take now $\{F_n\}_n\subset 
\mr span.\{ \al A.\,\al H._{\mt D.'}\mid D'\in\wh{\al G.}\}$
such that $F_n\to F\in\al F.$. Then for any 
$E_0\in\mr span.\{ \al I.\,\al H._\mt D.\mid D\in\wh{\al G.}\}$
we have $E_0F\in\al E._r$, because $\al E._r\ni E_0F_n\to E_0F$.
Similarly one can show that $EF\in\al E._r$ for all
$E\in\al E._r$, $F\in\al F.$, hence $\al E._r$ is a closed
right ideal in $\al F.$.
This implies that $\al E.:=\al E._r\cap\al E._r^*\lhd\al F.$ is a nonzero
closed 2-sided ideal in $\al F.$, which is proper since 
$\1\not\in\al E.$. Finally, $\al E.$ is also $\al G.$-invariant,
because $g(\al I.)=\al I.\subseteq\al A.$ and 
$g(\al H._\mt D.)=\al H._\mt D.$, $g\in\al G.$, $D\in\wh{\al G.}$.
Summing up we have shown that if $\al A.$ is not simple, then
$\al F.$ is not $\al G.$-simple.

ii) To show the last implication, assume that $\al F.$ is not
$\al G.$-simple: let $\al E.\lhd\al F.$ be a nontrivial,
$\al G.$-invariant and closed 2-sided ideal. According to 
Lemma~\ref{IdealGInv} we have that $\al E.$ is also 
$\Pi_\mt D.$-invariant for $D\in\wh{\al G.}$ and we 
define the following closed two-sided ideal in $\al A.$:
\[
 \al I.:=\Pi_\iota(\al E.)=\al E.\cap\al A.\lhd\al A.\,.
\]
We still need to show that $\al I.$ is $\rho$-invariant and
nontrivial. Since $\al E.$ is a two-sided ideal in $\al F.$ 
we have for any $D\in\dg$ and any $X\in\al I.=\al E.\cap\al A.$
\[
 \rho_\mt D.(X)=\sum_{k}\Phi_{\mt {D,k}.} \;X\; \Phi_{\mt {D,k}.}^*
               \in\al E.\cap\al A.\,,
\]
where $\{\Phi_{\mt {D,k}.}\}_\mt k.$ is an orthonormal basis
of $\al H._\mt D.$. Thus $\al I.$ is $\rho$-invariant.
Moreover, $\al I.$ is proper because $\al E.$ is proper:
$\1\not\in\al I.\subset\al E.$. To conclude the proof 
we have to show that $\al I.\not=\{0\}$. For this choose an
element $E'\in\al E.$ with $E'\not=0$. Since $\al E.$ is 
$\Pi_\mt D.$-invariant (cf.~Lemma~\ref{IdealGInv}) we have 
$\Pi_\mt D.(E')\in\al E.$ for all $D\in\dg$ and according
to Corollary~\ref{Compo}~(i) there is a $D\in\dg$ such 
that $E:=\Pi_\mt D.(E')\not=0$. Then we can write
\[
E =\sum_{\mt {k}.} A_{\mt k.} \Phi_{\mt {D,k}.} 
   \quad\mr for~some.\quad A_\mt k.\in\al A.
\]
and at least one of the coefficients does not vanish, say
$A_{\mt {k_0}.}\not=0$. Since $\al E.$ is a 
two-sided ideal in $\al F.$ we have
$E\Phi^*_{\mt {D,k_0}.}\in\al E.$. Then we compute
\begin{eqnarray*}
  \Pi_\iota(E\Phi^*_{\mt {D,k_0}.}) 
       &=& \int_\al G. g(E\Phi^*_{\mt {D,k_0}.})\;dg \\
       &=& \sum_k A_k\int_\al G. 
           g(\Phi_{\mt {D,k}.}) g(\Phi^*_{\mt {D,k_0}.})\;dg \\
       &=& \sum_k A_k \cdot\sum_{\mt {k',i'}.}\left(\int_\al G. 
           U_{\mt {k',k}.}(g)\overline{U_{\mt {i',k_0}.}(g)}\;dg \right)
           \Phi_{\mt {D,k'}.}\;\Phi^*_{\mt {D,i'}.} \\
       &=& \sum_k A_k \cdot\sum_{\mt {k',i'}.} \frac{1}{d(D)}\;
           \delta_{\mt {k',i'}.}\,\delta_{\mt {k,k_0}.}
           \Phi_{\mt {D,k'}.}\;\Phi^*_{\mt {D,i'}.} \\
       &=& \frac{1}{d(D)}\; A_{\mt {k_0}.} \sum_{\mt {k'}.} 
           \Phi_{\mt {D,k'}.}\;\Phi^*_{\mt {D,k'}.} \\
       &=& \frac{1}{d(D)}\; A_{\mt {k_0}.} \;\;\not=\;\; 0\,,
\end{eqnarray*}
where $d(D)$ is the dimension of the representation $U\in D$.
In the previous equations we have used that 
$\mr supp.\al H._\mt D.=\1$, the orthogonality
of the matrix elements $U_{\mt {k',k}.}(g)$ (recall Peter-Weyl's 
Theorem \cite[Theorem~27.40]{bHewittII}) as well as the transformation
$g(\Phi_{\mt {D,k}.})=\sum_{k'}U_{\mt {k',k}.}(g)\Phi_{\mt {D,k'}.}$.
We have thus shown that $\al I.$ is nonzero, since
\[ 
 0\not= \frac{1}{d(D)} A_{k_0}
      =\Pi_\iota(E\Phi^*_{\mt {D,k_0}.})\in\al I.
\]
and the proof is concluded.
\end{beweis}
\begin{cor}
If $\al I.\lhd\al A.$ is $\rho$-invariant,
then the closed right ideal defined by
\[
 \al E.:=\mr clo._{\|\cdot\|}\;
         \mr span.\{ \al I.\,\al H._\mt D.\mid D\in\wh{\al G.}\}
\]
satisfies $\al E.^*=\al E.$, hence is a closed two-sided ideal in 
$\al F.$.
\end{cor}
\begin{beweis}
First we show that $\al E.^*\subseteq\al E.$ by using the conjugation
structure of $\al T._\al G.$: $\al E.^*$ is generated by 
$\al H._\mt D.^* \,\al I.$, $D\in\dg$. From Remark~\ref{ConstrER}~(i)
it follows that $\Phi_{\mt {D,k}.}^*=R_\mt D.^*\,\Phi_{\mt {\ol D.,k}.}$
and from this we obtain
\[
 \al H._\mt D.^* \,\al I.=R_\mt D.^*\,\al H._{\mt {\ol D.}.}\,\al I.
        = R_\mt D.^*\, \rho_{\mt {\ol D.}.}(\al I.)\,\al H._{\mt {\ol D.}.}
        \subseteq\al I.\,\al H._{\mt {\ol D.}.}\subseteq\al E.\,.
\]
This shows the inclusion $\al E.^*\subseteq\al E.$ and from
$\al E.^*\subseteq\al E.=(\al E.^*)^*\subseteq\al E.^*$ we get
the equality $\al E.^*=\al E.$.
\end{beweis}

\section{Minimal and regular Hilbert C*-systems}
\label{MinRegSect}

The DR-theorem associates with a DR-category
$\al T.$ an essentially
unique compact group $\al G.$ \cite[Theorem~6.1]{Doplicher89b}.
In the context of Hilbert C*-systems 
we have a bijective correspondence between
\[ 
    \{\al A.,\al T.\} \quad\mr and.\quad \HS\,,
\]
where $\al A.$ is a unital C*-algebra with trivial center 
$\al Z.=\al A.'\cap\al A.=\C\1$ (and satisfying Property~B) and
$\al T.$ is a DR-category realized as unital endomorphisms of $\al A.$.
$\HS$ is a Hilbert extension of $\al A.$ having trivial relative 
commutant, i.e.~$\al A.'\cap\al F.=\C\1$   (see 
\cite{Doplicher89b,Doplicher90,Doplicher89a,Baumgaertel97}).
This correspondence is connected with the second part of
Proposition~\ref{Isom} which requires a trivial center of $\al A.$. 
The DR-theorem says that in the case of
Hilbert extensions of $\al A.$
with trivial relative commutant, the category $\al T. $
of all canonical endomorphisms can be indeed characterized
intrinsically by their abstract algebraic properties
as endomorphisms of $\al A.$
and a corresponding bijection can be established.

In this section we want to extend such a bijective
correspondence to C*-algebras $\al A.$ with nontrivial center 
$\al Z.\supset\C\1$ and satisfying Property~B. A
first step in this direction is given in \cite{Lledo01a}. 
In this context and due to Proposition~\ref{Isom} one has to 
face the problem that the category $\al T._\al G.$ and
$\al T.$ {\em can not} be isomorphic anymore, since now we have
\[
 \C\1\subset\al Z.\subseteq\al A.'\cap\al F.\,.
\]
We will investigate in the following the class
of Hilbert extensions $\HS$ with compact group $\al G.$ 
and where the relative commutant satisfies the following 
{\em minimality} condition
\begin{equation}\label{MinRC}
  \al A.'\cap\al F.=\al Z.\,.
\end{equation}

In items (i)-(iv) of the introduction we gave several motivations
that justify this choice. Therefore we define
\begin{defi}
A Hilbert C*-system
$\{\al F.,\al G.\}$
is called {\bf minimal} if the condition
\begin{equation}
\label{CZ}
\al A.'\cap\al F.=\al Z. 
\end{equation}
is satisfied.
\end{defi}

\begin{rem}
\begin{itemize}
\item[(i)]
The adjective minimal comes from the property of the relative commutant.
Note that one always has $\al Z.\subseteq\al A.'\cap\al F.$. In the context
of the DR-theory one has also minimal Hilbert C*-systems, because there
$\al Z.=\C\1$ and $\al A.'\cap\al F.=\C\1$.

\item[(ii)] Let $\{\al F.,\al G.\}$ be a C*-dynamical system with 
fixed point algebra $\al A.$ having trivial center 
$\al Z.(\al A.)=\C\1$ and relative commutant satisfying
\[
\al A.'\cap\al F.=\al Z.(\al F.)
    \quad\mr with.\quad\al Z.(\al F.)=\al F.'\cap\al F.\,.
\]
Then $\{\al F.,\al G.\}$ can be obtained by inducing up from an
essentially unique C*-dynamical system
$\{\al F._0,\al G._0\}$, where $\al G._0$ is a closed subgroup
of $\al G.$, the fixed point algebra coincides with $\al A.$ and
the relative commutant is trivial, i.e.
$\al A.'\cap\al F._0=\C\1$ (cf.~\cite[Theorem~1]{Doplicher86}).
For a generalization of this result in the case where 
$\al Z.(\al A.)$ is nontrivial and the corresponding relative
commutant satisfies
\[
 \al A.'\cap\al F.=\al Z.(\al A.)\vee \al Z.(\al F.)
\]
see \cite{pVasselli03a}.
\end{itemize}
\end{rem}

\begin{pro}\label{disj}
Let
$\HS$
be a given Hilbert C*-system. Then
$\al A.'\cap\al F.=\al Z. $ iff
$(\rho_{{\mt D.}},\rho_{{\mt D.}'})=\{0\}$
for
$D\neq D'$,
i.e.~iff the set
$\{\rho_{{\mt D.}}\;\big|\;D\in\wh{\al G.}\}$
is mutually disjoint.
\end{pro}
\begin{beweis}
First note that
$F\in\al A.'\cap\al F.$
iff
$\Pi_{{\mt D.}}F\in (\rho_{{\mt D.}},\iota)\al H._{{\mt D.}}$
for all
$D\in\wh{\al G.}$.
Therefore
$\al A.'\cap\al F.=\al Z. $
iff
$(\rho_{{\mt D.}},\iota)=0$
for all
$D\neq\iota$.
But if this is true then also
$(\rho_{{\mt D.}},\rho_{{\mt D.}'})=0$
follows for all
$D\neq D'$
(see e.g.~\cite[p.~193]{bBaumgaertel92}).
\end{beweis}

Observe that in any Hilbert C*-system,  
for each $\tau\in
\ob\al T. $ the space ${\got H}_\tau:=
\al H._\tau \al Z.$, (where $\al H._\tau $
is a $\al G.$-invariant algebraic Hilbert space)
is a $\al G.$-invariant free right Hilbert
$\al Z.$-module with inner product given as usual by
\[
\langle H_1,H_2\rangle:=H_1^*H_2\in\al Z.\;,\quad H_1,H_2\in\ot H._\tau\,.
\]
Moreover, since for any $\tau\in\ob\al T.$, we have that
$\al Z.\subset (\tau,\tau)$, it is easy to see that
there is a canonical left action of $\al Z.$ on $\ot H._\tau$.
Concretely, there is a natural *-homomorphism
$\al Z.\to\al L._\al Z.(\ot H.)$
(see \cite[Sections~3 and 4]{Lledo97b} for more details).
Hence $\ot H._\tau$ becomes a $\al Z.$-bimodule. Next we state
the isomorphism between the category of canonical endomorphisms
and the corresponding category of $\al Z.$-bimodules.

\begin{pro}
\label{prop0}
Let $\HS$ be a given minimal Hilbert C*-system, where the fixed 
point algebra $\al A.$ has center $\al Z.$. 
Then the category $\al T.$ of all canonical endomorphisms of
$\HS$ is isomorphic to the subcategory $\al M._\al G.$ 
of the category of free Hilbert 
$\al Z.$-bimodules with objects $\ot H.=\al H.\al Z.$, 
$\al H.\in\ob\al T._\al G.$, and 
arrows given by the corresponding $\al G.$-invariant module morphisms
$\al L._\al Z.({\got H}_{1},{\got H}_{2};\al G.)$.

The bijection of objects is given by 
$\rho_\al H.\leftrightarrow \ot H.=\al H.\al Z.$
which satisfies the conditions
\begin{eqnarray*}
\rho_\al H.=(\Ad V)\circ\rho_1+(\Ad W)\circ\rho_2
&\longleftrightarrow&
{\got H}=V{\got H}_1+W{\got H}_2\\[1mm]
\rho_1\circ\rho_2 &\longleftrightarrow&
{\got H}_1\cdot{\got H}_2    \;,
\end{eqnarray*}
where the latter product is the inner tensor product of the 
the Hilbert $\al Z.$-modules w.r.t.~the *-homomorphism
$\al Z.\to\al L._\al Z.(\ot H._2)$.
The bijection on arrows is defined by
\[
\al J.\colon\ \al L._{\al Z.}({\got H}_{1},{\got H}_{2};\al G.)
\to (\rho_1,\rho_2) \quad\mr with.\quad
\al J.(T):=\sum_{j,k}\Psi_{j}Z_{j,k}\Phi_{k}^{\ast}=:\wwh T.\;.
\]
Here
$\{\Psi_{j}\}_{j},\;\{\Phi_{k}\}_{k}$
are orthonormal basis of
${\got H}_{2},{\got H}_{1}$,
respectively, and
$(Z_{j,k})_{j,k}$
is the matrix of the right $\al Z.$-linear operator $T$ from
${\got H}_{1}$
to
${\got H}_{2}$ which intertwines the $\al G.$-action.
\end{pro}
\begin{beweis} 
Note first that the minimality condition and 
Proposition~\ref{HSnonunique} guarantee that the bijection on 
objects is independent of the choice of the algebraic Hilbert
spaces within $\ot H.$, provided these define the same canonical
endomorphism. The rest of the proof 
is in Proposition~3.1 and Section~4 of \cite{Lledo97b}.
\end{beweis}

\begin{rem}
Note that the category $\al M._\al G.$ is a tensor C*-category.
This follows from the fact that 
\[
 \wh T \,Z=Z\,\wh T \;,\quad Z\in\al Z.\,,\; 
                   T\in\al L._\al Z.({\got H}_{1},{\got H}_{2};\al G.)\,,
\]
where $\wh T\in\al A.$ is defined in Remark~\ref{remark1}~(ii).
The previous equation implies 
\[
   T(Z\cdot) =Z\,T(\cdot) \;,\quad Z\in\al Z.\,,\; 
                   T\in\al L._\al Z.({\got H}_{1},{\got H}_{2};\al G.)\,,
\]
and by \cite[p.~268]{Doplicher98} this condition guarantees that
$\al M._\al G.$ is a tensor C*-category. Note that in general
the category of Hilbert $\al Z.$-bimodules with the larger arrow
sets $\al L._\al Z.({\got H}_{1},{\got H}_{2})$ is only a 
{\em semitensor C*-category} (cf.~\cite[Section~2]{Doplicher98}
for the definition of this notion and further 
details).\footnote{We would like to acknowledge an
anonymous referee for mentioning reference 
\cite{Doplicher98} to us.}
\end{rem}

The following result recalls the useful decomposition for a general 
Hilbert $\al Z.$-module $\ot H.\in \ob\,\al M._\al G.$ in terms of 
irreducible ones $\ot H._\mt D.=\al H._\mt D.\al Z.$, $D\in\wh{\al G.}$. 
From Proposition~\ref{prop0} one has equivalently a decomposition of
endomorphism $\rho_\al H.\in\ob\,\al T.$ in terms of the corresponding
irreducibles $\rho_\mt D.\in\mr Irr.\al T.$.

\begin{pro}\label{DecompEnd}
Let $\ot H.\in\mr Ob.\al M._\al G.$ be a $\al G.$-invariant free 
Hilbert $\al Z.$-module in $\al F.$. Then
$\ot H.$ can be decomposed into the following orthogonal direct
sum:
\[
 \ot H.=\mathop{\oplus}_{\KIn{$D$}}\,(\rho_{\KIn{$D$}},\,\rho_\al H.)\,
       {\ot H.}_{\KIn{$D$}}\,.
\]
If $\{ W_{\KIn{$D$},\,l} \}_{l=1}^{m(\KIn{$D$})}$ denotes an \onb{} of
$(\rho_{\KIn{$D$}},\,\rho_\al H.)$, where $m(\In{$D$})$ is the 
multiplicity of $D\in\widehat{\cal G}$ in the decomposition of 
$U_{\KIn{$\cal H$}}$, then the isotypical projection can be written as
\[
 P_{\KIn{$D$}}:=\sum_{l=1}^{\mt m(D).} W_{{\mt D.},\,l}\,W_{{\mt D.},\,l}^*\,.
\]
The canonical endomorphism associated with $\al H.$ is given by
\[
 \rho_\al H.(A) =\sum_{{\mt D.},\,l} W_{{\mt D.},\,l}\cdot \rho_{{\mt D.}}(A) 
                \cdot W_{{\mt D.},\,l}^*\,.
\]
\end{pro}

Recalling the notion of $\rho$ invariant ideals in $\al A.$ stated
in Definition~\ref{GSimple}~(ii) we have the following immediate 
consequence of the previous decomposition result for canonical 
endomorphisms.
\begin{cor}
Let $\HS$ be a minimal Hilbert C*-system with fixed point algebra
$\al A.$ and let $\al I.$ be a closed two-sided ideal in $\al A.$.
Then $\al I.$ is $\rho$-invariant iff $\rho(\al I.)\subseteq\al I.$
for all $\rho\in\mr Ob.\al T.$.
\end{cor}

The Proposition~\ref{prop0} above shows that the canonical endomorphisms
uniquely fix the corresponding $\al Z.$-modules, but not the choice
of the generating algebraic Hilbert spaces. The assumption
of the minimality condition (\ref{CZ}) is crucial here.
From the point of view of the $\al Z.$-modules it is
natural to consider next the following property of Hilbert 
C*-systems:
the existence of a special choice of algebraic Hilbert spaces
within the modules that define
the canonical endomorphisms and which is compatible with products.

\begin{defi}\label{ReguCondi}
A Hilbert C*-system
$\HS$
is called {\bf regular} if there is an assignment
\[
\mr Ob.\al T.\ni\sigma\rightarrow\al H._{\sigma}\in\mr Ob.\al T._\al G.
\]
such that $\sigma=\rho\s{\al H._\sigma}.$,
i.e.~$\sigma$ is the canonical endomorphism of 
the algebraic Hilbert space $\al H._\sigma$,
and which is compatible with products:
\[
\sigma\circ\tau\to\al H._\sigma\cdot\al H._\tau\;.
\]
\end{defi}

\begin{rem}
\begin{itemize}
\item[(i)]
In a minimal Hilbert C*-system regularity means that there is a
``generating" Hilbert space
$\al H._{\tau}\subset{\got H}_{\tau}$
for each $\tau$ (with
${\got H}_{\tau}=\al H._{\tau}\al Z.$)
such that the compatibility relation for products
stated in Definition~\ref{ReguCondi} holds.
If a Hilbert C*-system is minimal and
$\al Z.=\C\un$
then it is necessarily regular.
\item[(ii)] Note that the minimality condition 
and the compactness of the group imply that
the Hilbert modules associated with objects in $\al T.$ are free.
From the point of view of crossed products by endomorphisms 
considered in \cite{pVasselli03b} this corresponds to a special
case (cf.~\cite[Example~4.1]{pVasselli03b}). 
Nevertheless, even in this particular 
situation with a compact group, there are cases where still 
one can not associate a symmetry to the larger category $\al T.$.
In the context of Hilbert C*-modules this means that the 
left action {\it does not} coincide with the corresponding 
right action (in contrast with the situation considered 
in \cite[Section~4]{pVasselli03b}). The existence of a symmetry is
related to the nontriviality of the chain group homomorphism
to be introduced in the next section (see also
Proposition~\ref{ABasesZ}, Remark~\ref{NotSym}~(i)
and Section~\ref{Conclu}).
\end{itemize}
\end{rem}

\begin{lem}\label{ONBNonIrr}
Let $\HS$ be a minimal and regular Hilbert C*-system.
For $\sigma,\tau\in\mr Ob.\al T.$ and $\rho,\rho'\in\mr Irr.\al T.$
put
\[
  (\sigma,\,\tau)_\c
       :=\al J.\left(\al L.\s{\al G.}.(\al H._\sigma,\al H._\tau)\right)
         \subseteq (\sigma,\tau)\;,\quad
   (\rho,\,\tau)_\c
       :=\al J.\left(\al L.\s{\al G.}.(\al H._\rho,\al H._\tau)\right)
       \subseteq (\rho,\,\tau)\;,\;\mr etc..,
\]
(cf. Remark~\ref{remark3}~(iii)).
If $\{ W_{\tau,\,\rho,\,k} \}_{k}\subset (\rho,\,\tau)_\c\,$,
$\{ W_{\sigma,\,\rho',\,k'} \}_{k'}\subset (\rho',\,\sigma)_\c$
are \onb{}, then 
\[
\{ W_{\tau,\,\rho,\,k}\cdot W_{\sigma,\,\rho,\,k'}^* \}_{\rho,k,k'}
\]
is an \onb{} of $(\sigma,\,\tau)_\c$, where the $\rho$'s are 
those irreducibles appearing in the decomposition of $\sigma$ 
{\em and} $\tau$ (cf.~Proposition~\ref{DecompEnd}).
\end{lem}
\begin{beweis}
From the orthonormality relations of 
$\{ W_{\tau,\,\rho,\,k} \}_{k}$ and 
$\{ W_{\sigma,\,\rho',\,k'} \}_{k'}$ 
and from the disjointness relation for irreducible 
objects in Proposition~\ref{disj},
it follows directly that
$\{ W_{\tau,\,\rho,\,k}\cdot W_{\sigma,\,\rho,\,k'}^* \}_{\rho,k,k'}$
is an orthonormal system in $(\sigma,\,\tau)_\c$. For any 
$A\in(\sigma,\,\tau)_\c$ put
\[
 \lambda_{\rho,k,k'}
     := W_{\tau,\,\rho,\,k}^*\,A\,W_{\sigma,\,\rho,\,k'}
     \in (\rho,\rho)_\c=\C\1\,,
\]
where the last equation follows from the fact that the algebraic
Hilbert space corresponding to $\rho$ carries an irreducible 
representation. Then it is immediate to verify that
\[
 A=\sum_{\rho,k,k'}\lambda_{\rho,k,k'} 
         \,W_{\tau,\,\rho,\,k}\,W_{\sigma,\,\rho,\,k'}^*\,,
\]
hence $\{ W_{\tau,\,\rho,\,k}\cdot W_{\sigma,\,\rho,\,k'}^* \}_{\rho,k,k'}$
is an  \onb{} in $(\sigma,\,\tau)_\c$.
\end{beweis}

\begin{teo}
\label{Teo1}
Let $\HS$ be a minimal and regular Hilbert C*-system
(where the fixed point algebra $\al A.$ satisfies Property~B). Then
$\al T.$ contains a DR-subcategory $\al T._\c$ with
the same objects, $\ob\al T._\c=\ob\al T.,$ and arrows
$(\sigma,\,\tau)_\c\subseteq(\sigma,\tau)$ such that:
\begin{equation}\label{GeneratingArrows}
\quad(\sigma,\tau)=(\sigma,\tau)\s\c.\;\sigma(\al Z. )
             =\tau(\al Z.)\,(\sigma,\tau)\s\c.\;\sigma(\al Z.)\;.
\end{equation}
\end{teo}
\begin{beweis}
Let
$\al H.,\al K.$
are $\al G.$-invariant algebraic Hilbert spaces. Recall that the
isometry
$\al J.:\al L.(\al H.,\al K.)\to\al F.$
has the property
\[
(\al H.,\al K.)=\al J.(\al L._{\al G.}(\al H.,\al K.))
                \subseteq (\rho_{\al H.},\rho_{\al K.})\,,
\]
(cf. Remark~\ref{remark3}~(iii)). Now let
$\sigma\to\al H._{\sigma}$
be the 
assignment given in Definition~\ref{ReguCondi}
for regular Hilbert C*-systems, and put
\[
  (\sigma,\,\tau)\s\c.:=(\al H._\sigma,\al H._\tau)\;.
\]
Then the definitions of the symmetry
$\epsilon(\cdot,\cdot)$,
the conjugates
$\ol\sigma.$
and their intertwiners
$R_\sigma,\,S_\sigma$
are as follows (cf.~Remark~\ref{ConstrER}~(i)):
\[
\epsilon(\sigma,\tau):=\sum_{j,k}\Psi_k\Phi_j\Psi^{\ast}_{k}\Phi^{\ast}_{j},
\]
where
$\{\Phi_j\}_j$ resp.~$\{\Psi_k\}_k$
are orthonormal basis of
$\al H._\sigma$
resp.~$\al H._\tau$.
We also define:
$\ol\sigma.
:=\rho_{\al K.},$
where $\al K.$ carries the representation of $\al G.$ which is conjugated
to the representation on
$\al H._\sigma$.
Recall that all finite-dimensional representations of $\al G.$
can be realized by some $\al G.$-invariant algebraic Hilbert space.
Moreover, $\al K.$ is chosen as the ``distinguished"
Hilbert space according to the assumption of regularity. Let
$\{\Omega_j\}_j$
be an orthonormal basis of $\al K.$ and put
\[
R_\sigma:=\sum_{j}\Omega_j\Phi_j,\qquad S_\sigma:=
\epsilon(\ol\sigma.,\sigma)R_\sigma.
\]
With these choices it is easy to verify that $\al T._\c$ is indeed 
a DR-subcategory of $\al T.$ (cf.~Definition~\ref{DRCat}).

It remains to show Eq.~(\ref{GeneratingArrows}). The inclusion
$\quad(\sigma,\tau)\supseteq(\sigma,\tau)\s\c.\;\sigma(\al Z. )$
follows immediately from the fact that $(\sigma,\tau)$ is a 
right $\sigma(\al Z.)$-module.
To show the reverse inclusion let
$\{ W_{\tau,\,\rho,\,k} \}_{k}\subset (\rho,\,\tau)_\c\,$,
$\{ W_{\sigma,\,\rho',\,k'} \}_{k'}\subset (\rho',\,\sigma)_\c$
be \onb{} as in Lemma~\ref{ONBNonIrr}.
Take $A\in(\sigma,\tau)$ and define
\[
 \rho(Z_{\rho,k,k'})
     := W_{\tau,\,\rho,\,k}^*\,A\,W_{\sigma,\,\rho,\,k'}
     \in (\rho,\rho)=\rho(\al Z.)\,.
\]
Then 
\begin{eqnarray*}
 A &=&\sum_{\rho,k,k'} 
   \,W_{\tau,\,\rho,\,k}\;\rho(Z_{\rho,k,k'})\;W_{\sigma,\,\rho,\,k'}^*\\
   &=&\sum_{\rho,k,k'} 
   \,W_{\tau,\,\rho,\,k}\,W_{\sigma,\,\rho,\,k'}^* \,\sigma(Z_{\rho,k,k'})
   \in(\sigma,\tau)\s\c.\;\sigma(\al Z. )\,,
\end{eqnarray*}
and the proof of Eq.~(\ref{GeneratingArrows}) is completed.
\end{beweis}

\begin{rem}\label{Admi}
\begin{itemize}
\item[(i)]
Recall that
the category $\al T._\al G.$ introduced in Subsection~\ref{GCanEnd}
is an example of a DR-category (recall Remark~\ref{remark3})
and in fact it plays the role of the 
subcategory $\al T._{\c}$ as a subcategory of the
in general larger category 
$\al T.$ of canonical endomorphisms of the minimal and regular
Hilbert C*-system.

\item[(ii)]
The assumptions of the previous theorem imply that each basis of
$(\sigma,\tau)_{\c}$
is simultaneously a module basis of
$(\sigma,\tau)$
modulo
$\sigma(\al Z. )$
as a right module,
i.e.~the module
$(\sigma,\tau)$
is free.

\item[(iii)]
For simplicity we will sometimes call
a DR-subcategory $\al T._\c$ of $\al T.$ satisfying the properties
of Theorem~\ref{Teo1} {\bf admissible}.

\item[(iv)] Note that the properties (P.2)-(P.4) (with the exception 
of P.2.6) in \cite[Section~2]{Lledo01a} are contained in the 
assumption that $\al T._\c$ is a DR-subcategory of $\al T.$ 
(cf.~Definition~\ref{DRCat}).
\end{itemize}
\end{rem}

Let
$\ob\al T.\ni\rho\to V_\rho\in (\rho,\rho)$
be a choice of unitaries such that
\begin{equation}
\label{Vrs}
V_{\rho\circ\sigma}=V_\rho\times V_{\sigma}.
\end{equation}
Note that (\ref{Vrs}) implies that
$V_\iota =\un,$
because
$V_\sigma=V_{\iota\sigma}=V_\iota\times V_\sigma=V_\iota V_\sigma.$
This choice allows one to define from the subcategory
$\al T._\c$
of $\al T.$ another subcategory
$\al T.'_\c$ of $\al T.$
with the same properties as $\al T._\c.$ To do this, put
\begin{equation}
\label{rsc}
(\rho,\sigma)'_\c:=V_\sigma(\rho,\sigma)_{\c}V^{\ast}_\rho\subset (\rho,\sigma)
\end{equation}
and for the corresponding permutation structure
$\epsilon'(\rho,\sigma)$ for $\al T.'_\c$
one takes:
\begin{equation}
\label{ers}
\epsilon'(\rho,\sigma)
:=(V_\sigma\times V_{\rho})\cdot\epsilon(\rho,\sigma)\cdot(V_\rho\times 
V_\sigma)^{\ast}.
\end{equation}
It is easy to check that $\epsilon'$ defines a permutation 
structure (cf.~Proposition~\ref{PermuStructu}). The
corresponding conjugates $R'_\rho$ are defined by
\begin{equation}
\label{Rpr}
R'_\rho:=V_{\ol\rho.\rho}R_\rho,\qquad S'_\rho
:=\epsilon'(\rho,\ol\rho.)R'_\rho
\end{equation}
(recall also Proposition~\ref{ConjuStructu}). Then
it is straightforward to verify that the new subcategory $\al T.'_\c$
also satisfies the assumptions of Theorem~\ref{Teo1}.

This suggests the following definition of equivalence between
subcategories.
\begin{defi}
Two admissible DR-subcategories 
$\al T._\c$ and $\al T.'_\c$ of $\al T.$
are called {\bf equivalent}
if there is an assignment
\[
\ob\al T.\ni\rho\to V_\rho\in (\rho,\rho),
 \quad\mr with.~V_\rho~\mr unitary~satisfying.
\quad V_{\rho\circ\sigma}=V_\rho\times V_\sigma,
\]
and such that the equations (\ref{rsc}), (\ref{ers}) and 
(\ref{Rpr}) hold.
\end{defi}

The converse of Theorem~\ref{Teo1} gives our main duality theorem.
The proof, which is constructive, is divided into several steps
(see Subsection~\ref{StepProof} below).
\begin{teo}
\label{Teo2}
Let $\al T.$ be a tensor C*-category of unital endomorphisms
of $\al A.$ and let $\al T._\c$ be an admissible DR-subcategory.
Then there is a minimal and regular Hilbert
extension $\HS$ of $\al A.$, where $\al G.$ is the compact
group assigned to the DR-category $\al T._\c$ and
$\al T.$ is isomorphic 
to the category of all canonical endomorphisms of $\HS$.
Moreover, if $\al T.\s\c.,\;\al T.'_\c$ are two admissible 
DR-subcategories of
$\al T.$, 
then the corresponding Hilbert
extensions are $\al A.$-module isomorphic iff $\al T.\s\c.$
is equivalent to $\al T.'_\c.$
\end{teo}

The previous result says that for minimal and regular Hilbert C*-systems  
there is an intrinsic characterization of the category of all
canonical endomorphisms in terms of 
$\al A.$
only. Moreover,
up to $\al A.$-module isomorphisms, there is a bijection
between minimal and regular Hilbert extensions 
and tensor C*-categories $\al T.$ of unital endomorphisms of $\al A.$
with admissible DR-subcategories.  

Note that Theorem~\ref{Teo2}
is an immediate generalization of the
DR-theorem (mentioned at the beginning of this section) for
the case that
$\al Z. \supset\C\un$,
i.e.~if
$\al Z. =\C\un$
then $\al T.$
itself is admissible (hence a DR-category) and the corresponding
Hilbert extensions have trivial relative commutant.
Notice also that from the assumption of the existence of an
admissible DR-subcategory it follows that $\al A.$ satisfies 
Property~B.

\begin{rem}\label{commut}
{\bf (Abelian groups)} 
In the case that $\al G.$
is abelian and compact the preceding structure simplifies radically
(see Remark~\ref{ChangeAb}). Specifically,
$\wh{\al G.}$
is a discrete abelian group (the character group),
each $\al H._{\mt D.}$, $D\in\wh{\al G.}$, 
is one-dimensional with a generating
unitary $U_{\mt D.},$ hence the canonical endomorphisms 
$\rho_{\al H._{{\mt D.}}}$ 
(denoted briefly by
$\rho_{{\mt D.}}$,
see Remark~\ref{remark2}(ii))
are in fact automorphisms
(necessarily outer on $\al A.).$
Since
$\rho_{{\mt D.}_{1}}\circ\rho_{{\mt D.}_{2}}=\rho_{{\mt D.}_{1}{\mt D.}_{2}}$
in this case the set
$\Gamma$
of all canonical endomorphisms
$\rho_{\al H._{{\mt D.}}}$
is a group with the property
\[
 \wh{\al G.}\cong \Gamma/\hbox{int}\,\al A.\,.
\]
Therefore, it is not necessary to take into account
direct sums or subobjects in this case and one can drop Property~B 
as an assumption on $\al A.$.

If in addition $\al Z. =\C\un$
the permutators $\epsilon$ (restricted to
$\wh{\al G.}\times\wh{\al G.}$)
are elements of the second cohomology group
$H^{2}(\wh{\al G.})$
and
\[
U_{{\mt D.}_{1}}\cdot U_{{\mt D.}_{2}}=\omega({\ms D.}_{1},{\ms D.}_{2})
U_{{\mt D.}_{1}\circ{\mt D.}_{2}},
\]
where
\[
\epsilon(D_{1},D_{2})=\frac{\omega(D_{1},D_{2})}
{\omega(D_{2},D_{1})}
\]
and $\omega$ is a corresponding 2-cocycle. The field algebra $\al F.$
is just the $\omega$-twisted discrete crossed product
of $\al A.$ with $\wh{\al G.}$
(see e.g.~\cite[p.~86 ff.]{bBaumgaertel95} for details).
For the case
$\al Z. \supset\C\un$
see \cite{Lledo02}. (The minimal case is not specially mentioned there.)
\end{rem}

\subsection{Proof of Theorem~\ref{Teo2}}\label{StepProof}

The construction of the Hilbert extension $\HS$ associated with the 
pair $\al T._\c<\al T.$, where $\al T._\c$ is an admissible DR-subcategory,
is done in several steps which are adapted from 
\cite[Sections~3-6]{Baumgaertel97} and \cite{Lledo01a,Lledo97b}. 
The strategy is to define first
a left $\al A.$-module $\al F._0$ using the dimension function defined
in (\ref{DimFkt}). The structure of $\al F._0$ is then gradually enriched
by making use of the properties of $\al T._\c$.

\paragraph{1.~Step:} {\em The bimodule $\al F._0$, 
  algebraic Hilbert spaces and free Hilbert $\al Z.$-modules.}

To each $\rho\in\mathrm{Irr}_{0}\,\mathcal{T}$ (cf.~Remark~\ref{AssocD})
we assign a Hilbert space
$\mathcal{H}_{\rho}$ with
$\dim\,\mathcal{H}_{\rho}=d(\rho)$
and, using orthonormal bases $\{\Phi_{\rho j}\}_{j}$ of
$\mathcal{H}_{\rho}$, we define the left $\mathcal{A}$-module
\[
\mathcal{F}_{0}:=\Big\{\sum_{\rho, j}A_{\rho j}\Phi_{\rho j}\mid
               A_{\rho j}\in\mathcal{A},\,\mbox{finite sum}\Big\}\;,
\]
where the $\{\Phi_{\rho j}\}_{\rho j}$ form an
$\mathcal{A}$-module basis of $\mathcal{F}_{0}$.
$\mathcal{F}_{0}$ is independent of the special choice of the bases
$\{\Phi_{\rho j}\}_{j}$ of $\mathcal{H}_{\rho}$ and 
putting $\Phi_{\rho j}A:=\rho(A)\Phi_{\rho j}$,
$\mathcal{F}_{0}$ turns out to be a bimodule.

Next we introduce Hilbert spaces for any object in the category. 
For this purpose recall that 
$\rho < \alpha$ means that $\rho$ is a subobject of $\alpha$
and that $(\rho,\alpha)_\c$ is an algebraic Hilbert space
in $\al A.$ whose dimension coincides with the multiplicity
of $\rho$ in the decomposition of $\alpha$ 
(cf.~Proposition~\ref{Dec}~(iv)). Then we have
\begin{equation}\label{HS}
\mathcal{H}_{\alpha}:=\mathop{\oplus}_{\rho<\alpha}(\rho,\alpha)_{\c}
\al H._\rho
\quad\mathrm{and}\quad \mathcal{H}_{\alpha}\subset 
\mathcal{F}_{0},\,\alpha \in\mathrm{Ob}\,\mathcal{T}\,,
\end{equation}
as well as the right-$\mathcal{Z}$-Hilbert modules
\[
{\ot H.}_{\rho}:=\mathcal{H}_{\rho}\mathcal{Z}=\mathcal{Z}\mathcal{H}_{\rho}
\quad\mbox{and}\quad
{\ot H.}_{\alpha}:=\mathop{\oplus}_{\rho<\alpha}
                   (\rho,\alpha){\ot H.}_{\rho}
                  =\mathcal{H}_{\alpha}\mathcal{Z},
\]
with the corresponding $\mathcal{Z}$-scalar product
\[
\langle X,Y\rangle_{\alpha}:= \sum_{\rho,j}\rho^{-1}(X_{\rho j}^{\ast}
Y_{\rho j})\,,\;\mathrm{where}
\]
\[
X=\sum_{\rho,j}X_{\rho j}\Phi_{\rho
j},\,X_{\rho j}\in(\rho,\alpha),\,Y=\sum_{\rho,j}Y_{\rho
j}\Phi_{\rho j},\,Y_{\rho j}\in(\rho,\alpha).
\]
The preceding comments show that we have established the
following functor $\ot F.$ between the categories
$\mathcal{T}$ (resp.~$\mathcal{T}_{\c}$)
and the corresponding category of Hilbert $\mathcal{Z}$-modules
(resp.~Hilbert spaces); (cf.~e.g.~\cite[Section~4]{Lledo97b} and 
\cite[Corollary~3.3]{Baumgaertel97}).

\begin{lem}\label{Funct}
The functor $\ot F.$ given by
\[
\mathrm{Ob}\,\mathcal{T}\ni\alpha\mapsto {\ot H.}_{\alpha}\subset
\mathcal{F}_{0}\quad\mathrm{and}\quad(\alpha,\beta)\ni A\mapsto
\ot F.(A)\in\mathcal{L}_\mathcal{Z}({\ot H.}_{\alpha}
\rightarrow {\ot H.}_{\beta})\,,
\]
where $\ot F.(A)X:=AX,\,X\in{\ot H.}_{\alpha}$,
defines an  isomorphism between the corresponding categories and
$\ot F.(A^{\ast})$ is the module adjoint w.r.t.~$\langle\cdot,\cdot\rangle_{\alpha}$. 
Similarly, one can apply $\ot F.$ to $\mathcal{T}_{\c}$
in order to obtain the associated subcategory of algebraic
Hilbert spaces $\mathcal{H}_{\alpha}$ and arrows
$\ot F.((\alpha,\beta)_{\c})\subset \mathcal{L}(\mathcal{H}_{\alpha}
\rightarrow \mathcal{H}_{\beta})$.
\end{lem}
\begin{beweis}
Similar as in \cite[p.~791~ff]{Lledo97b}.
\end{beweis}

\paragraph{2.~Step:} {\em Product and *-structure on $\al F._0$.}
 
We can now apply results in \cite{Baumgaertel97} to the category
$\ot F.(\al T._\c)$

\begin{lem}\label{Prod}
There exists a product structure on
$\mathcal{F}_{0}$ with the properties
\[
\mathrm{span}\,\{\Phi\cdot\Psi\mid\Phi\in\mathcal{H}_{\alpha},\,\Psi\in
\mathcal{H}_{\beta}\}=\mathcal{H}_{\alpha\beta},
\]
\[
\epsilon(\alpha,\beta)\Phi\Psi=\Psi\Phi,\quad \Phi\in\mathcal{H}_{\alpha},\,
\Psi\in\mathcal{H}_{\beta},
\]
\[
\langle XY,X'Y'\rangle_{\alpha\beta}=\langle X,X'\rangle_{\alpha}\cdot
\langle Y,Y'\rangle_{\beta},\quad X,X'\in
 \mathcal{H}_{\alpha},\,Y,Y'\in\mathcal{H}_{\beta}.
\]
\end{lem}

Note that for orthonormal bases
$\{\Phi_{j}\}_{j},\,\{\Psi_{k}\}_{k}$ of
$\mathcal{H}_{\alpha},\,\mathcal{H}_{\beta}$,
respectively, we obtain from Lemma~\ref{Prod} that
\[
\epsilon(\alpha,\beta)=\sum_{j,k}\Psi_{k}\Phi_{j}
\Psi_{k}^{\ast}\Phi_{j}^{\ast}\,.
\]

As in \cite[Section 5]{Baumgaertel97} we introduce the notion of
a conjugated basis $\Phi_{\overline{\alpha}j}$
of $\mathcal{H}_{\overline{\alpha}}$ w.r.t.~an orthonormal basis
$\Phi_{\alpha j}$ of $\mathcal{H}_{\alpha}$ such that
$R_{\alpha}=\sum_{j}\Phi_{\overline{\alpha}j}\Phi_{\alpha j}$.
This is necessary in order to put a compatible *-structure on 
$\mathcal{F}_{0}$.

\begin{lem}
Let $\Phi_{\overline{\rho},j}$ be a conjugated basis corresponding 
to the basis $\Phi_{\rho,j},\,\rho\in\mathrm{Irr}_{0}\,\mathcal{T}$, and
define $\Phi_{\rho j}^{\ast}:=R_{\rho}^{\ast}\Phi_{\overline{\rho}j},
\,j=1,2,...,d(\rho)$. Then $\mathcal{F}_{0}$
turns into a *-algebra. The Hilbert spaces $\mathcal{H}_{\alpha}$
and the corresponding modules ${\ot H.}_{\alpha}$
are algebraic,  i.e.
\[
\langle X,Y\rangle_{\alpha}= X^{\ast}Y,\quad X,Y\in{\ot H.}_{\alpha}.
\]
The objects $\alpha\in\mathrm{Ob}\,\mathcal{T}$
are identified as canonical endomorphisms
\[
\alpha(A)=\sum_{j=1}^{d(\alpha)}\Phi_{\alpha j}A\Phi_{\alpha j}^{\ast}.
\]
\end{lem}

\paragraph{3.~Step:} {\em C*-norm and completion of $\al F._0$.}

In $\mathcal{F}_{0}$ one has natural projections
$\Pi_{\rho}$ onto the $\rho$-component of the decomposition:
\[
\Pi_{\rho}\Big(\sum_{\sigma,j}A_{\sigma j}\Phi_{\sigma j}\Big):=
\sum_{j=1}^{d(\rho)}A_{\rho j}\Phi_{\rho j},\quad
\rho\in\mathrm{Irr}_0\,\mathcal{T}.
\]
To put a C*-norm
$\Vert\cdot\Vert_{\ast}$
we argue as in \cite[Section~6]{Baumgaertel97}. Its construction 
is essentially based on the following $\mathcal{A}$-scalar product on
$\mathcal{F}_{0}$
\[
\langle F,G\rangle :=\sum_{\rho,j}\frac{1}{d(\rho)}A_{\rho
j}B_{\rho j}^{\ast},\,\mathrm{for}\; F:=\sum_{\rho,j}A_{\rho j}\Phi_{\rho j},
\,G:=\sum_{\rho,j}B_{\rho j}\Phi_{\rho j},
\]

\begin{lem}
The scalar product $\langle\cdot,\cdot\rangle$ satisfies
$\langle F,G\rangle =\Pi_{\iota}FG^{\ast}$ and
$\Pi_{\rho}$ is selfadjoint w.r.t.~$\langle\cdot,\cdot\rangle$.
The projections $\Pi_{\rho}$
and the scalar product have continuous extensions to
$\mathcal{F}:=\mathrm{clo}_{\Vert\cdot\Vert_{\ast}}\mathcal{F}_{0}$ and
$\Pi_{\rho}\mathcal{F}=\mathrm{span}\,\{\mathcal{A}\mathcal{H}_{\rho}\}.$
\end{lem}

\paragraph{4.~Step:} {\em The compact group $\al G.$.}

Finally, the symmetry group w.r.t.~$\langle\cdot,\cdot\rangle$
is defined by the subgroup of all automorphisms
$g\in\mathrm{aut}\,\mathcal{F}$ satisfying
$\langle gF_{1},gF_{2}\rangle=\langle F_{1},F_{2}\rangle$.
It leads to

\begin{lem} 
The symmetry group coincides with the stabilizer
$\mathrm{stab}\,\mathcal{A}$ of $\mathcal{A}$ and the modules
${\ot H.}_{\alpha}$ are invariant w.r.t.~$\mathrm{stab}\,\mathcal{A}$.
\end{lem}
\begin{beweis}
Use \cite[Lemma~7.1]{Lledo97b}
(cf.~also with \cite[Section~6]{Baumgaertel97}).
\end{beweis}

This suggests to consider the subgroup
$\mathcal{G}\subseteq\mathrm{stab}\,\mathcal{A}$
consisting of all elements of $\mathrm{stab}\,\mathcal{A}$
leaving even the Hilbert spaces $\mathcal{H}_{\alpha}$
invariant. Then it turns out that the pair
$\{\mathcal{F},\mathcal{G}\}$ is a Hilbert extension of
$\mathcal{A}$. 

\begin{lem}\label{LastLem}
$\mathcal{G}$ is compact and the spectrum $\mathrm{spec}\,\mathcal{G}$ on
$\mathcal{F}$ coincides with the dual $\hat{\mathcal{G}}$. For
$\rho\in\mathrm{Irr}\,\mathcal{T}$ the Hilbert spaces
$\mathcal{H}_{\rho}$ are irreducible w.r.t.~$\mathcal{G},$
i.e.~there is a bijection
$\mathrm{Irr}_0\,\mathcal{T}\ni\rho\leftrightarrow D\in\hat{\mathcal{G}}$.
Moreover $\mathcal{A}$ coincides with the fixed point algebra 
of the action of $\mathcal{G}$ in $\mathcal{F}$ and
$\mathcal{A}'\cap \mathcal{F}=\mathcal{Z}$.
\end{lem}

\begin{rem}
The DR-Theorem shows that any DR-category $\al T._\c$ has an 
unique (modulo isomorphisms) compact group $\al G._\mrt DR.$ associated with
it. We will show here that $\al G._\mrt DR.$ coincides with the 
compact group $\al G.$ obtained as a subgroup of $\mr stab.\al A.$ in
the previous lemma.

For any $\alpha\in\mr Ob.\al T._\c=\mr Ob.\al T.$ we can assign 
a Hilbert space as in the first step above:
\[
\alpha\to {\al H.}_{\alpha}:=\mathop{\oplus}_{\rho<\alpha}
                   (\rho,\alpha)_\c{\al H.}_{\rho}\,.
\]
These Hilbert spaces together with the corresponding arrows
(cf.~Lemma~\ref{Funct}) defines a TK-category $\al T._\mrt TK.$
for $\al G._\mrt DR.$ and by construction it is clear that we
have the isomorphism $\al T._\mrt TK.\cong\al T._\c$.
But at the same time $\al T._\mrt TK.$
is a Tannaka-Krein category for $\al G.$, since by
Lemma~\ref{LastLem}, $\al G.$ acts on $\al H._\rho$ irreducibly,
$\rho\in\mr Irr.\al T.$, and invariantly on $\al H._\alpha$,
$\alpha\in\mr Ob.\al T.$ (recall that $(\rho,\alpha)_\c\subset\al A.$
and $\al G.\subset\mr stab.\al A.$). Therefore
$\al G._\mrt DR.$  
and $\al G.$ are isomorphic, because they have the same 
TK-category.
\end{rem}

\paragraph{5.~Step:} {\em Uniqueness of the Hilbert extension.}

First assume that the subcategories
$\mathcal{T}_{\c}$ and $\mathcal{T}_{\c}'$
are equivalent. We consider the Hilbert extension $\mathcal{F}$
assigned to $\mathcal{T}_{\c}$.
The corresponding invariant Hilbert spaces are given by (\ref{HS}). 
Now we change these Hilbert spaces by
\[
\mathcal{H}_{\alpha}\rightarrow V_{\alpha}\mathcal{H}_{\alpha}=:\al H._\alpha'\,.
\]
Using the function $\ot F.$ of Lemma~\ref{Funct} so that
$\mathcal{L}_\mathcal{G}(\mathcal{H}_{\alpha}\rightarrow\mathcal{H}_{\beta})
:=\ot F.((\alpha,\beta)_{\c})\cong (\alpha,\beta)_{\c}$
we obtain
\begin{equation}\label{Intw}
\mathcal{L}_\mathcal{G}(V_{\alpha}\mathcal{H}_{\alpha}\rightarrow
V_{\beta}\mathcal{H}_{\beta})=
V_{\beta}\mathcal{L}_\mathcal{G}(\mathcal{H}_{\alpha}\rightarrow
\mathcal{H}_{\beta})V_{\alpha}^{\ast}\cong
(\alpha,\beta)_{\c}'.
\end{equation}
Further, w.r.t.~the ``new Hilbert spaces" we obtain the `primed'
permutators and conjugates of the second subcategory. This means, it
is sufficient to prove that if the subcategory
$\mathcal{T}_{\c}$
is given, then two Hilbert extensions, assigned to
$(\mathcal{T},\mathcal{T}_\c)$
according to the first part of the theorem, are always
$\mathcal{A}$-module isomorphic. Now let
$\mathcal{F}_{1},\mathcal{F}_{2}$
be two Hilbert extensions assigned to
$\mathcal{T}_{\c}$.
For $\rho\in\mathrm{Irr}_{0}\,\mathcal{T}$ let
$\{\Phi_{\rho j}^{1}\}_j$,$\{\Phi_{\rho j}^{2}\}_k$
be orthonormal bases of the Hilbert spaces
$\mathcal{H}_{\rho}^{1},\mathcal{H}_{\rho}^{2}$,
respectively. Then
\[
\Phi_{\rho j}^{r}\cdot\Phi_{\sigma k}^{r}=
\sum_{\tau, l}K_{\rho j\sigma k}^{\tau l}\Phi_{\tau l}^{r},\quad
K_{\rho j\sigma k}^{\tau l}\in(\tau,\rho\sigma)_{\c},\,r=1,2.
\]
Therefore the definition
\[
\al J.(\sum_{\rho,j}A_{\rho j}\Phi_{\rho j}^{1}):=
\sum_{\rho,j}A_{\rho j}\Phi_{\rho j}^{2}
\]
is easily seen to extend to an
$\mathcal{A}$-module isomorphism from
$\mathcal{F}_{1}$ onto $\mathcal{F}_{2}$ (see 
\cite[p.~203~ff.]{bBaumgaertel92}).

Second, we assume that the Hilbert extensions
$\mathcal{F}_{1},\mathcal{F}_{2}$ assigned to
$\mathcal{T}_{\c}^{1},\mathcal{T}_{\c}^{2}$,
respectively, are 
$\mathcal{A}$-module isomorphic. The 
$\mathcal{G}$-invariant Hilbert spaces are given by (\ref{HS}). 
Now let $\al J.$ be an $\mathcal{A}$-module isomorphism 
$\al J.\colon\ \mathcal{F}_{1}\rightarrow \mathcal{F}_{2}$ so that
\[
\al J.(\mathcal{H}_{\alpha}^{1})=\mathop{\oplus}_{\rho<\alpha}
(\rho,\alpha)_{\c}^{1}\,\al J.(\mathcal{H}_{\rho}^{1})
\]
and again the $\al J.(\mathcal{H}_{\alpha}^{1})$
form a system of $\mathcal{G}$-invariant Hilbert spaces in 
$\mathcal{F}_{2}$. Further we have the system
$\mathcal{H}_{\alpha}^{2}$ in $\mathcal{F}_{2}$.
That is, to each $\alpha$ we obtain two
$\mathcal{G}$-invariant Hilbert spaces $\mathcal{H}_{\alpha}^{2}$
and $\al J.(\mathcal{H}_{\alpha}^{1})$
that are contained in the Hilbert module $\ot H._\alpha^2$.
Let $\{\Phi_{\alpha,j}\}_j$,$\{\Psi_{\alpha,j}\}_j$ be orthonormal
bases of $\al J.(\mathcal{H}_{\alpha}^{1})$,$\al H._\alpha^2$,
respectively. Then obviously 
$V_\alpha:=\sum_j \Psi_{\alpha,j}\Phi_{\alpha,j}^*$ is a unitary
with $V_\alpha\in(\alpha,\alpha)$ and $\al H._\alpha^2=
V_\alpha\,\al J.(\mathcal{H}_{\alpha}^{1})$. Further, 
for $X\in\al H.^1_\alpha$,$Y\in\al H.^1_\beta$ (hence 
$XY\in\al H.^1_{\alpha\beta}$) we have 
\[
 V_\alpha \al J.(X) V_\beta \al J.(Y)
    = V_\alpha\alpha(V_\beta)\,\al J.(XY) 
    = V_{\alpha\circ\beta}\,\al J.(XY)\,,
\]
and this implies $V_{\alpha\circ\beta}=V_{\alpha}\times V_{\beta}$.
Finally, we argue as in (\ref{Intw}) to obtain
\[
 V_\beta\, (\alpha,\beta)_\C^1 \,V_\alpha=(\alpha,\beta)_\C^2.
\]
and the latter equation implies Eqs.~(\ref{rsc})-(\ref{Rpr}).

\section{Minimal Hilbert C*-systems for nonabelian groups}

Let
$\{\al F.,\al G.\}$
be a minimal Hilbert C*-system with 
$\al G.$
nonabelian and such that the fixed point algebra
$\al A.$
satisfies Property~B. Recall from Remark~\ref{remark3}~(i)
that $\rho\in\mr Irr.\al T.$ if $(\rho,\rho)=\rho(\al Z.)$.
In this case one has trivially that $\al Z.\subseteq\rho(\al Z.)$.
We will next show that for irreducible endomorphisms the
previous inclusion is actually an equality.

\begin{pro}\label{ProAutZS}
For any  $\rho\in\mr Irr.\al T.$ we have 
$\rho(\al Z.)=\al Z.$ and $\rho\rest\al Z.$
is an automorphism of $\al Z.$, i.e.~$\rho\in\mr aut.\al Z.$.
\end{pro}
\begin{beweis}
Choose an irreducible endomorphisms $\rho\in\mr Irr.\al T.$ and
recall that in this case the conjugated endomorphism
is also irreducible, i.e.~$(\ol\rho.,\ol\rho.)=\ol \rho.(\al Z.)$.
From Proposition~\ref{ConjuStructu} we have that the conjugates
$R_{\rho}$ and $S_{\rho}:=\epsilon(\ol\rho.,\rho)R_{\rho}$ 
satisfy the relations
\[
R_{\rho}^{\ast}\ol\rho.(S_{\rho})=\un
 \quad\mr and.\quad 
S_{\rho}^{\ast}\rho (R_{\rho})=\un\;,
\]
and we can define in terms of these the following two vector space isomorphisms
(see e.g.~\cite{Longo97}):
\[
(\ol\rho.,\ol\rho.)\ni Y\rightarrow Y^{\ast}R_{\rho}\in (\iota,\ol\rho.\rho)
\quad\mr and.\quad 
(\iota,\ol\rho.\rho)\ni X\rightarrow S_{\rho}^{\ast}\rho(X)\in(\rho,\rho).
\]
Composing these isomorphisms in the case of an irreducible pair
$\rho,\ol \rho.$ we obtain
\[
(\ol\rho.,\ol\rho.)\ni\ol\rho.(Z)\rightarrow S_{\rho}^{\ast}\rho
(\ol\rho.(Z)^{\ast}R_{\rho})=Z^{\ast}\in(\rho,\rho)
\]
i.e.
\[
(\ol\rho.,\ol\rho.)\ni\ol\rho.(Z)\rightarrow Z^{\ast}\in(\rho,\rho)
\]
is a vector space isomorphism from
$(\ol\rho.,\ol\rho.)$
onto
$(\rho,\rho)$.
Therefore $\rho(\al Z.)=\al Z.$ and
$\rho\rest\al Z.$
is an automorphism of
$\al Z.$.
\end{beweis}

From the previous proposition it follows immediately:
\begin{cor}\label{FinProdIrr}
  \begin{itemize}
  \item[(i)] Let $\lambda=\rho_{1}\circ\rho_{2}\circ\dots\circ\rho_{n}$, 
             where $\rho_i\in\mr Irr.\al T.$, $i=1,\dots,n$. Then
             $\lambda(\al Z.)=\al Z.$ and 
             $\lambda\rest\al Z.\in\mr aut.\al Z.$.
  \item[(ii)] For any unitary $U\in\al A.$ we have 
              $(\mr ad. U\circ\rho)\rest\al Z.=\rho\rest\al Z.$,
              $\rho\in\mr Irr.\al T.$.
  \end{itemize}
\end{cor}

From the previous result we can now introduce the following automorphism
on $\al Z.$ which only depends on the class $D\in\wh{\al G.}$:
\begin{equation}\label{alphaD}
 \wh{\al G.}\ni D\to \alpha_\mt D.:=\rho_\mt D.\rest\al Z.\in\mr aut.\al Z.\,.
\end{equation}

We will introduce next
an equivalence relation\footnote{We would like to acknowledge an 
observation of P.A.~Zito concerning the equivalence relation that
served to simplify the presentation below.} 
in $\wwh {\al G.}.$
which, roughly speaking, relates elements $D,D'\in\wwh {\al G.}.$ if there 
is a ``chain of tensor products'' of elements in $\wwh {\al G.}.$
containing $D$ and $D'$. 
This equivalence relation appears naturally 
when considering the action of irreducible canonical endomorphisms 
on $\al Z.$ (see Theorem~\ref{EndoChain} and Remark~\ref{CompletePicture}
below).

To make the preceding notion precise recall the 
algebraic structure of $\wwh {\al G.}.$: 
denote by $\times$ the 
natural operation on subsets of $\wwh {\al G.}.$ associated with the 
decomposition of the tensor products of irreducible representations: 
for $D_1,D_2\in
\wwh {\al G.}.$ the set $D_1\times D_2$ contains the corresponding classes
that appear in the decomposition of 
$U_\mt D_1. \otimes U_\mt D_2.$,
where $U_\mt D_i.\in D_i$, $i=1,2$, is any representant of the 
corresponding class. That means that if 
\[
 U_\mt D_1. \otimes U_\mt D_2. = \mathop{\oplus}_{k=1}^s m_k\; U_\mt {D_k'}.
     \quad (m_k\equiv\mr multiplicities.)\,,
\]
then $D_1\times D_2=\{D_1',\dots,D_s'\}$.
For $\Gamma,\Gamma_1,\Gamma_2\subset\wg$ put
\[
 \Gamma_1\times \Gamma_2=\cup\{D_1\times D_2\mid D_i\in \Gamma_i,\,i=1,2\}
 \quad\mr and.\quad D\times\Gamma =\{D\}\times \Gamma\,.
\]
Further if $\ol D.\in\wg$ denotes the conjugate class to
$D\in\wg$ denote by $\ol \Gamma.=\{\ol D.\mid D\in \Gamma\}$.
Recall in particular that if $D\in D_0\times D_1$,
$D'\in D_0'\times D_1'$, then 
$D\times D'\subset D_0\times D_1\times D_0'\times D_1'$
or that $1\in\Gamma\times\ol \Gamma.$
(cf.~\cite[Definition~27.35]{bHewittII} for further details).

We can now make precise the previous idea:
\begin{defi}\label{DefChain}
$D,D'\in\wwh {\al G.}.$ are called equivalent, $D\approx D'$, 
if there exist $D_1,\ldots, D_n \in\wwh {\al G.}.$ such that
\[
D,D'\in D_1\times\dots\times D_n\,.
\]
\end{defi}

\begin{rem}\label{CommentChain}
  \begin{itemize}
  \item[(i)]
    It is easy to check that the relation $\approx$ is an equivalence
    relation in $\wwh {\al G.}.$. Indeed, reflexivity and symmetry 
    are trivial. To show transitivity let $D\approx D'$ as above
    and $D'\approx D''$, i.e.~$D',D''\in D_1'\times\dots\times D_k'$
    for some $D_1',\ldots, D_k'\in\wwh {\al G.}.$. We will show that
    $D,D''$ is contained in the larger chain
   \[
     \Gamma:=D_1\times\dots\times D_n\times\ol {D_1\times\dots\times D_n}.
             \times D_1'\times\ldots\times D_k'\,.
   \]
    Indeed, $D''\in\Gamma$, because 
    $1\in D_1\times\dots\times D_n\times\ol {D_1\times\dots\times D_n}.$ 
    and $D''\in D_1'\times\dots\times D_k'$. Further, since
    $\ol D'.\in\ol {D_1\times\dots\times D_n}.$ and 
    $D'\in D_1'\times\dots\times D_k'$ we also have that
    $1\in\ol {D_1\times\dots\times D_n}.\times D_1'\times\ldots\times D_k'$,
    and therefore $D\in\Gamma$ (cf.~\cite[Theorem~27.38]{bHewittII}).

    Denote by square brackets $[\cdot]$
    the corresponding chain equivalence classes and 
    by $\ot C.(\al G.)$ the set of these equivalence classes, i.e.
    \[ 
     \ot C.(\al G.):=\{[D] \mid D\in\wwh {\al G.}.\}\,.
    \]

\item[(ii)] 
Note that {\em any} pair $D,D'\in D_0\times D_1$ satisfies
by definition $D\approx D'$. Therefore for
$D_0,D_1\in\wwh {\al G.}.$ we have that $D_0\times D_1$
also specifies an element of $\ot C.(\al G.)$ and
we can simply put
\[
 [D_0\times D_1]:=[D]\,,
\]
where $D$ is any element in $D_0\times D_1$.

\item[(iii)] It is also possible to formulate the previous equivalence
relation entirely in terms of characters.
\end{itemize}
\end{rem}

We will define on $\ot C.(\al G.)$ a product $\bt$ 
(see Eq.~(\ref{ChainProduct}) below) so that 
$(\ot C.(\al G.),\bt)$ becomes an abelian group which for simplicity
we call {\bf chain group}.
Moreover, the chain group 
can be related to the character group of the
center $\al C.$ of $\al G.$. 
For this recall also the notion of 
conjugacy class of a representation
(cf.~\cite{bFuchs97}): 
let $D\in\wwh {\al G.}.$ and $U_\mt D.$ any representant in 
$D$. By Schur's Lemma we have 
\begin{equation}\label{RepCC}
 U_\mt D.\rest\al C.=\Upsilon_\mt D.\cdot\1\;,
\end{equation}
and it can be easily seen that $\Upsilon_\mt D.$ is a character
on the center $\al C.$ of $\al G.$
which only depends on $D$, i.e.~$\Upsilon_\mt D.\in\wwh {\al C.}.$. 

\begin{teo}\label{TeoChain}
Let $\al G.$ be a compact nonabelian group and denote 
its center by $\al C.$.
  \begin{itemize}
  \item[(i)] The set $\ot C.(\al G.)$ becomes an abelian group w.r.t.~the
             following multiplication: for $D_0,D_1\in\wwh {\al G.}.$ put
\begin{equation}\label{ChainProduct}
 [D_0]\;\bt\; [D_1]:= [D_0\times D_1]
\end{equation}
 (recall Remark~\ref{CommentChain}~(ii)).          
  \item[(ii)] The conjugacy classes $\Upsilon_\mt D.$ 
              (cf.~Eq.~(\ref{RepCC})) depend  
              on the chain equivalence class $[D]$.
              The chain group and the character group of the center of 
              $\al G.$ are isomorphic. The isomorphism is given by
\[
 \eta:\ot C.(\al G.)\to\wh{\al C.} 
 \quad\mr with.\quad \eta([D]):=\Upsilon_\mt [D].\;,
\]
where $\Upsilon_\mt [D].$ is the conjugacy class associated with 
$[D]\in\ot C.(\al G.)$. 
\end{itemize}
\end{teo}
\begin{beweis}
(i) Recall first that by Remark~\ref{CommentChain}~(ii) the r.h.s.~of 
(\ref{ChainProduct}) is well defined.
We need to verify next that the l.h.s.~of 
(\ref{ChainProduct}) is independent of the representants 
$D_0,D_1\in\wh{\al G.}$: suppose that $D_0\approx D_0'$
as well as $D_1\approx D_1'$ and we will show that 
any $D\in D_0\times D_1$ is related to any $D'\in D_0'\times D_1'$,
i.e.~$D\approx D'$, so that they specify the same equivalence class.
Let 
\[
 D_i,D_i'\in
 D_{\mt i1.}\times\dots\times D_{\mt ik_i.}\;,\quad i=0,1\,,
\]
for some $D_{\mt i1.},\dots, D_{\mt ik_i.}\in\wg$. Then
\[
 D\in D_0\times D_1\subset
      D_{\mt 01.}\times\dots\times D_{\mt 0k_0.}\times
      D_{\mt 11.}\times\dots\times D_{\mt 1k_1.}
                   \supset D_0'\times D_1'\ni D'\,.
\]
This shows that 
$D\approx D'$ and the product $\bt$ is well defined.

The neutral element of $\ot C.(\al G.)$ is given by the class generated by 
the trivial representation $1$. The inverse element to $[D]$ is given
by $[\ol D.]$, because $1\in D\times\ol D.$. Finally the commutativity
of the product is guaranteed by the equation $D_0\times D_1=D_1\times D_0$
(see~\cite[Theorem~27.38]{bHewittII} for details).

(ii) To show that the mapping $\eta$ is well defined note that
for the tensor product of {\em irreducible} representations $U_\mt D_i.$,
$i=1,\dots, n$, we have
\begin{equation}\label{TPC}
 U_\mt D_1.\otimes\dots\otimes U_\mt D_n.\rest\al C.
       =\Upsilon_\mt D_1.\cdot\ldots\cdot\Upsilon_\mt D_n. \1\;.
\end{equation}
Even more, any irreducible representation 
appearing in the decomposition of 
$U_\mt D_1.\otimes\dots\otimes U_\mt D_n.$ defines the same
conjugacy class $\Upsilon_\mt D_1.\cdot\ldots\cdot\Upsilon_\mt D_n.$. 
This shows that if $D\approx D'$, then $\Upsilon_\mt D.=\Upsilon_\mt D'.$ 
and that the $\eta$ is a group homomorphism from 
$\ot C.(\al G.)$ to $\wh{\al C.}$:
\[
 \Upsilon_{\mt [D_1].\,\bt\,\mt [D_2].}
           =\Upsilon_{\mt {[D_1\bt D_2]}.}
           =\Upsilon_\mt [D_1].\cdot\Upsilon_\mt [D_2].\,.
\]
Similarly we can consider the group homomorphism 
$\wh{\eta}\colon\ \al C. \to \wh{\ot C.(\al G.)}$ given by
\[
 \al C.\ni c\mapsto \wh{\eta}_c\,,\quad\mr where.\quad
 \wh{\eta}_c([D]):= \Upsilon_\mt [D].(c)\;,\quad [D]\in\ot C.(\al G.)\,.
\]
The homomorphism $\wh{\eta}$ is injective because if
$\wh{\eta}_{c_0}([D])=1$ for all $[D]\in\ot C.(\al G.)$,
then $U_\mt D.(c_0)=\1$ for all $D\in\wg$, hence $c_0=e$
(recall from Gelfand-Raikov's theorem \cite[Theorem~22.12]{bHewittI}
that the continuous irreducible unitary representations of $\al G.$ 
separate the points of the group).

Finally, the surjectivity of $\wh{\eta}$ is an application of Tannaka's 
duality theorem (cf.~Theorem~30.40 in \cite{bHewittII}). To sketch
the argument we need to recall the following facts. Let $\al T._\mrt TK.$
be a Tannaka-Krein category associated with the compact group $\al G.$.
We denote by $\{\al H._\mt D.\}_{\mt {D\!\in\!\dg}.}$ 
a complete set of irreducible objects. 
A representation $\ot r.$ of $\al T._\mrt TK.$ is an assignment
\[
 \mr Ob.\al T._\mrt TK.\ni
 \al H.\to \ot r.(\al H.)\in\al U.(\al H.)
 \quad (\mr Unitary~operators~on.~\al H.)\,,
\]
which is compatible with the direct sums, the tensor products, the conjugation
structure and the arrows of $\al T._\mrt TK.$ (see properties 
$T_1-T_6$ in 30.34 of \cite{bHewittII} for further details).
Now any character $\chi\in\wh{\ot C.(\al G.)}$ 
specifies the following assignment
\[
 \al H._\mt D.\to c_\chi(\al H._\mt D.)
         :=\chi([D])\cdot\1_{\al H._\mt D.}\in\al U.(\al H._\mt D.)
           \;,\quad D\in\dg\,.
\]
Taking into account direct sums, tensor products, the conjugation
structure and the arrows of $\al T._\mrt TK.$ we may extend 
$c_\chi(\cdot)$ to a representation of $\al T._\mrt TK.$.
By Tannaka's duality theorem $c_\chi(\cdot)$ specifies an 
element of $\al G.$. Even more, $c_\chi(\cdot)\in\al C.$, because
\[
 c_\chi(\al H._\mt D.) \cdot g(\al H._\mt D.)
     = g(\al H._\mt D.)\cdot c_\chi(\al H._\mt D.)\;,
       \quad D\in\dg
\]
implies
\[
 c_\chi(\al H.) \cdot g(\al H.)
     = g(\al H.)\cdot c_\chi(\al H.)\;,
       \quad \al H.\in\mr Ob.\al T._\mrt TK.\,.
\]
Therefore $\wh{\eta}_{c_\chi}=\chi$ and
we have shown that $\al C.$ and 
$\wh{\ot C.(\al G.)}$ are isomorphic. 
Pontryagin's duality theorem concludes the proof.
\end{beweis}

\begin{rem}
\begin{itemize}
\item[(i)]
The injectivity of the mapping $\eta$ in the previous theorem 
was stated as a conjecture in the first version of this paper 
(see Conjecture~5.10 in \cite{pLledo03v1}).
This conjecture was then proved by M.~M\"uger in
\cite{pMueger03}. For the sake of completeness we have
included a simple proof of this property. We also refer to
\cite{pMueger03} for further nice consequences of the isomorphism
$\ot C.(\al G.)\cong\wh{\al C.}$ in the context of fusion 
categories.

\item[(ii)]
We will leave for the end of this section the computation of 
chain groups associated with several finite and compact Lie groups.
In all the examples we will show explicitly
that the chain group $\ot C.(\al G.)$ is isomorphic 
to the character group  $\wwh {\al C.(\al G.)}.$ 
of the center of $\al G.$.
\end{itemize}
\end{rem}

We will now concentrate on the relation of the chain group $\ot C.(\al G.)$,
associated with the group $\al G.$ of a Hilbert C*-system $\HS$,
with the irreducible canonical endomorphisms restricted to $\al Z.$.
In particular recall the automorphisms on $\al Z.$ given in 
Eq.~(\ref{alphaD}) by
$\alpha_\mt D.:=\rho_\mt D.\rest\al Z.\in\mr aut.\al Z.$
which are associated with any class $D\in\wh{\al G.}$.

\begin{teo}\label{EndoChain}
\begin{itemize}
\item[(i)]
Let $D,D'\in\wh{\al G.}$ be equivalent, i.e.$~D\approx D'$. Then
$\alpha_\mt D.=\alpha_\mt D'.$ and we can
associate the automorphism $\alpha_\mt [D].\in\mr aut.\al Z.$  
with the chain group element $[D]\in\ot C.(\al G.)$.

\item[(ii)] There is a natural group homomorphism between the 
chain group and the automorphism group 
generated by the irreducible endomorphisms restricted to 
$\al Z.$ (cf.~Proposition~\ref{ProAutZS}):
\begin{equation}\label{HomoCA}
\ot C.(\al G.)\ni [D]\mapsto \alpha_\mt [D].\in\mr aut.\al Z.\,.
\end{equation}
\end{itemize}
\end{teo}
\begin{beweis}
(i)
First we show that if $\lambda=\rho_1\circ\dots\circ\rho_n$
is a finite product of irreducibles whose dimensions are larger
than one, then $\lambda\rest\al Z.=\alpha_\mt D.$ for all 
$D\in\wh{\al G.}$ 
appearing in the decomposition of $\lambda$: for this recall
from Proposition~\ref{DecompEnd} that $\lambda$ can be decomposed as
\begin{equation}\label{Ldecomp}
\lambda(\cdot)=\sum_{{\mt D.},j}W_{{\mt D.},j}\;\rho_{\mt D.}(\cdot)
               \;W_{{\mt D.},j}^{\ast}\,,
\end{equation}
where the $\rho_{\mt D.}$'s, $D\in\wwh {\al G.}.$, are irreducible and
$\{W_{{\mt D.},j}\}_{j}$ is an orthonormal basis of
$(\rho_{\mt D.},\lambda)$. 
Then by Corollary~\ref{FinProdIrr}~(i) and Proposition~\ref{ProAutZS} we 
have $\lambda(Z)\in\al Z.$ for all $Z\in\al Z.$ as well as 
\[
\lambda(Z)
   =\sum_{{\mt D.},j}W_{{\mt D.},j}\;\alpha_{\mt D.}(Z)\;W_{{\mt D.},j}^*
   =\sum_{{\mt D.}}\alpha_{\mt D.}(Z)\sum_{j}W_{{\mt D.},j}W_{{\mt D.},j}^*
   =\sum_{{\mt D.}}\alpha_{\mt D.}(Z)E_{\mt D.}\,,
\]
where
\begin{equation}\label{Erho}
E_{\mt D.}:=\sum_{j}W_{{\mt D.},j}W_{{\mt D.},j}^*
\end{equation}
is the so-called isotypical projection w.r.t.~$D\in\wwh {\al G.}.$.
Therefore
$\sum_{{\mt D.}}\alpha_{\mt D.}(Z)E_{\mt D.}=\lambda(Z)$
or
$\sum_{{\mt D.}}(\alpha_{\mt D.}(Z)-\lambda(Z))E_{\mt D.}=0$
and this implies
$(\alpha_\mt D.(Z)-\lambda(Z))E_{\mt D.}=0$
for all $D$ appearing in the decomposition of $\lambda$. 
Using (\ref{Erho}) and the orthogonality relations of the 
$W_{{\mt D.},j}$'s we obtain finally
$\alpha_{\mt D.}=\lambda\rest\al Z.$, hence all
irreducibles $\rho_\mt D.$ occurring in the decomposition of $\lambda$
coincide on $\al Z.$. 

Second, let $D\approx D'$, 
i.e.~$D,D'\in D_1\times\dots\times D_n$
for some $D_1,\ldots, D_n \in\wwh {\al G.}.$ 
(cf.~Definition~\ref{DefChain}).
For the endomorphisms this implies that 
$\rho_\mt {D}.$ and $\rho_\mt {D'}.$ appear in the decomposition of
$\rho_\mt {D_{\mt 1.}}.\circ\dots\circ\rho_\mt {D_{\mt n.}}.$
The arguments of the first part of the proof show that
\[
 \alpha_{\mt D.}=\alpha_{\mt D.'}\,.
\]

(ii) The homomorphism property follows immediately from the arguments
of the proof of part~(i) since for any $D_1,D_2\in\wh{\al G.} $ 
we have
\[
\alpha_\mt [D_1].\,\circ\,\alpha_\mt [D_2].
     =\alpha_{\,\mt [D_1\!\times\! D_2].}
     =\alpha_\mt {[D_1]\,\bt\, [D_2]}.
\]
(recall Theorem~\ref{TeoChain}~(i)).
\end{beweis}

\begin{rem}\label{CompletePicture}
Note that the chain group and in particular Theorem~\ref{EndoChain}~(i)
completes the picture of the action of the irreducible canonical 
endomorphisms on the center $\al Z.$ of the fixed point algebra $\al A.$
(recall also Eq.~(\ref{alphaD})). Indeed, we may summarize this 
action by means of the following diagram
\begin{eqnarray*}
  \wh{\al G.}&\to& \ot C.(\al G.)\;\; \to\;\; \mr aut.\,\al Z. \\
    D        &\mapsto& \; [D] \;\;\;\mapsto\;\; \alpha_\mt [D].
\end{eqnarray*}
\end{rem}

\begin{teo}\label{GeneralZMap}
Let $\lambda\in\ob\,\al T.$. Then its action on $\al Z.$ can be described
by means of the following formula
\[
 \lambda(Z)=\sum_{[\mt D.]} \alpha_\mt [D].(Z)\cdot
            \Big(\sum_{\mt {D'\in [D]}.} E_\mt D'. \Big)\;,
            \quad Z\in\al Z.\,,
\]
where $E_\mt D'.$ is the isotypical projection w.r.t.~$D'\in\wh{\al G.}$.
\end{teo}
\begin{beweis}
First note that for a general $\lambda\in\ob\,\al T.$ the equation
(\ref{Ldecomp}) is still valid (cf.~Proposition~\ref{DecompEnd}).
From this we have for $Z\in\al Z.$
\[
\lambda(Z)
   =\sum_{{\mt D.}}\alpha_{\mt D.}(Z)\sum_{j}W_{{\mt D.},j}W_{{\mt D.},j}^*
   =\sum_{{\mt D.}}\alpha_{\mt D.}(Z) E_\mt D.
   =\sum_{[\mt D.]} \alpha_\mt [D].(Z)\cdot
            \Big(\sum_{\mt {D'\in [D]}.} E_\mt D'. \Big)\;,
\]
where for the last equation we have used Theorem~\ref{EndoChain}~(i).
\end{beweis}

Next we show how a nontrivial chain group homomorphism (\ref{HomoCA})
acts as an obstruction to the existence of a symmetry associated with
the larger category $\al T.$.

Let $\al H.,\wt {\al H.}.\in\mr Irr.\al T._\al G.$ be irreducible 
algebraic Hilbert spaces and
$\ot H.=\al H.\al Z.,\wt {\ot H.}.=\wt {\al H.}.\al Z.
\in\mr Irr.\al M._\al G.$ the 
corresponding free $\al Z.$-modules. By Proposition~\ref{prop0}
we associate with them the irreducible endomorphisms 
$\rho,\wt \rho.\in\mr Irr.\al T.$ and denote the automorphisms
of their restriction to $\al Z.$ by
\[
\alpha:=\rho\rest\al Z.\quad\mr and.\quad\wt \alpha.:=\wt \rho.\rest\al Z.\,.
\]

If $\{ \Phi_i \}_i$ and $\{ \wt \Phi._j \}_i$
are orthonormal basis of $\al H.$ resp.~$\wt {\al H.}.$,
then
\begin{equation}\label{BaseChange}
 \Psi_i:=\sum_{i'}\Phi_{i'}\,Z_{i'i}\quad\mr and.\quad
 \wt \Psi._j:=\sum_{j'}\wt \Phi._{j'}\,\wt Z._{j'j}
 \quad Z_{i'i},\wt Z._{j'j}\in{\cal Z}\,,
\end{equation}
are arbitrary orthonormal basis of the corresponding modules
$\ot H.$ resp.~$\wt {\ot H.}.$, where  
\[
 {\got Z}:=(Z_{i'i})_{i',i}\quad\mr and.\quad 
 \wt {\got Z}.:=(\wt Z._{j'j})_{j',j}  
 \in\mr {Mat}.(\cal Z)\quad\mr satisfy.\quad
 {\got Z}^*{\got Z}={\got Z}{\got Z}^*=\1
    =\wt {\got Z}.^*\wt {\got Z}.=\wt {\got Z}.\wt {\got Z}.^*\,.
\]

\begin{pro}\label{ABasesZ}
With the previous notation define
\[
  \epsilon(\al H.,\wt {\al H.}.) 
     := \sum_{i,j} \wt \Phi._j\,\Phi_i\,\wt \Phi._j^*\,\Phi_i^*
 \quad\mr and.\quad
  \epsilon(\ot H.,\wt {\ot H.}.) 
     := \sum_{i,j} \wt \Psi._j\,\Psi_i\,\wt \Psi._j^*\,\Psi_i^* \,.
\]
Then we have 
\[
\epsilon(\al H.,\wt {\al H.}.) =
\epsilon(\ot H.,\wt {\ot H.}.)
\quad\mr iff.\quad
\wt \alpha.\left({\got Z}\right){\got Z}^*=\1
\quad\mr and.\quad
\wt {\got Z}.\;\alpha(\wt {\got Z}.^*)=\1\,,
\]
where $\wt \alpha.({\got Z}):=(\wt \alpha.(Z_{i'i}))_{i',i}
\in\mr {Mat}.(\cal Z)$ 
(similarly for $\alpha(\wt {\got Z}.^*)$).
\end{pro}
\begin{beweis}
Using Eq.~(\ref{BaseChange}) we have 
\begin{eqnarray*}
\epsilon(\ot H.,\wt {\ot H.}.)
  &=& \sum_{i,j} \wt \Psi._j\,\Psi_i\,\wt \Psi._j^*\,\Psi_i^*
      =\sum_{i,j,i',j',i'',j''}
        \wt \Phi._{j'}\,\wt Z._{j'j}\,\Phi_{i'}\,Z_{i'i}
       \wt Z._{j''j}^*\wt \Phi._{j''}^*\,Z_{i''i}^*\,\Phi_{i''}^*\\
  &=&   \sum_{i,j,i',j',i'',j''}
        \wt \Phi._{j'}\,\Phi_{i'}\, 
        \alpha^{-1}(\wt Z._{j'j})\,Z_{i'i}\,
       \wt Z._{j''j}^*\,\wt\alpha.^{-1}(Z_{i''i}^*)\,
        \wt \Phi._{j''}^*\,\Phi_{i''}^*\\ 
  &\mathop{=}\limits^{!}&
       \sum_{i,j} \wt \Phi._j\,\Phi_i\,\wt \Phi._j^*\,\Phi_i^* 
       =\epsilon(\al H.,\wt {\al H.}.)
\end{eqnarray*}
Multiplying the previous equation from the left with 
$\Phi_{i_0'}^* \wt \Phi._{j_0'}^*$
and from the right with $\Phi_{i_0} \wt \Phi._{j_0}$ 
we obtain
\[
 \sum_{i,j}\alpha^{-1}(\wt Z._{j'_0j})\,\wt Z._{j_0^{\phantom{*}}j}^*
           \,Z_{i'_0i}\,\wt\alpha.^{-1}(Z_{i^{\phantom{*}}_0i}^*)
   =\delta_{j_0'j_0}\,\delta_{i_0'i_0}\,.
\]
That means 
\[
 \sum_{j}\wt Z._{j'_0j}\,\alpha(\wt Z._{j_0^{\phantom{*}}j}^*)
       =\delta_{j_0'j_0}\quad\mr and.\quad
 \sum_{i} \wt\alpha.(Z_{i'_0i})\,Z_{i^{\phantom{*}}_0i}^*
   =\delta_{i_0'i_0}\,.
\]
and the proof is concluded.
\end{beweis}
\subsection{Examples of chain groups for some finite and compact Lie
            groups} \label{SomeChainGroups}

We will give next several examples of chain groups associated with 
nonabelian finite and compact Lie groups. We will also show 
that in all the examples considered the chain group is 
isomorphic to the character group of the center.
Note also that if the group
is abelian, then one can identify the chain group with the 
corresponding character group.

If $\al G.$ is the group we will denote its center by $\al C.(\al G.)$
and the corresponding chain group by $\ot C.(\al G.)$.

\paragraph{Compact Lie groups:}

We begin with the case $\al G.=\mr SU.\!(2)$. Denote by 
\[
 l\in\{0,\ms \frac12.,1,\ms \frac32.,\dots\}=\wwh {\mr SU.\!(2)}.
\]
the class specified by the usual
representation $T^{(l)}$ of $\mr SU.\!(2)$ on the space of complex polynomials
of degree $\leq 2l$ which has dimension $2l+1$. Then the 
decomposition theory for the tensor products $T^{(l)}\otimes T^{(l')}$
(cf.~\cite[Theorem~29.26]{bHewittII}) gives
\[
 l\times l'=\left\{\, |l-l'|,|l-l'|+1,\dots,l+l' \,\right\}\,,\quad 
             l,l'\in\{0,\ms \frac12.,1,\ms \frac32.,\dots\}\,.
\]
This decomposition structure implies that
\[
 l\approx l'\quad\mr iff. \quad 
         l,l'\;\; \mr are~both~integers~or~both~\mbox{half\,-\,integers}.\,.
\]
We can finally conclude that
\[
 \ot C.(\mr SU.\!(2))=\left\{[0],[\ms \frac12.]\right\}
                     \cong \Z_2\cong \wwh {\al C.(\mr SU.\!(}.2))\,.
\]

Using Brauer-Weyl theory one can directly establish 
for $\al G.=\mr SU.\!(N)$ the isomorphism 
between the corresponding chain group and the 
character group of the 
center.\footnote{Christoph Schweigert, private communication.}

Similarly one can proceed in other examples.
Using well-known results on the decomposition of the tensor
product of irreducible representations
(see e.g.~\cite[Section~29]{bHewittII})
we list the following further examples of chain groups.
\begin{itemize}
\item[(i)] If $\al G.=\mr U.\!(2)$ we have that its dual is given
           by the following labels:
 \[
   \wwh {\mr U.\!(2)}.=\left\{\, (m,l)\mid 
                              l\in\{0,\ms \frac12.,1,\ms \frac32.,\dots\}
                              \,,\;m\in\Z\quad\mr and.\quad
                              m+2l~\mr even.\right\}\,.
 \]
Then we compute
  \[
 \ot C.\left(\mr U.\!(2) \right)
            =\left\{\;[(m_+,0)]\;,\;[(m_-,\ms \frac12.)]\;\mid
                      m_+/m_-~\mr is~even/odd.\right\}
            \cong \Z\cong \wwh {\al C.(\mr U.\!(2}.))\,.
\]
\item[(ii)] If $\al G.=\mr O.\!(3)$ we have 
   \[
   \wwh {\mr O.\!(3)}.=\left\{\, (0,l)\,,\, (1,l)\mid 
                              l\in\{0,\ms \frac12.,1,\ms \frac32.,\dots\}
                       \right\}
 \]
and 
\[
\ot C.\left(\mr O.\!(3)\right)
            =\left\{\;[(0,0)]\;,\;[(1,0)]\;
              \right\}
            \cong \Z_2\cong \wwh {\al C.(\mr O.\!(3}.))\,.
\]
\item[(iii)] For $\al G.=\mr SO.\!(3)$ recall that
\[
 \wwh {\mr SO.\!(3)}.=\{0,1,2,\dots\}\,.
\]
In this case the corresponding center as well as the chain group are trivial:
\[
 \ot C.(\mr SO.\!(3))=\left\{[0]\right\}
                     \cong \wwh {\al C.(\mr SO.\!(}.3))\,.
\]
\end{itemize}

\paragraph{Finite groups:} We consider first the family of 
dihedral groups $\al G.=\D_{2m}$, $m\geq 2$. The group $\D_{2m}$
has order $2m$ and is generated by two elements $a,b$ that 
satisfy the relations
\[
 a^m=b^2=e\quad\mr and.\quad bab=a^{m-1}\,.
\] 
We consider first the case where $m=2l$, $l\in\N$, is even. Then the 
center of $\D_{2m}$ is $\al C.(\D_{2m})=\{e,a^l\}\cong\Z_2$ and
its dual is given by 
\[
 \wwh {\D_{2m}}.=\{1,\chi_1,\chi_2,\chi_3\}\cup\{D_k
                  \mid k=1,\dots,\ms \frac{m-2}{2}.\}\,,
\]
where $\{1,\chi_1,\chi_2,\chi_3\}$ are 1-dimensional representations
and $\{D_k\}_{\mt k=1.}^{\mt \frac{m-2}{2}.}$ 
are 2-dimensional representations. From the results 
concerning the decomposition of tensor products of 
irreducible representations stated in \cite[\S27.62~(d)]{bHewittII} 
we conclude:
\begin{equation}\label{D4m}
 \ot C.(\D_{2m})=\left\{\,[1]\,,\,[D_1]\,\right\}\cong\Z_2
                \cong\wwh {\al C.(\D_{2m})}.\,.
\end{equation}
 To check the details of the previous example $\D_{2m}$, $m=2l$, $l\in\N$,
 it is useful to distinguish further between the cases $l$ even or odd. 
 Indeed, if $l$ is even, then 
\[
 1\approx \chi_1\approx\chi_2\approx\chi_3\approx D_k\,,\; k~\mr even.
 \quad\mr and.\quad D_1\approx D_{k'}\,,\; k'~\mr odd\,.\,.
\]
If $l$ is odd, then the corresponding chain classes have
slightly different representants
\[
 1\approx \chi_1\approx D_k\,,\; k~\mr even.
 \quad\mr and.\quad \chi_2\approx\chi_3\approx D_1\approx D_{k'}
 \,,\; k'~\mr odd.\,.
\]

Similarly we can use the results in 
\cite[Section~27]{bHewittII} to compute the following 
family of examples:

\begin{itemize}
\item[(iv)] $\al G.=\D_{2m}$ with $m$ odd. Then the center is trivial, 
 $\al C.(\D_{2m})=\{e\}$, and
 \[
 \wwh {\D_{2m}}.=\{1,\chi_1\}\cup\{D_k
                  \mid k=1,\dots,\ms \frac{m-1}{2}.\}\,,
\]
where $\{1,\chi_1,\chi_2,\chi_3\}$ are 1-dimensional representations
and $\{D_k\}_{\mt k=1.}^{\mt \frac{m-1}{2}.}$ are 2-dimensional 
representations. As before we compute
\[
 \ot C.(\D_{2m})=\left\{\,[1]\,\right\}
                \cong\wwh {\al C.(\D_{2m})}.\,.
\]
\item[(v)] Let $\al G.=\Q_{4m}$ be the generalized quaternion groups
 which is a group of order $4m$ generated by two elements $a,b$ that 
 satisfy the relations $a^{2m}=b^4=e$, $b^2=a^m$ and $bab^{-1}=a^{2m-1}$.
 Its center is given by $\al C.(\Q_{4m})=\{e,a^m\}\cong\Z_2$ and
\[
 \wwh {\Q_{4m}}.=\{1,\chi_1,\chi_2,\chi_3\}\cup\{D_k
                  \mid k=1,\dots,m-1\}\,,
\]
where $\{1,\chi_1,\chi_2,\chi_3\}$ are 1-dimensional representations
and $\{D_k\}_{\mt k=1.}^{\mt m\!-\!1.}$ are 2-dimensional 
representations. Distinguishing again between
the cases $m$ even or odd we obtain using \cite[\S27.62~(e)]{bHewittII}
\begin{equation}\label{Q4m}
 \ot C.(\Q_{4m})=\left\{\,[1]\,,\,[D_1]\,\right\}\cong\Z_2
                \cong\wwh {\al C.(\Q_{4m})}.\,.
\end{equation}
\item[(vi)] We conclude this list of examples mentioning the cases of 
 the permutation groups $\P_3,\P_4$ and the the alternating group 
 $\A_4$ which have trivial center. It is straightforward to 
 verify that the corresponding chain groups are also trivial.
\end{itemize}

\begin{rem}
As stated in Corollary~3.2 of \cite{pMueger03} the isomorphism
between the chain group of $\al G.$ and the character group of its
center shows 
that the center of a compact group depends only
on the representation ring of $\al G.$. This is in fact 
explicitly verified for the
groups $\D_{8l}$ and $\Q_{8l}$,
$l\in\N$, which are particularly interesting for this question. 
Recall that these groups
are non-isomorphic but have isomorphic duals 
(cf.~\cite[\S27.62~(f)]{bHewittII}) and therefore isomorphic chain groups.
Therefore the centers of $\D_{8l}$ and $\Q_{8l}$ must also be isomorphic
(compare with the Eqs.~(\ref{D4m}) and (\ref{Q4m}) above).
\end{rem}

\section{A family of examples}
\label{RealInclu}

In this section we will give a
family of examples of pairs of
categories $\al T._\c<\al T.$,
where $\al T._\c$ is admissible
(recall Remark~\ref{Admi}~(iii)).

Let $\al A._\c$ be a unital C*-algebra with trivial center,
$\al Z.(\al A._\c)=\C\1$, and satisfying Property~B. 
Denote by $\ot Z.$
a unital abelian C*-algebra and define
\[
 \al A.:=\al A._\c\otimes\ot Z.\,,
\]
which is again a C*-algebra with unit $\1\otimes\1$ and center
$\al Z.=\al Z.(\al A.)=\1\otimes\ot Z.$. Let $\al T._\mrt DR.$ be
a DR-category (recall Definition~\ref{DRCat}) realized as 
endomorphisms of $\al A._\c$. The objects of $\al T._\mrt DR.$
are denoted by $\rho, \sigma$ etc.~and the 
corresponding arrows by $(\rho,\sigma)$.
Let $\al G.$ be the compact group associated with $\al T._\mrt DR.$
and denote by $\ot C.$ its corresponding chain group. We consider
also a fixed group homomorphism (recall Theorem~\ref{EndoChain})
\begin{equation}\label{HomoCA1}
\ot C.\ni [D]\mapsto \alpha_\mt [D].\in\mr aut.\al Z.\,.
\end{equation}

We can now start defining the C*-category $\al T.$ realized as endomorphism
of the larger algebra $\al A.$ with nontrivial center $\al Z.$.
To identify the new objects we
proceed in two steps: first we extend irreducible endomorphisms in
$\al T._\mrt DR.$ to endomorphisms of $\al A.$. Second,
we use the decomposition result in Proposition~\ref{DecompEnd}
to extend general objects in $\al T._\mrt DR.$ to endomorphisms
of $\al A.$. The extended endomorphisms of the larger algebra
$\al A.$ are interpreted as new objects of the category $\al T.$.

\begin{itemize}
\item[(a)] If $\rho$ is irreducible, $\rho\in\mr Irr.\al T._\mrt DR.$, we 
 define
\[
  \wt \rho.:=\rho\otimes\alpha_\mt [D].\in\mr end.\al A.\,,
\]
where $D\in\widehat{\al G.}$ is the corresponding class associated
with $\rho\in\mr Irr.\al T._\mrt DR.$ (cf.~Remark~\ref{AssocD}).

\item[(b)] Let $\tau\in\ob\,\al T._\mrt DR.$. According to
Proposition~\ref{DecompEnd}, the endomorphism
$\tau$ can be decomposed in terms 
of irreducible objects as
\[
 \tau(\cdot)=\sum_{\rho,\,l} W_{{\rho},\,l}\,\rho(\cdot) 
                \, W_{\rho,\,l}^*\,,
\]
where $\{ W_{\rho,\,l} \}_{l}$ denotes an \onb{} of
$(\rho,\,\tau)$ and $\rho\in\mr Irr.\al T._\mrt DR.$. 
We assign to $\tau$ the following endomorphism of $\al A.$
\begin{eqnarray}
 \wt \tau.(\cdot)
   &:=& \sum_{\rho,\,l} \left(W_{{\rho},\,l}\otimes\1\right)
        \,\wt \rho.\,(\cdot) \,\left( W_{\rho,\,l}\otimes\1 \right)^*
        \nonumber \\
   &=& \sum_{\rho,\,l} (W_{{\rho},\,l}\;\rho\,(\cdot)\, W_{\rho,\,l}^* )
       \otimes\alpha_\mt [D].(\cdot)\,,   \label{rewritetau}
\end{eqnarray}
where for the second equation we have used the previous item (a).

\item[(c)] The arrows in $\al T.$ are defined as usual
\[
 (\,\wt \rho.,\wt \tau.\,):=\{A\in\al A.\mid A \,\wt \rho.(X)=\wt \tau.(X)\,A
                          \;,\quad X\in\al A.\}\,.
\]
\end{itemize}

\begin{pro}\label{ExtIrr2}
Let $\al T.$ be the C*-category defined by means of (a),(b) and (c)
above. Then the objects $\wt \rho.,\wt \rho._1,\wt \rho._2$ 
defined in part (a) satisfy:
  \begin{itemize}
  \item[(i)] Irreducibility: $\wt \rho.\in\mr Irr.\al T.$, 
              i.e.~$(\,\wt \rho.,\wt \rho.\,)=\wt \rho.(\al Z.)=\al Z.$.
  \item[(ii)] Pairwise disjointness: $(\,\wt \rho._1,\wt \rho._2\,)=\{0\}$
              if $\wt \rho._1\not=\wt \rho._2$.
  \end{itemize}
\end{pro}
\begin{beweis}
For the proof it is convenient to apply Gelfand's theorem: any $Z\in\ot Z.$
can be identified with a continuous function over the compact space 
$\mr spec.(\ot Z.)$, $Z(\cdot)\in C(\mr spec.\ot Z.)$, and 
\[
 \al A.=\al A._\c\otimes\ot Z.\cong C(\mr spec.\ot Z.\to\al A._\c)\,.
\]
In particular, we will need below that 
any elementary tensor $A=A_0\otimes Z\in\al A._\c\otimes\ot Z.$
can be expressed as the function $\mr spec.\ot Z.\ni\mu\to A(\mu)=Z(\mu) A_0$.

(i) Let $\wt \rho.:=\rho\otimes\alpha_\mt [D].\in\ob\al T.$
with $D\in\widehat{\al G.}$ associated
with $\rho\in\mr Irr.\al T._\mrt DR.$
by means of the DR-Theorem 
(cf.~Remark~\ref{AssocD}). It is clear that 
$\wt \rho.(\al Z.)\subseteq (\,\wt \rho.,\wt \rho.\,)$, since
\[
 \wt \rho.\,(\al Z.)=\1\otimes\alpha_\mt [D].(\ot Z.)
                  =\1\otimes\ot Z. 
                  =\al Z.\subseteq (\,\wt \rho.,\wt \rho.\,)\,.
\]
For the converse inclusion let $A\in (\,\wt \rho.,\wt \rho.\,)$.
In particular, this implies
\[
 A\;\wt \rho.\,(X\otimes\1)=\wt \rho.\,(X\otimes\1)\,A\,,\quad X\in\al A._\c\,.
\]
The previous equation can be rewritten using Gelfand's theorem in 
terms of functions over $\mr spec.\ot Z.$ as
\[
 A(\mu) \; \rho\,(X)=\rho\,(X)\,A(\mu)\,,\quad X\in\al A._\c\;,\;
                   \mu\in\mr spec.\ot Z.\,.
\]
Since $\rho\in\mr Irr.\al T._\mrt DR.$, i.e.~$(\rho,\rho)=\C\1$, we conclude
that $A(\mu)=\lambda(\mu)\,\1$, where $\lambda$ is a continuous scalar
function on $\mr spec.\ot Z.$. Applying once more Gelfand's theorem we have
\[
 A\in\1\otimes\ot Z.=\1\otimes\alpha_\mt [D].(\ot Z.)=\wt \rho.\,(\al Z.)
\]
and we have shown that $\wt \rho.\in\mr Irr.\al T.$.

(ii) Let $\wt \rho._i=\rho_i\otimes\alpha_\mt [D_i].$,
$\rho_i\in\mr Irr.\al T._\mrt DR.$, $i=1,2$, and 
$\wt \rho._1\not=\wt \rho._2$. To show the disjointness relation
$(\,\wt \rho._1,\wt \rho._2\,)=\{0\}$ choose
$A\in(\,\wt \rho._1,\wt \rho._2\,)\subset\al A.=\al A._\c\otimes\ot Z.$.
The intertwiner element $A$ must satisfy in particular 
\[
 A\;(\rho_1\otimes \alpha_\mt [D_1].)\,(X\otimes\1)
    =(\rho_2\otimes \alpha_\mt [D_2].)\,(X\otimes\1)\,A
     \,,\quad X\in\al A._\c\,,
\]
which in terms of functions over $\mr spec.\ot Z.$ means
\[
 A(\mu) \; \rho_1\,(X)=\rho_2\,(X)\,A(\mu)\,,\quad X\in\al A._\c\;,\;
                   \mu\in\mr spec.\ot Z.\,.
\]
Since $\rho_1\not=\rho_2$ are irreducible in $\al T._\mrt DR.$ we
have $(\rho_1,\rho_2)=\{0\}$, hence $A(\mu)=0$, $\mu\in\mr spec.\ot Z.$,
and we conclude that $A=0$.
\end{beweis}

We can now state the main result of this section, namely the specification
of a family of  examples satisfying the hypothesis of Theorem~\ref{Teo2}.
Note that from the construction prescription in (a) and (b) above,
there is a bijective correspondence between $\ob\al T._\mrt DR.$
(which are realized as endomorphisms of $\al A._\c$) and 
$\ob\al T.$ (which are realized as endomorphisms of the larger 
algebra $\al A.=\al A._\c\otimes\ot Z.$).
\begin{teo}\label{AdCatTP}
Let $\al A._\c$ be a unital C*-algebra with trivial center
(and satisfying Property~B). Denote by $\ot Z.$
a unital abelian C*-algebra and define
\[
 \al A.:=\al A._\c\otimes\ot Z.\quad\mr so~that.\quad
  \al Z.:=\al Z.(\al A.)=\1\otimes\ot Z.\,.      
\]
If $\al T._\mrt DR.$ is a DR-category (cf.~Definition~\ref{DRCat})
realized as endomorphisms of $\al A._\c$, let $\al T.$ be the C*-category
specified by (a),(b),(c) above. Define the C*-category 
$\al T._\c$ as follows: $\ob\al T._\c:=\ob\al T.$ and
\[
 (\,\wt \sigma.,\wt \tau.\,)_\c:=(\sigma,\tau)\otimes\1\subset\al A.\,,
\]
where $\sigma,\tau$ are the objects in $\al T._\mrt DR.$ corresponding
to $\wt \sigma.,\wt \tau.\in\ob\al T.$. Then $\al T._\c$ is an admissible
subcategory of $\al T.$, i.e.~the arrows satisfy the 
equation 
\begin{equation}\label{GeneratingArrows2}
\quad(\wt \sigma.,\wt \tau.)=(\wt \sigma.,\wt \tau.)\s\c.\;\wt \sigma.(\al Z.)
             \,,
\end{equation}
(cf.~Theorem~\ref{Teo1}).
\end{teo}
\begin{beweis}
We need to show Eq.~(\ref{GeneratingArrows2}), that means that we can 
generate the ``larger'' arrows set of $\al T.$ with the ``smaller''
arrow set of $\al T._\c$ and the center $\al Z.$. From the decomposition
result in Proposition~\ref{DecompEnd} it is sufficient to prove the 
special case
\[
 (\wt \rho.,\wt \tau.)
         =(\wt \rho.,\wt \tau.)\s\c.\;\wt \rho.(\al Z.)
           \,,\quad \wt \rho.\in\mr Irr.\al T.\,,\,\wt \tau.\in\ob\al T.\,.
\]
First we show $(\wt \rho.,\wt \tau.)
\supseteq(\wt \rho.,\wt \tau.)\s\c.\;\wt \rho.(\al Z.)$. 
For this, take
an orthonormal basis $\{W_{\rho,l}\}_{l=1}^n \subset (\rho,\tau)$,
where $\rho\in\mr Irr.\al T._\mrt DR.$, $\tau\in\ob\al T._\mrt DR.$.
It is enough to show that for any $l=1,\dots,n$ and any $Z_0\in\ot Z.$
the following equation holds:
\[
 \Big((W_{\rho,\,l}\otimes\1)\,\,\wt \rho.\,(\1\otimes Z_0)\Big)
                                                \wt \rho.\,(X\otimes Z)
  =\wt \tau.\,(X\otimes Z)
              \Big((W_{\rho,\,l}\otimes\1)\,\,\wt \rho.\,(\1\otimes Z_0)\Big)
   \;,\quad X\in\al A._\c\,,\;Z\in\ot Z.\,.
\]
Using the definition of irreducible $\wt \rho.$ in part (a) above
we can rewrite the last equation as
\begin{equation}\label{IntwExtEq}
W_{\rho,\,l}\,\rho\,(X)\otimes \alpha_\mt [D].(Z_0Z)
   =\wt \tau.\,(X\otimes Z)\,(W_{\rho,l}\otimes \alpha_\mt [D].(Z_0))
      \;,\quad X\in\al A._\c\,,\;Z\in\ot Z.\,,
\end{equation}
where $D\in\wg$ corresponds to $\rho\in\mr Irr.\al T._\mrt DR.$ according
to DR-Theorem. We consider now the expression $\wt \tau.\,(X\otimes Z)$
separately and use Eq.~(\ref{rewritetau}) to obtain 
\[
  \wt \tau.\,(X\otimes Z)
    = \sum_{\rho',\,l'} (W_{{\rho'},\,l'}\;\rho'(X)\, W_{\rho',\,l'}^* )
       \otimes\alpha_\mt [D'].(Z)\,,
\]
where $\{W_{\rho',l'}\}_{l'} \subset (\rho',\tau)$
is an orthonormal basis 
and $D'\in\wg$ corresponds to $\rho'\in\mr Irr.\al T._\mrt DR.$.
Inserting this in the r.h.s.~of Eq.~(\ref{IntwExtEq})
and using the orthogonality relations
$W_{\rho',l'}^*W_{\rho,l}=\delta_{\rho\rho'}\delta_{ll'}$
we obtain
\begin{eqnarray*}
\wt \tau.\,(X\otimes Z)\,(W_{\rho,l}\otimes \alpha_\mt [D].(Z_0))
   &=& \sum_{\rho',\,l'} (W_{{\rho'},\,l'}\;\rho'\,(X)\, W_{\rho',\,l'}^* 
        W_{\rho,\,l})
       \otimes(\alpha_\mt [D'].(Z) \alpha_\mt [D].(Z_0)) \\
   &=& W_{\rho,\,l}\,\rho\,(X)\otimes \alpha_\mt [D].(Z_0Z)
\end{eqnarray*}
which coincides with the l.h.s.~of Eq.~(\ref{IntwExtEq}). This concludes
the proof of the inclusion $(\wt \rho.,\wt \tau.)
\supseteq(\wt \rho.,\wt \tau.)\s\c.\;\wt \rho.(\al Z.)$.

To show the reverse inclusion choose $A\in (\wt \rho.,\wt \tau.)$
so that
\begin{eqnarray*}
A\, \wt \rho.\,(X\otimes Z) 
     &=& \wt \tau.\,(X\otimes Z) \,A 
        \;,\quad X\in\al A._\c\,,\;Z\in\ot Z.\\
A\;( \rho\,(X)\otimes \alpha_\mt [D].(Z))   
     &=& \Big(\sum_{\rho',\,l'} (W_{{\rho'},\,l'}\;\rho'\,(X)\, 
          W_{\rho',\,l'}^* ) \otimes\alpha_\mt [D'].(Z)\Big)\,A\,,
\end{eqnarray*}
where $\{W_{\rho',l'}\}_{l'} \subset (\rho',\tau)$
is an orthonormal basis as before. Multiplying the previous equation 
with $W_{\rho,l}^*\otimes\1$ from the left we get
\[
((W_{\rho,l}^*\otimes\1) A)\;( \rho\,(X)\otimes \alpha_\mt [D].(Z))   
     =  (\rho\,(X) \otimes\alpha_\mt [D].(Z))\,
        ( ( W_{\rho,\,l}^* \otimes \1) \,A)
        \;,\quad X\in\al A._\c\,,\;Z\in\ot Z.\,,
\]
and this shows that 
\[
 (W_{\rho,l}^*\otimes\1) A\in(\wt \rho.,\wt \rho.)
      =\wt \rho.(\1\otimes\ot Z.)\,,
\]
where for the last equation we have used Proposition~\ref{ExtIrr2}~(i).
Therefore, for any $W_{\rho,l}$ there is a $Z_{\rho,l}\in\ot Z.$
such that 
$(W_{\rho,l}^*\otimes\1) A=\1\otimes \alpha_\mt [D].(Z_{\rho,l})$.
Multiplying this relation from the left
by $W_{\rho,l}\otimes\1$ and summing up 
w.r.t.~$l$ we obtain
\[
 \sum_l (E_\rho\otimes\1)\,A
    = \sum_l W_{\rho,l}\otimes \alpha_\mt [D].(Z_{\rho,l})\,,
\]
where $E_\rho:=\sum_l W_{\rho,l} \,W_{\rho,l}^*\in (\tau,\tau)$
is the isotypical projection w.r.t.~$\rho\in\mr Irr.\al T._\mrt DR.$.
To conclude the proof
recall the disjointness relation in 
Proposition~\ref{ExtIrr2}~(ii) which implies
\[
 (W_{\rho',l'}^*\otimes\1) A\in(\wt \rho.,\wt \rho'.)=\{0\}\quad
 \mr for~all.\quad \wt \rho.\not=\wt \rho'.\,,
\]
hence $(E_{\rho'}\otimes\1)\,A=0$ for all $\wt \rho.\not=\wt \rho'.$.
Now from the property $\sum_\rho E_\rho=\1$ of the isotypical
projections we obtain
\[
 A=\sum_l W_{\rho,l}\otimes \alpha_\mt [D].(Z_{\rho,l})
  =\sum_l (W_{\rho,l}\otimes\1) \,(\1\otimes\alpha_\mt [D].(Z_{\rho,l}))
   \in(\wt \rho.,\wt \tau.)\s\c.\;\wt \rho.(\al Z.)
\]
and the proof is concluded.
\end{beweis}

We can now apply Theorem~\ref{Teo2} to the pair of categories
$\al T._\c<\al T.$ constructed in this section to obtain the following
result:

\begin{pro}
Let $\al T._\c<\al T.$ be the pair of C*-categories constructed in
Theorem~\ref{AdCatTP}, where $\al T._\c$ is an admissible
subcategory of $\al T.$. Then there exists an essentially unique
minimal and regular Hilbert extension $\HS$ of $\al A.$.
\end{pro}

\begin{rem}
Note that construction of the inclusion of C*-categories
$\al T._\c<\al T.$ in Theorem~\ref{AdCatTP} depends crucially
on the choice of the chain group homomorphism in Eq.~(\ref{HomoCA1}).
Therefore different choices of this homomorphism will produce
different minimal and regular Hilbert extensions.
\end{rem}

\section{The case of a trivial chain group homomorphism}
\label{TrivialChainHom}

We will assume in this section that 
the chain group homomorphism given in Eq.~(\ref{HomoCA}) is trivial.
We will see that in this case the
analysis of minimal Hilbert C*-systems $\HS$ simplifies considerably.
In fact, this assumption implies that any irreducible endomorphism
acts trivially on the center $\al Z.$ and by Proposition~\ref{DecompEnd}
we finally obtain 
\begin{equation}\label{TrivialAct}
 \rho\rest\al Z.=\mr id.\rest\al Z. 
                 \quad\mr for~any.~\rho\in\mr Ob.\al T.\,.
\end{equation}
For example, the chain group homomorphism is trivial
if the chain group $\ot C.(\al G.)$ itself is trivial 
(see the examples in (iii),(iv) and (vi) of Subsection~\ref{SomeChainGroups}). 
This means that any $D\in\wh{\al G.}$
lies in the chain equivalence class of the trivial representation.

\begin{pro}\label{ZF}
Let $\HS$ be minimal Hilbert C*-system with fixed point algebra
$\al A.$ satisfying Property~B
and center $\al Z.=\al A.\cap\al A.'$. Then the
center of $\al F.$
coincides with
$\al Z.$, i.e.~$\al F.\cap\al F.'=\al Z.$.
\end{pro}
\begin{beweis}
Let
$\rho=\rho_{\al H.}$
be irreducible,
$\rho(A)=\sum_{j}\Phi_{j}A\Phi_{j}^{\ast}$
with an orthonormal basis
$\{\Phi_{j}\}_{j}$
of
$\al H.$.
Since
$\rho(Z)=Z$
we get
$\Phi_{j}Z=Z\Phi_{j}$
for all $j$. Further
$\al F.=C^{\ast}(\al A.,\{\al H.\})$,
where
$\al H.$
runs through all irreducible Hilbert spaces. This implies
$ZF=FZ$
for all
$F\in\al F.$. Therefore
$\al Z.\subseteq\al F.'\cap\al F.\subseteq\al A.'\cap\al F.=\al Z.$,
hence
$\al F.'\cap\al F.=\al Z.$
follows.
\end{beweis}

Next we show that in the case of a trivial chain group 
homomorphism one can still
associate a symmetry with the larger category $\al T.$.
For this purpose, recall from Proposition~\ref{prop0} that
to any $\rho\in\mr Ob.\al T.$ there exists a unique
free $\al Z.$-bimodule 
$\ot H._\rho\in\mr Ob.\al M._\al G.$.

\begin{cor}\label{EndZMod}
Each canonical endomorphism of a minimal Hilbert C*-system
$\{\al F.,\al G.\}$
is a $\al Z.$-module endomorphism. The symmetry
$\epsilon(\rho,\sigma)$, $\rho,\sigma\in\mr Ob.\al T.$, defined 
for $\{\Phi_i\}_i$,$\{\Psi_j\}_j$ orthonormal basis in 
$\ot H._\rho$, $\ot H._\sigma$ by
\[
 \epsilon(\rho,\sigma):=\sum_{i,j} \Psi_j\Phi_i\Psi_j^*\Phi_i^*
                        \in (\rho\sigma,\sigma\rho)\,,
\] 
satisfy the corresponding properties of Proposition~\ref{PermuStructu}.
\end{cor}
\begin{beweis}
From the result $\al F.'\cap\al F.=\al Z.$ and the definition of
a canonical endomorphisms we get immediately
that $\rho(AZ)=\rho(A)Z=Z\rho(A)$ for
any $\rho\in\mr Ob.\al T.$ and any $Z\in\al Z.$, 
hence the objects in $\al T.$ are 
$\al Z.$-module endomorphism.  
In particular, this also implies 
$\alpha:=\rho\rest\al Z.=\mr id.\rest\al Z.$ for all 
$\rho\in\mr Ob.\al T.$. Therefore by Proposition~\ref{ABasesZ}
we get that the definition of $\epsilon(\cdot,\cdot)$ is
independent of the module basis chosen in the 
Hilbert $\al Z.$-module 
assigned to the objects of $\al T.$.
The additional properties of $\epsilon(\cdot,\cdot)$ are
then verified easily (cf.~Proposition~\ref{PermuStructu}).
\end{beweis}

\begin{rem}\label{NotSym}
\begin{itemize}
\item[(i)]
 As already mentioned in \cite[Remark~6.4]{Lledo97b} it is not 
 possible in general to associate a symmetry 
 $\epsilon(\cdot,\cdot)$ with the larger 
 category of canonical endomorphisms $\al T.$. The reason is 
 that, in general, the formula in the previous corollary is not 
 independent of the module basis chosen
 (recall Proposition~\ref{ABasesZ}). 
 Therefore the existence 
 of a symmetry in the present context suggests that the 
 nontriviality of the chain group 
 homomorphism given in Eq.~(\ref{HomoCA})
 is an obstruction to the 
 existence of a well defined $\epsilon$ within the category
 $\al T.$.
\item[(ii)]
  The present section is also related to the notion of extention 
  of C*-categories $\al C.$ by abelian C*-algebras $C(\Gamma)$
  studied in \cite{Vasselli03}. In this reference it is shown
  that the DR-algebra associated with
  an object of the extension category $\al C.^{\Gamma}$
  is a continuous field of DR-algebras corresponding to
  the initial category $\al C.$.
  (For the construction of DR-algebras associated with
   suitable C*-categories see \cite{Doplicher87}.)
\end{itemize}
\end{rem}

The previous corollary means that, 
in the present situation, the category $\al T.$
of all canonical endomorphisms is ``almost'' a DR-category
(cf.~Definition~\ref{DRCat}): in fact, there
is a permutation and a conjugation structure. The only difference 
is that we have $(\iota,\iota)=\al Z.\supset\C\1$.

The next theorem shows that, using the central decomposition
w.r.t.~the common center $\al Z.$ (cf.~Proposition~\ref{ZF}),
a minimal Hilbert C*-system $\HS$
satisfying $\rho\rest\al Z.=\mr id.\rest\al Z.$,
$\rho\in\ob\al T.$,
can be considered as a direct integral
of Hilbert C*-systems $\{\al F.(\lambda),\al G.\}_\lambda$,
$\lambda\in\mr spec.\al Z.$,
with trivial relative commutant and
a fiber-independent compact group. Moreover, the Hilbert C*-systems 
$\{\al F.(\lambda),\al G.\}$ associated with the base points $\lambda$
are Hilbert extensions of their fixed point algebras $\al A.(\lambda)$
which carry DR-categories given by the fibre decomposition of the 
category $\al T.$ of $\{\al F.,\al G.\}$.

We assume for the rest of this section that the C*-algebras
$\al A.$ and $\al F.$ are {\em separable} and 
{\em faithfully represented}
in some Hilbert space $\ot h.$. We put $\al Z.\cong C(\Gamma)$,
$\Gamma:=\mr spec.\al Z.$ (cf.~Section~\ref{AppenDecom}).

\begin{teo}\label{Zerleg}
Let
$\{\al F.,\al G.\}$ 
be a minimal Hilbert C*-system 
with fixed point algebra $\al A.$ satisfying Property~B.
$\al Z.$ is 
the common center of
$\al A.$ and
$\al F.$ (cf.~Proposition~\ref{ZF}).
The fiber C*-algebras corresponding to the 
central decomposition w.r.t.~$\al Z.$ (cf.~Section~\ref{AppenDecom})
are denoted by $\al A.(\lambda)$, 
$\al F.(\lambda)$, $\lambda\in\Gamma$. 
Then, there
is an exceptional Borel set $\Gamma_0\subset\Gamma$, with $\mu(\Gamma_0)=0$,
such that for all $\lambda\in\Gamma\setminus\Gamma_0$,
\begin{itemize}
\item[(i)]
$\al A.(\lambda)\subset\al F.(\lambda)$ 
and
$\al A.(\lambda)$
satisfies Property~B (cf.~Subsection~\ref{GCanEnd}).
\item[(ii)]
Let $\al H.\subset\al F.$ be a $\al G.$-invariant algebraic Hilbert 
space with support $\1$. Then the
fiber spaces $\al H.(\lambda)\subset\al F.(\lambda)$ are again
$\al G.$-invariant algebraic Hilbert spaces
satisfying $\mr supp.\al H.(\lambda)=\1$.
If $\ot H.=\al H.\al Z.$ is the free $\al Z.$-bimodule generated
by $\al H.$, then $\ot H.(\lambda)=\al H.(\lambda)$.
\item[(iii)]
$\{\al F.(\lambda),\al G.\}$
is a Hilbert C*-system.
Let $\rho,\sigma\in\mr Ob.\al T.$ be canonical endomorphisms and denote
by $\rho_\lambda,\sigma_\lambda$ their fibre decomposition
(cf.~Proposition~\ref{FibreRho}). Then $\rho_\lambda,\sigma_\lambda$
are canonical endomorphisms associated with
$\{\al F.(\lambda),\al G.\}$
and the intertwiner space
$(\rho_{\lambda},\sigma_{\lambda})$
is given by
$(\rho_{\lambda},\sigma_{\lambda})=(\rho,\sigma)_{\lambda}$.
\item[(iv)]
The category 
$\al T._{\lambda}$
of canonical endomorphisms associated with
$\{\al F.(\lambda),\al G.\}$
is a DR-category, that means in particular
$(\iota_{\lambda},\iota_{\lambda})=\C\un_{\lambda}$
and the Hilbert C*-system
$\{\al F.(\lambda),\al G.\}$
has a trivial relative commutant. 
\end{itemize}
\end{teo}
\begin{beweis}
(i) The inclusion $\al A.(\lambda)\subset\al F.(\lambda)$ 
is obvious from Theorem~\ref{C*Zdecom} in the appendix.
The C*-algebra $\al A.$ satisfies Property~B,
if there exist isometries $V,W\in\al A.$ satisfying
$VV^*+WW^*=\1$. By Theorem~\ref{C*Zdecom} it follows
that the representation of $V,W$ on the fibre spaces
satisfy analogous properties.

(ii) Let $\{\Phi_i\}_i$ be an orthonormal basis of the
$\al G.$-invariant algebraic Hilbert $\al H.\in\mr Ob.\al T._\al G.$.
It transforms according to a unitary representation 
$U$ of $\al G.$. By
Theorem~\ref{C*Zdecom} we have that, for 
$\lambda\in\Gamma\setminus\Gamma_0$,
$\{\Phi_i(\lambda)\}_i$ is an 
orthonormal basis of $\al H.(\lambda)$ transforming 
according to the {\it same} representation $U$ 
(hence $\al H.(\lambda)$ is $\al G.$-invariant) 
and $\mr supp.\al H.(\lambda)=\1$. Finally, let
$\ot H.=\al H.\al Z.\subset\al F.$ be the free 
$\al Z.$-bimodule generated by $\al H.$. Any $H\in\ot H.$ can 
be written as $\sum_i \Phi_i\,Z_i$ for some $Z_i\in\al Z.$,
hence its fibre component becomes 
\[
 H(\lambda)=\sum_i \Phi_i(\lambda)\,Z_i(\lambda)\in\al H.(\lambda).
\]
This shows $\ot H.(\lambda)=\al H.(\lambda)$ for 
all $\lambda\in\Gamma\setminus\Gamma_0$.

(iii) The first part follows already from (ii). Eq.~(\ref{TrivialAct}) implies
that all canonical endomorphisms of $\{\al F.,\al G.\}$ are 
$\al Z.$-module endomorphisms, hence from
Proposition~\ref{FibreRho} we have
\[
 \rho_\lambda(A(\lambda))=
       \sum_i \Phi_i(\lambda)\,A(\lambda)\,\Phi_i(\lambda)^*\,,
\]
i.e.~$\rho_\lambda$ is canonical w.r.t.~$\{\al F.(\lambda),\al G.\}$.
Finally, it is straightforward to show that
$A\in(\rho,\sigma)$
iff
$A(\lambda)\in(\rho_{\lambda},\sigma_{\lambda})$
for all $\lambda\in\Gamma\setminus\Gamma_0$.

(iv) The permutation and conjugation structures of $\al T._\al G.$
(recall~Subsection~\ref{Conjperm}) can be ``disintegrated'' 
and the fibre components define the 
permutation and conjugation structures of $\al T._\lambda$.
Finally,
$(\iota,\iota)=\al Z.$
implies
$(\iota_{\lambda},\iota_{\lambda})=\C\un_{\lambda}$
and
$\al A.'\cap\al F.=\al Z.$
gives
$\al A.(\lambda)'\cap\al F.(\lambda)=\C\un_{\lambda}$
for all $\lambda\in\Gamma\setminus\Gamma_0$.
\end{beweis}

\begin{rem}\label{UniqueDisconn}
Note the item (iii) in the previous theorem implies
that the compact group $\al G.$ is unique for all 
$\al F.(\lambda)$, $\lambda\in\Gamma\setminus\Gamma_0$,
even if $\mr spec.\al Z.$ is disconnected.
\end{rem}

Furthermore we have the following inverse theorem.

\begin{teo}\label{Constrfromfibre}
Let
$\al A.$
be a unital C*-algebra with center $\al Z.$ and satisfying 
Property~B. Suppose that $\Gamma:=\mr spec.\al Z.$ is 
connected and let
$\al T.$ be
a tensor C*-category realized as unital endomorphisms of
$\al A.$ and
equipped with the following properties:
\begin{itemize}
\item[(i)]
All
$\rho\in\ob\,\al T.$
are $\al Z.$-module endomorphisms.
\item[(ii)]
$\al T.$
is closed w.r.t.~direct sums and subobjects.
\item[(iii)]
$\al T.$
is equipped with a permutation  and a
conjugation structure (cf.~Propositions~\ref{PermuStructu}
and \ref{ConjuStructu}).
\end{itemize}
Then there is a minimal Hilbert extension
$\{\al F.,\al G.\}$
of
$\al A.$ with $\al F.'\cap\al F.=\al Z.$ and
such that
$\al T.$
is isomorphic to the category of all canonical endomorphisms of
$\{\al F.,\al G.\}.$
\end{teo}
\begin{beweis}
Let
$\al A.=\int_{\Gamma}\al A.(\lambda)\mu(d\lambda)$
be the central decomposition of
$\al A.$ over the direct integral
$L^2(\Gamma,\mu,\ot f._\lambda)$ 
(cf.~Theorem~\ref{C*Zdecom}).
Denote by $\al T.(\lambda)$, 
$\lambda\in\Gamma\setminus\Gamma_0$,
the C*-category associated with $\al A.(\lambda)$
given by the fibre decomposition of $\al T.$ and
whose objects $\rho_\lambda$ are 
realized as endomorphisms of $\al A.(\lambda)$
(cf.~Proposition~\ref{FibreRho} and recall that
$\Gamma_0$ is the corresponding exceptional set).
First we show that
$\al T.(\lambda)$
is a DR-category (recall Definition~\ref{DRCat}): indeed, 
$(\iota_\lambda,\iota_\lambda)=\C\1_\lambda$ follows from 
the fact that that the C*-algebras $\al A.(\lambda)$ have 
trivial center (cf.~Theorem~\ref{C*Zdecom}). Note also that
$\rho\in\mr Irr.\al T.$ iff 
$(\rho_{\lambda},\rho_{\lambda})=\C\un_{\lambda}$,
$\lambda\in\Gamma\setminus\Gamma_0$,
i.e.~$\rho_{\lambda}\in\mr Irr.\al T.(\lambda)$.
Closure under direct sums and
subobjects follows from (ii). Recall that from (i), (iii) and 
Corollary~\ref{EndZMod} we have a well defined permutation
structure on $\al T.$ which can be carried over to $\al T.(\lambda)$.
Similarly we can ``disintegrate'' the conjugation structure, 
e.g.~we have for the conjugate
\[
 (\iota,\ol\rho.\rho)\ni R_\rho
  =\int_{\Gamma} R_{\rho_\lambda}\;\mu(d\lambda)  
\]
and $R_{\rho_\lambda}$ satisfies the corresponding properties
on the fibre. This shows that $\al T.(\lambda)$
is a DR-category for all
$\lambda\in\Gamma\setminus\Gamma_0$.

Now, by the DR-theory we have on the one hand that
\[
  R_{\rho_\lambda}^*R_{\rho_\lambda}=d_{\rho_{\lambda}}\un_{\lambda}\,,
\]
where $d_{\rho_{\lambda}}\in\N$. On the other hand 
$R_{\rho}^*R_{\rho}\in\al Z.\cong C(\Gamma)$ and therefore
$\lambda\to d_{\rho_{\lambda}}$
is continuous on $\Gamma$.
That means that
$d_{\rho_{\lambda}}=d_{\rho}\in\N$
is constant over $\Gamma$. We use this result to analyze 
the Hilbert extension
$\{\al F.(\lambda),\al G.(\lambda)\}$ of 
$\al A.(\lambda)$ which satisfies
$\al A.(\lambda)'\cap\al F.(\lambda)=\C\un_{\lambda}$
(cf.~\cite{Baumgaertel97})).
The existence of this fibre Hilbert C*-system for 
all $\lambda\in\Gamma\setminus\Gamma_0$
is guaranteed by the DR-Theorem. 
For any $\rho_{\lambda}\in\mr Irr.\al T.(\lambda)$ we consider a 
$d_\rho$-dimensional
algebraic Hilbert space $\al H._{\rho_\lambda}\equiv\al H._\rho$,
which is constant over $\Gamma$ and 
generated by the an orthonormal 
basis $\{\Phi_{\rho,i}\}_{i=1}^{d_\rho}$, i.e.
\begin{equation}\label{GeneAHS}
 \al H._\rho=\mr span.\{\Phi_{\rho,i}\mid i=1,\dots, d_\rho \}\,.
\end{equation}
We may assume that the isometries $\Phi_{\rho,i}$ are represented 
on a fixed Hilbert space $\ot h._0$. For any arbitrary
$\tau_\lambda\in\ob\,\al T.(\lambda)$ we associate the algebraic 
Hilbert space
\[
\al H._{\tau}(\lambda)
   :=\mathop{\oplus}\limits_{\rho_\lambda\in\mr Irr.\al T.(\lambda)}
     (\rho_{\lambda},\tau_{\lambda})\;\al H._{\rho}\,.
\]
We have
\[
 \al H._{\tau}(\lambda)\subset\al F._{\mrt fin.}(\lambda)
 :=\left\{\sum\limits_{\mt {\rho_\lambda\!\in\!\mr Irr.\al T.(\lambda)}.}
   A_{\rho_\lambda,j}(\lambda)\,\Phi_{\rho,j}\right\}
   \subset\al F.(\lambda)
 :=\mr clo._{\|\cdot\|_\lambda}\al F._{\mrt fin.}(\lambda)
\]
and the elements of the algebraic Hilbert space are bounded operators
on $\ot f._\lambda\otimes\ot h._0$. Put
\[
\al F.:=\int_{\Gamma}\al F.(\lambda)\mu(d\lambda)
    \subset\al L.\left(L^2(\Gamma,\mu,\ot f._\lambda\otimes\ot h._0)\right)
\]
with C*-norm given by 
$\Vert\cdot\Vert:=\hbox{ess sup}_{\lambda}\,\Vert\cdot\Vert_{\lambda}$.

Finally we have to define a compact group action on $\al F.$. 
For this, recall that on 
each fibre over
$\lambda\in\Gamma\setminus\Gamma_0$,
the compact group $\al G.(\lambda)$ acts
as follows: each $A(\lambda)\in\al A.(\lambda)$ is left pointwise
invariant under the group action. Moreover, any 
$\al H._\rho$, $\rho\in\mr Irr.\al T.$, carries an irreducible 
representation of $\al G.(\lambda)$ which does not depend 
on $\lambda$ (cf.~Eq.~(\ref{GeneAHS})). Therefore 
the action $\al G.(\lambda)$ is independent of $\lambda$, hence
we put $\al G.(\lambda)\equiv\al G.$, where $\al G.$ is compact.
Since $\{\al F.(\lambda),\al G.\}$ is a Hilbert C*-system it follows
immediately that $\{\al F.,\al G.\}$ is also a Hilbert C*-system. 
We still need to show that it is also minimal.
For this take for any $D\in\al G.$ a canonical endomorphism
$\rho_{\mt D.,\lambda}\in\mr Irr.\al T.(\lambda)$, 
$\lambda\in\Gamma\setminus\Gamma_0$.
Since $\rho_{\mt D.,\lambda}$ is disjoint from $\iota_\lambda$ for
any $D\not=\iota$, i.e.~$(\rho_{\mt D.,\lambda},\iota_\lambda)=0$,
we have that the corresponding integrated 
endomorphism satisfies the same property 
$(\rho_{\mt D.},\iota)=0$ for any $D\not=\iota$. 
From \cite[Lemma~10.1.8]{bBaumgaertel92} we have that
$\al A.'\cap\al F.=\al Z.$ and the proof is completed.
\end{beweis}

\begin{cor}\label{Mintoreg}
Each minimal Hilbert C*-system satisfying $\rho\rest\al Z.=\mr id.\rest\al Z.$,
$\rho\in\ob\al T.$, with $\Gamma:=\mr spec.\al Z.$ 
connected, is regular.
\end{cor}
\begin{beweis}
First recall from the proof of the previous theorem
that the fiber Hilbert C*-systems have trivial relative
commutant. This means that for all
$\lambda\in\Gamma\setminus\Gamma_0$
there is a one-to-one correspondence between the fiber
endomorphisms and the generating Hilbert spaces,
\[
\rho_{\lambda}\leftrightarrow \al H._{\rho_{\lambda}}\subset 
\al F.(\lambda)
\]
such that
\[
\rho_{\lambda}\circ\sigma_{\lambda}\leftrightarrow 
\al H._{\rho_{\lambda}}\cdot\al H._{\sigma_{\lambda}}.
\]
If one chooses for all
$\lambda\in\Gamma\setminus\Gamma_0$
a fixed Hilbert space
$\al H._{\rho_{\lambda}}$,
then these fiber spaces define a Hilbert space
$\al H._{\rho}\subset\al F.$
such that
\[
\rho\circ\sigma\leftrightarrow \al H._{\rho}\cdot\al H._{\sigma}.
\]
Therefore
$\{\al F.,\al G.\}$
is regular.
\end{beweis}

\begin{rem}
We will extend in this remark the inverse result stated 
in Theorem~\ref{Constrfromfibre}
to the following situation: let $\{\al A.,\al T.\}$ satisfy the 
assumptions of Theorem~\ref{Constrfromfibre} except that now
$\Gamma$ is a disjoint union
of (in general infinite) connected components $\Gamma_a$, $a\in \A$, i.e.
\[
 \Gamma=\mathop{\cup}\limits^{\cdot} \Gamma_a \,.
\] 
With each $\Gamma_a$, $a\in \A$, we can associate a central orthoprojection 
$P_a$, which is defined by means of the following continuous
function over $\Gamma$:
\[
   P_a(\lambda):=\left\{
                 \begin{array}{l}
                  1\quad\mr if.\lambda\in\Gamma_a \\
                  0\quad\mr otherwise. \,.
                 \end{array}
                 \right. 
\]
The projections in $\{P_a\}_{a\in \A}$ are mutually disjoint 
and satisfy $\sum_a P_a=\1$ (strong operator convergence;
to define the previous sum in the infinite case
consider a net of projections indexed by the class of all
finite subsets of $\A$ partially ordered by inclusion
$\subseteq$, cf.~\cite[Sections~2.5 and 2.6]{bKadisonI}). 
The Hilbert space $\ot h.$, on which the algebra $\al A.$ is 
represented, decomposes as
$\ot h.=\oplus_a\ot h._a$ and
$(x_a)_{a\in\A}\in\oplus_{a\in \A}\ot h._a$
if $\sum_a \|x_a\|^2<\infty$.
Therefore we can decompose $\al A.$ as a direct sum
\[
 \al A.=\oplus_{a\in \A}\al A._a\quad\mr and~in~particular.\quad
        \al Z.=\oplus_{a\in \A}\al Z._a 
\]
where $\al A._a:=\al A.\, P_a$ has center~$\al Z._a:=\al Z.\, P_a$.
Recall that $(A_a)_{a\in\A}\in\oplus_{a\in \A}\al A._a$ if
$\mr sup.\{\|A_a\|\mid a\in\A\}<\infty$.

From property (i) we can consistently define a family of 
$\al Z._a$-module endomorphisms
$\al T._a:=\{\rho_a\}\subset\mr end.\al A._a$ by means of
\[
 \rho_a(AP_a):=\rho(AP_a)=\rho(A)P_a\in\al A._a\,,\quad A\in\al A.\,.
\]
Moreover, since $\al A.$ satisfies Property~B, then $\al A._a$ also satisfies
this property on $\ot h._a$, 
i.e.~$\al A._a$ contains are isometries $V_a,W_a$
satisfying $V_aV_a^*+W_aW_a^*=P_a$. Similarly we can adapt the assumptions
(i)-(iii) to the pair $\{\al A._a,\al T._a\}$, $a\in \A$.
By the proof of Theorem~\ref{Constrfromfibre}
we can construct Hilbert C*-systems
$\{\al F._a,\al G._a\}$, with $\al G._a$ compact and satisfying 
\[
 \al A._a'\cap \al F._a=\al Z._a
    \quad\mr as~well~as.\quad
 \al F._a'\cap \al F._a=\al Z._a\;,\quad a\in \A\,.
\]
Now, in order to be able to built up from these systems a minimal Hilbert 
C*-system with a compact group we need to make the following additional
assumption (recall Remark~\ref{UniqueDisconn}):
\begin{itemize}
\item[]{\bf Assumption:}
  The compact groups $\al G._a$ are mutually isomorphic, i.e.
\begin{equation}\label{AssunBuilt}
  \al G._a\cong\al G.\,,\quad\mr for~some~compact~group. \al G.\,.
\end{equation}
\end{itemize}
Under this assumption put finally
$\al F.:=\oplus_a\al F._a$, where $\al G.$ acts on each component.
The Hilbert C*-system $\{\al F.,\al G.\}$ satisfies
\[
 \al A.'\cap \al F.=\oplus_a\al Z._a=\al Z.\,,
\]
hence it is minimal. Finally, note that Corollary~\ref{Mintoreg} can
be also adapted to the present more general situation satisfying
the assumption (\ref{AssunBuilt}).
\end{rem}

\begin{rem}
The reason why we need to make the assumption (\ref{AssunBuilt})
is that we want to reconstruct Hilbert C*-systems $\{\al F.,\al G.\}$
of the type studied in Theorem~\ref{Zerleg}. A similar situation that
considers more general groups, where e.g.~the mapping 
$\Gamma\ni\lambda\to\al G._\lambda$ is not constant, is studied 
in \cite[Section~3, Example~3.1]{pVasselli03b}.
\end{rem}

\section{Conclusions}\label{Conclu}

In the present paper we have described  a generalization of the
DR-duality theory of compact groups, to the case where the underlying 
C*-algebra $\al A.$ has a nontrivial center. The abstract characterization
of minimal and regular Hilbert C*-systems with a compact group
$\al G.$ 
is now given by the inclusion of C*-categories 
$\al T._\c<\al T.$, where $\al T._\c$ is an admissible 
DR-subcategory of $\al T.$, the latter category being
realized as endomorphisms of $\al A.$. A crucial new 
entity that appears when the center $\al Z.$ of $\al A.$ is
nontrivial is the chain group $\ot C.(\al G.)$, which is 
an abelian group constructed from a suitable equivalence
relation in $\wh{\al G.}$ (the dual object of $\al G.$)
and which is isomorphic to the 
character group of the center of 
$\al G.$. Our results suggest the following considerations:

\begin{itemize}
\item
As far as the symmetry $\epsilon$ is concerned, the special case 
studied in Section~\ref{TrivialChainHom}
was also addressed in the context of
vector bundles and crossed products by endomorphisms
(see e.g.~\cite[Eq.~(3.7) and Section~4]{pVasselli03b}). 
In the mentioned reference,
the existence of a symmetry is guaranteed by the 
fact that the left and right $\al Z.$-actions on 
$(\iota,\al E.)$ coincide, where $\al E.$ is a vector bundle
over the compact Hausdorff space $\mr spec.\al Z.$ and 
$(\iota,\al E.)$ denotes the $\al Z.$-bimodule
vector bundle morphism from $\iota:=\mr spec.\al Z.\times\C$
into $\al E.$. (Vasselli studies also bundles in \cite{pVasselli03b},
where left and right actions do not coincide.)

However, the situation analyzed in the present paper
can not be fully compared with 
the case studied in \cite{pVasselli03b}. 
In the latter paper
much more general groups are considered
and, in fact, many of 
them are not even locally compact. For this reason  
no decomposition theory in terms of irreducible objects is mentioned
in that context. It is therefore not clear how the notion of a 
nontrivial chain group should be extended to the general framework of 
vector bundles. Recall that the notion of chain group was 
suggested by the decomposition theory of canonical endomorphisms
and their restriction to $\al Z.$ (cf.~Theorem~\ref{EndoChain}~(ii)).
The nontriviality of the chain group homomorphism
Eq.~(\ref{HomoCA}) gives an obstruction to the existence of a 
symmetry associated with the larger category $\al T.$
(see Proposition~\ref{ABasesZ} and Remark~\ref{NotSym}~(i)).

\item
In lower dimensional quantum field theory models
(see e.g.~\cite{Buchholz88,Mack90,Fredenhagen92}
or \cite[Chapter~8]{bEvans98}), a 
nontrivial center appears when one constructs the
so-called universal algebra. In the case of nets of 
C*-algebras indexed by open intervals of $S^1$,
the universal algebra replaces 
the notion of quasi-local algebra (inductive limit).
(Recall that in this case the 
index set is not directed. See
\cite{FredenhagenIn90} or \cite[Chapter~5]{bBaumgaertel92}.)
Although these models do not fit completely within
the frame studied in this paper (there is no
DR-Theorem and a nontrivial monodromy in two dimensions) 
we still hope that some pieces of the analysis
considered here can be also applied in that situation.
E.g.~the generalization of the notion of irreducible 
objects and the analysis of their restriction to
the center $\al Z.$ that in our context led to the 
definition of the chain group.
\end{itemize}

\section{Appendix: Decomposition of a C*-algebra 
         w.r.t.~its center}\label{AppenDecom}

For convenience of the reader we recall the following facts: let
$\al A.$
be a unital and separable C*-algebra,
$\al Z.$
its center and
$\pi$
a faithful representation of
$\al A.$
on a separable Hilbert space
$\ot h.,\,\pi(\al A.)\subset\al L.(\ot h.).$
According to Gelfand's theorem we have
$\al Z.\cong C(\Gamma)$, where $\Gamma:=\hbox{spec}\,\al Z.$
is a compact second countable Hausdorff space. Then
$\pi\rest\al Z.$
defines a distinguished spectral measure
$E_{\pi}(\cdot)$
on the Borel sets
$\{\Delta\}\subset\Gamma$
such that
\begin{equation}\label{SpecPi}
\pi(Z)=\int_{\Gamma}Z(\lambda)E_{\pi}(d\lambda) \,,
\end{equation}
where $Z(\cdot)\in C(\Gamma)$ is the continuous function 
corresponding to $Z\in\al Z.$
(see e.g.~M.A.~Neumark \cite[p.~278]{bNeumark90}). 
Since $\al Z.$ is the center of $\al A.$ we obtain from
(\ref{SpecPi})
\begin{equation}\label{SpecVertausch}
E_{\pi}(\Delta)\pi(A)=\pi(A)E_{\pi}(\Delta),\quad A\in\al
A.,\,\Delta\subset \Gamma\,.
\end{equation}
Let
$\Phi:\ot h.\rightarrow\wh{\ot h.}:=L^{2}(\Gamma,\mu,
\ot f._{\lambda})$ 
be a unitary spectral transformation assigned to
$E_{\pi}$, where $\mu$ is
a corresponding regular Borel measure on
$\Gamma$ and $\ot f._{\lambda}$
are the fibre Hilbert spaces (cf.~\cite[Chapter~14]{bWallach92}). 
(The spectral representation space $\wh{\ot h.}$
(direct integral) is also denoted in the literature as
$\int_{\Gamma}\ot f._{\lambda}\,\mu(d\lambda)$.)
The transformed projections $E(\Delta)$ on $\wh{\ot h.}$ 
act as multiplication by the corresponding characteristic function
$\chi_\Delta(\cdot)$.

Applying the spectral transformation 
we obtain from the equations
(\ref{SpecPi}) and (\ref{SpecVertausch})
the following inclusions 

\begin{equation}\label{ContInclu}
C(\Gamma)\subset\mr ad.\Phi\circ\pi(\al A.)
         \subset L^\infty\left(\Gamma,\mu,\al L.(\ot f._\lambda) \right)\,,
\end{equation}
where $L^\infty\left(\Gamma,\mu,\al L.(\ot f._\lambda) \right)$
denotes the von Neumann algebra on $\wh{\ot h.}$ of all decomposable
operators $B\colon\ \wh{\ot h.}\to \wh{\ot h.}$ given for almost all
$\lambda\in\Gamma$ by 
\[
 (B\,g)(\lambda):=B(\lambda)\,g(\lambda)\;,\quad g\in\wh{\ot h.}\;,
                  B(\lambda)\in\al L.(\ot f._\lambda)\,.
\]
These operators $\lambda\mapsto B(\lambda)$ are called ``admissible'' 
(see e.g.~\cite[Chapter~4]{bBaumgaertel83}) 
and the measurable Borel function
$\lambda\to \|B(\lambda)\|_\lambda$ satisfies 
$\mr ess~sup.\|B(\lambda)\|_\lambda<\infty$, where $\|\cdot\|_\lambda$
denotes the operator norm in $\al L.(\ot f._\lambda)$.
Then $\|B\|=\mr ess~sup.\|B(\lambda)\|_\lambda$ follows.

The equation (\ref{ContInclu}) means that we have an isomorphism
\[
 \al A.\ni A\mapsto A(\cdot) \;,\quad A(\lambda)\in\al L.(\ot f._\lambda)\,,
\]
and that 
\[
\Big((\hbox{ad}\,\Phi\circ\pi(A))\;\wh{x}\Big)(\lambda)
          =\Big(A(\lambda)\Big)\;\wh{x}(\lambda)\;,
          \quad\mr for.\; \wh{\ot h.}\ni\wh{x}:=\Phi(x)\;,\,x\in\ot h.\,.
\]
Moreover $C(\Gamma)$ is the center of the C*-algebra 
$\mr ad.\Phi\circ\pi(\al A.)$ and if $A(\lambda)=0$ for almost all 
$\lambda\in\Gamma$, then $A=0$.

\begin{teo}\label{C*Zdecom}
Let $(\Gamma,\mu)$ be the measure space mentioned above. Then, there
is an exceptional 
Borel set $\Gamma_0\subset\Gamma$, with $\mu(\Gamma_0)=0$,
such that for all $\lambda\in\Gamma\setminus\Gamma_0$, 
\begin{itemize}
\item[(i)] the set
  \[
   \al A._\lambda:=\{A(\lambda)\mid A\in\al A.\}
  \]
 is a well defined C*-subalgebra of $\al L.(\ot f._\lambda)$ and 
 $\pi_\lambda$, given by
  \[
   \al A.\ni A\mapsto\pi_\lambda(A):= A(\lambda)\in\al A._\lambda\,,
  \]
 is a representation of $\al A.$ on $\ot f._\lambda$,
\item[(ii)] $\al A._\lambda$ has trivial center, 
 i.e.~$\al Z.(\al A._\lambda)=\al A.'_\lambda\cap\al A._\lambda=\C\1_\lambda$.
\end{itemize}
\end{teo}
\begin{beweis}
Part~(i) follows from Eq.~(\ref{ContInclu}). For simplicity we omit in
following the explicit use of the representation $\pi$. 
To show (ii) let 
$\tz$ be a separable abelian C*-algebra containing $\al Z.\cong C(\Gamma)$
and strongly closed in $\al A.'\cap\al A.''$, i.e.
\[
 \al Z.\subset\tz\subset\tz''=\al A.'\cap\al A.''\,.
\]
Then by Gelfand's theorem we have $\tz\cong C(\tg)$, where $\tg$
is a compact second countable Hausdorff space. Moreover we have
$\Gamma\cong\tg / \sim$, where $\sim$ denotes the following equivalence
relation: $\wt{\lambda}._1 \sim \wt{\lambda}._2$ if 
$Z(\wt{\lambda}._1)=Z(\wt{\lambda}._2)$ for all $Z\in\al Z.$. In other words
the elements of $\al Z.$ can be identified with functions in $C(\tg)$
that are constant on the corresponding equivalence classes
(i.e.~let $\lambda\in\Gamma$ and denote by $[\lambda]$ the corresponding
equivalence class in $\tg$, so that for any $\wt{\lambda}.\in[\lambda]$
we have $Z(\wt{\lambda}.)=Z(\lambda)$).

According to M.A.~Neumark \cite[p.~278]{bNeumark90} we have also
for $\tz$ a distinguished spectral measure $\wt E.(\cdot)$ on the 
Borel sets $\wt\Delta.$ in $\tg$ such that
\[
 \wz=\int_{\tg} \wt Z.(\wt\lambda.)\wt E.(d\wt\lambda.).
\]
The relation of the previous decomposition with 
the one given in (\ref{SpecPi}) is specified by the 
following equation:
for $\Delta$ be a Borel set in $\Gamma$ we have
\[
 E(\Delta)=\int_{\wt\Delta.} \wt E.(d\wt\lambda.)\;,\quad\mr where.\quad
          \wt\Delta.=\mathop{\cup}_{\lambda\in\Delta} [\lambda]
          \quad\mr is~a~Borel~set~in.\;\tg\;.
\]

Now the central decomposition of the von Neumann
algebra $\al A.''$ is done over the 
space $\tg$ with regular Borel measure $\wt\mu.$
and the fibre von Neumann algebras $\al A.''(\wt\lambda.)$
are factors for all $\wt\lambda.\in\tg\setminus\wt{\Gamma_0}.$,
where $\wt{\Gamma_0}.=\mathop{\cup}_{\lambda\in\Gamma_0} [\lambda]$
and $\wt\mu.(\wt{\Gamma_0}.)=\mu(\Gamma_0)=0$. 
In the decomposition of $\al A.\subset\al A.''$ over $\wt\Gamma.$
we have that the algebras $\al A._{\wt\lambda.}$ coincide on 
the representatives of the equivalence class $[\lambda]$ with
the algebra $\al A._\lambda$, $\lambda\in(\Gamma\setminus\Gamma_0)$. 
Moreover the functions $\wt\lambda.\to A(\wt\lambda.)\in\al A._\lambda$
have constant values, $A(\wt\lambda.)=A(\lambda)$, for all
$\wt\lambda.\in [\lambda]$.

Finally, for $\wz\in\tz$ we have that $\wz(\wt\lambda.)$, 
$\wt\lambda.\in\tg\setminus\wt{\Gamma_0}.$, are 
scalar functions. If, in particular, $Z\in\al Z.$, then
$Z(\wt\lambda.)$ has a constant value $\zeta_\mt Z.$ for all
$\wt\lambda.\in [\lambda]$. Therefore we have
\[
  Z(\lambda)=\zeta_\mt Z.\1_\lambda
\]
and the proof is concluded.
\end{beweis}

\begin{rem}
We mention here the special case where the spectral measure
$E_{\pi}$
has homogeneous multiplicity. Then there is a unique
fiber Hilbert space
$\ot f.$
and
$\wh{\ot h.}=L^{2}(\Gamma,\mu,\ot f.)$.
Moreover, $\al A._\lambda$ is a C*-algebra on $\ot f.$ 
for $\mu$-almost all $\lambda\in\Gamma$ and for
$A(\lambda)\in\al A._\lambda$ we have
\[
\hbox{ess sup}_{\lambda\in\Gamma}\Vert A(\lambda)
\Vert_{\al L.(\ot f.)}<\infty\;.
\]
If we assume that all operator functions
$\Gamma\ni\lambda\rightarrow A(\lambda)\in\al L.(\ot f.)$
are continuous w.r.t.~the operator norm
$\Vert\cdot\Vert_{\al L.(\ot f.)}$,
then also
$\lambda\rightarrow\Vert A(\lambda)\Vert_{\al L.(\ot f.)}$
is continuous and
\[
 \wh{\al A.}:=\mr ad.\Phi\circ\pi(\al A.)
           \subset C(\Gamma,\al L.(\ot f.))
           \subset L^{\infty} (\Gamma,\mu,\al L.(\ot f.)). 
\]
No exceptional set is needed and $\al A._{\lambda}$
is a unital C*-subalgebra of
$\al L.(\ot f.),\,\{\al A._{\lambda},\wh{\al A.}\}$
is a {\it continuous field} of C*-algebras over $\Gamma$ and
$\wh{\al A.}$
is simultaneously the C*-algebra defined by this field (see
Dixmier \cite[p.~218 ff.]{bDixmier77}). If
$\al A._{\lambda}$
is independent of $\lambda$, for example
$\al A._{\lambda}=\al L.(\ot f.)$ 
(this is true if
$\wh{\al A.}=C(\Gamma,\al L.(\ot f.))=C(\Gamma)\otimes\al L.(\ot f.)$
consists of {\it all} continuous operator functions on
$\Gamma$), then the field is trivial (in the sense of Dixmier).
\end{rem}

\subsection{$\al Z.$-module endomorphisms}\label{SubZMod}

To keep notation simple we omit the 
explicit use of the representation $\pi$. 
Recall that a 
unital endomorphism $\rho$ of $\al A.$ is called a
$\al Z.$-module endomorphism if
\begin{equation}\label{RAZ}
\rho(AZ)=\rho(A)Z,\quad A\in\al A.,\,Z\in\al Z..
\end{equation}

The following proposition can be easily verified using the results
in this section.

\begin{pro}\label{FibreRho}
Let $\rho$ be a unital $\al Z.$-module endomorphism of $\al A.$
and let $\Gamma_0$ be the exceptional set of Theorem~\ref{C*Zdecom}.
Then the family of mappings 
$\{\rho_{\lambda}\colon\ \al A._\lambda\to\al A._\lambda
\}_{\lambda\in (\Gamma\setminus\Gamma_0)}$
defined by 
\[
 \al A._\lambda\ni A(\lambda)
         \to\rho_{\lambda}(A(\lambda)):=(\rho(A))(\lambda)\in\al A._\lambda
\]
is a family of unital endomorphism of $\al A._\lambda$.
\end{pro}

\begin{rem}
Note that the family of endomorphisms
$\{\rho_{\lambda}\}_{\lambda\in (\Gamma\setminus\Gamma_0)}$ 
introduced in the previous proposition satisfies
\[
\Vert \rho_{\lambda}(A(\lambda))\Vert_{\lambda}\leq
\Vert A(\lambda)\Vert_{\lambda}
\quad\mr and.\quad
\hbox{ess sup}_{\lambda\in (\Gamma\setminus\Gamma_0)}
  \Vert\rho_\lambda(A(\lambda))\Vert_{\lambda}<\infty\,,
  \quad A(\lambda)\in\al A.(\lambda) \,.
\]

If $\{\sigma_{\lambda}\}_{\lambda\in (\Gamma\setminus\Gamma_0)}$ 
is {\em any} family of unital endomorphism of 
$\{\al A._\lambda\}_{\lambda\in (\Gamma\setminus\Gamma_0)}$,
then it also satisfies the following boundedness condition:
\[
\hbox{ess sup}_{\lambda\in (\Gamma\setminus\Gamma_0)}
   \Vert\sigma_{\lambda}(A(\lambda))\Vert_{\lambda}
\leq
\hbox{ess sup}_{\lambda\in (\Gamma\setminus\Gamma_0)}
   \Vert A(\lambda)\Vert_{\lambda}<\infty\;,
   \quad A\in\al A.\;.
\]
However, the family 
$\{\sigma_{\lambda}\}_{\lambda\in (\Gamma\setminus\Gamma_0)}$
{\em does not} necessarily define a ``global'' endomorphism $\sigma$ of 
$\al A.$. But if this is the case, then $\sigma$ is also a 
$\al Z.$-module endomorphism, because
\[
(\sigma(AZ))(\lambda)=\sigma_{\lambda}((AZ)(\lambda))=
\sigma_{\lambda}(A(\lambda)Z(\lambda))=Z(\lambda)\sigma_{\lambda}
(A(\lambda))=(\sigma(A)Z)(\lambda)\,.
\]
\end{rem}

\paragraph{Acknowledgments}
First of all we want to acknowledge Hendrik Grundling
for a thorough reading of a first version of the 
manuscript and for many proposals
that improved the presentation of Hilbert C*-systems. We would
also like to thank Michael M\"uger and
Christoph Schweigert for some useful remarks
concerning the chain group. Finally, we are grateful to the
DFG-Graduiertenkolleg ``Hierarchie und Symmetrie in mathematischen Modellen''
for supporting two visits of H.B.~to the RWTH-Aachen.




\end{document}